\documentclass[12pt,reqno,a4paper,twoside]{amsbook}
\usepackage[frenchb]{babel}
\usepackage{pb-diagram,graphicx,epsfig}
\usepackage{amssymb}

\newtheoremstyle{remarks}{3pt}{3pt}{}{}{\bfseries}{}{ }{}

\newtheorem{theorem}[equation]{\bf T{\footnotesize H\'EOR\`EME}}
\newtheorem{definition}[equation]{\bf D{\footnotesize \'EFINITION}}
\newtheorem{prop}[equation]{\bf P{\footnotesize ROPOSITION}}
\newtheorem{postulat}[equation]{\bf P{\footnotesize OSTULAT}}
\newtheorem{lemma}[equation]{\bf L{\footnotesize EMME}}
\newtheorem{corollary}[equation]{\bf C{\footnotesize OROLLAIRE}}

\theoremstyle{definition}
\newtheorem{remark}[equation]{Remarque}

\theoremstyle{remarks}

\DeclareMathAlphabet{\doba}{U}{msb}{m}{n}
\gdef\mC{\doba{C}}

\gdef\mN{\doba{N}}

\gdef\mR{\doba{R}}

\gdef\mZ{\doba{Z}}

\def\qed{{\leavevmode\unskip\nobreak\hfil\penalty 50\hskip 1em%
  \hbox{}\nobreak\hfil\lower 1pt\hbox{$\Box$\kern-.5pt}\parfillskip 0pt
  \finalhyphendemerits 0\par\bigbreak}}
\def\qedmath#1{\setbox0\hbox{$\displaystyle #1$}\templaenge=\textwidth\advance\templaenge by -\wd0%
\setbox1\hbox{$\Box$}\advance\templaenge by -2\wd1%
$$#1\hbox to0pt{\kern.5\templaenge$\Box$\kern-.5pt\hss}$$\par\bigbreak}

\let\<\langle     
\let\>\rangle

\def\al{{\alpha}}
\def\be{{\beta}}

\def\om{{\omega}}
\def\Om{{\Omega}}
\def\la{{\lambda}}

\def\si{{\sigma}}
\def\Si{{\Sigma}}

\def\ep{{\varepsilon}}

\def\th{{\vartheta}}
\let\theta\vartheta
\def\Th{{\Theta}}

\def\phi{{\varphi}}

\def\ol{\overrightarrow}

\def\nummerarray#1#2{\par\noindent\setbox0\hbox{\rm (#1)}\setbox1\hbox{$#2$}\unhcopy0%
\dimen0=.5\textwidth \advance\dimen0 by -\wd0 \advance\dimen0 by -.5\wd1 \kern\dimen0 \unhcopy1}

\def\eref#1{{\rm (\ref{#1})}} 

\def\res#1#2{{#1}\lower .11ex\hbox{$|$}\lower .644ex\hbox{$\scriptstyle #2$}}

\begin{document}

\vfill

\title[]{\large{M\'ecanique
  Quantique}}

\vfill

\author{Emmanuel Humbert\\
{\large Laboratoire de Math\'ematiques et Physique Th\'eorique} \\
{\large Universit\'e F. Rabelais de Tours}\\
{\large Parc de Grandmont} \\
{\large 37000 Tours} \\
{\large humbert@lmpt.univ-tours.fr}\\
{\em Travail soutenu en partie par 
l'ANR-10-BLAN 0105}}
\maketitle

\noindent Ce texte est dans le m\^eme esprit que celui que j'avais \'ecrit sur la relativit\'e g\'en\'erale 
\cite{humbert:10}. 
Il a pour but de pr\'esenter la m\'ecanique en utilisant 
un langage de math\'ematicien tout en ne perdant pas de vue l'aspect physique du probl\`eme.
 La m\'ecanique quantique \'etant par
essence une th\'eorie math\'ematique, on peut \`a juste titre se dire que 
cette d\'emarche n'est pas originale. 
L'originalit\'e r\'eside plut\^ot dans le point de vue que je choisis d'adopter: 
j'essaye de pr\'esenter les intuitions qui ont conduit au formalisme de la th\'eorie 
en allant le plus rapidement possible au but. On ne trouvera donc ici ni exemples, 
ni exercices et seulement tr\`es peu d'applications, bien que la plupart du texte ne d\'epasse pas le niveau de licence.
Je me base principalement sur l'excellent livre de Jean-Louis Basdevant et Jean
Dalibard \cite{basdevant.dalibard:09}. J'utilise \'egalement l'incontournable \cite{messiah:95}. 
Mon travail a surtout consist\'e \`a en extraire les points qui me semblaient les plus importants pour comprendre 
la construction de la m\'ecanique quantique et en r\'epondant \`a certaines questions que, 
de par ma formation purement  math\'ematique, je me suis pos\'ees. Il va sans dire que les choix que j'ai 
faits sont discutables et certains les trouveront sans doute inadapt\'es \`a leur mani\`ere de penser. Je donne aussi 
tr\`es peu de r\'ef\'erences. 
Pour finir, je conseille \`a tout lecteur qui voudrait aborder la m\'ecanique quantique par un ouvrage de vulgarisation, 
de lire le livre d'A. Mouchet \cite{mouchet:10}, qui, sans rentrer dans les d\'etails techniques, pr\'esente de mani\`ere passionnante
les bases de la th\'eorie. \\

\tableofcontents


\chapter{La naissance de la m\'ecanique quantique} \label{naissance}
La m\'ecanique quantique est n\'ee un peu de la m\^eme mani\`ere que la
relativit\'e g\'en\'erale: des observations exp\'erimentales nous ont oblig\'es \`a revoir
compl\`etement notre mani\`ere de penser. Alors que  
la relativit\'e g\'en\'erale a boulevers\'e la vision que nous avions de l'univers \`a grande \'echelle, 
la m\'ecanique quantique remet en cause toute notre intuition concernant la physique des particules. 
En particulier, jusqu'en 1900, mati\`ere et ondes
\'etaient compl\`etement dissoci\'ees. La mati\`ere \'etait compos\'ee de
particules, c'est-\`a-dire de briques \'el\'ementaires, vues alors comme de
petites ``billes''. Les ondes \'etaient d\'ej\`a plus difficiles \`a concevoir,
plus abstraites  
car justement immat\'erielles, bien qu'\'etant omnipr\'esentes au quotidien:
lumi\`ere, son, vagues \`a la surface de la mer,... On pourrait les d\'efinir comme des propagations
de perturbations du milieu ambiant. 
La barri\`ere entre ces deux notions va tomber avec la m\'ecanique
quantique mais il faudra du temps pour accepter de tels changements de
conception du monde qui nous entoure. C'est en 1900 que tout a commenc\'e
lorsque Planck montre que des oscillateurs m\'ecaniques
charg\'es ne peuvent \'emettre ou absorber que des quantit\'es discr\`etes
d'\'energie lumineuse. Le travail de Planck n'est en fait qu'empirique:
pour \^etre pr\'ecis, il montre qu'en faisant cette hypoth\`ese
\'etonnante, on peut \'ecrire une formule simple qui mod\'elise
parfaitement le spectre des corps noirs. Les physiciens de l'\'epoque ne
voyaient l\`a qu'une astuce math\'ematique permettant de trouver des
r\'esultats correspondant \`a l'observation. En 1905, Einstein fut le
premier \`a interpr\'eter physiquement ces id\'ees: la lumi\`ere \'etait
elle-m\^eme compos\'ee de particules, ou tout au moins de ce qu'il appelait
{\em quanta d'action} et qui furent rebaptis\'es {\em photons} en 1926 par
Lewis, ce qui \'etait en compl\`ete contradiction avec le caract\`ere
ondulatoire de la lumi\`ere. Cette interpr\'etation physique des
r\'esultats de Planck fut l'objet de nombreuses
contreverses pendant plusieurs ann\'ees: qu'une onde comme la lumi\`ere
puisse avoir un comportement corpusculaire \'etait \`a l'\'epoque
inacceptable. En 1912, Bohr propose un mod\`ele de l'atome d'hydrog\`ene
convaincant en postulant l\`a-aussi que l'absorption et l'\'emission de lumi\`ere par
la mati\`ere  se faisait par quantit\'es discr\`etes.  \\

\noindent Ce n'est en fait qu'\`a partir de  1914 que ces id\'ees th\'eoriques vont
trouver leur confirmation exp\'erimentale avec les exp\'eriences de Franck
et Hertz qui corroborent parfaitement les pr\'edictions de Bohr sur la
quantification de l'absorption et \'emission d'\'energie lumineuse par les
syst\`emes atomiques ou mol\'eculaires. Plusieurs exp\'eriences viennent
confirmer cette th\'eorie naissante (la m\'ecanique quantique) mais nous
allons nous borner \`a \'etudier celle des deux fentes de Young (qui est en
fait r\'ealis\'ee par F. Shimizu, K. Shimizu et
H. Takuma telle que d\'ecrite ci-dessous) qui prouve de mani\`ere
\'eclatante la dualit\'e onde/mati\`ere.  \\

\noindent L'exp\'erience est la suivante: on lance des atomes de n\'eon
(l'exp\'erience peut aussi \^etre faite avec des \'electrons, neutrons, mol\'ecules)
perpendiculairement \`a une plaque opaque perc\'ee de deux fentes. Un peu
plus loin, on place un \'ecran parall\`ele \`a la plaque et on observe
l'endroit o\`u les atomes de n\'eons vont frapper l'\'ecran apr\`es \^etre
pass\'es dans l'une des deux fentes.

\begin{figure}[h]
\centerline{
\epsfig{file=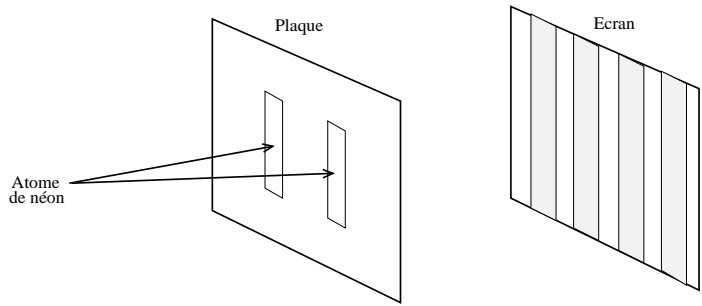, scale=0.5}
}
\caption{Exp\'erience des deux fentes de Young}
\end{figure}

\noindent Voici alors ce qu'on observe: 
\begin{itemize} 
\item Chaque atome laisse une marque sur l'\'ecran en un point
  pr\'ecis. Cela confirme ce que l'on savait d\'ej\`a \`a savoir l'aspect
  corpusculaire de la mati\`ere. 
\item Beaucoup plus \'etonnant: quand on envoie une succession d'atomes
  du m\^eme endroit, avec la m\^eme vitesse et dans la m\^eme direction,
  ils ne frappent pas tous l'\'ecran au m\^eme point. Quand on le fait un
  grand nombre de fois, on s'aper\c{c}oit qu'ils frappent l'\'ecran \`a peu
  pr\`es partout mais beaucoup plus souvent sur des franges (colori\'ees
  sur le dessin ci-dessus). Plus pr\'ecis\'ement, on peut v\'erifier que les atomes viennent
  frapper l'\'ecran de mani\`ere al\'eatoire avec une probabilit\'e dont la
  densit\'e est parfaitement d\'efinie et sinuso\"{\i}dale. Les franges correspondent aux "pics" de cette densit\'e. 
\item Plus \'etonnant encore, si l'on bouche l'une des
  fentes, d'autres franges apparaissent. Ce qui est frappant, c'est qu'il y a
  des r\'egions de l'\'ecran qui sont atteintes r\'eguli\`erement lorsque
  une seule des deux fentes est ouverte et qui ne le sont pas si on ouvre la
  deuxi\`eme fente. 
Ce type de comportement est caract\'eristique de ce qu'on observe pour des
ondes: sans rentrer dans les d\'etails, la somme de deux ondes est encore une onde mais les ``pics''
d'oscillations ne sont plus les m\^emes.  Cette exp\'erience met en \'evidence le comportement ondulatoire de
la mati\`ere.  
\item Le point pr\'ec\'edent sugg\`ere alors une autre exp\'erience: on
ouvre les deux fentes et on  observe les atomes qui passent par la fente 1
et on note leur point d'impact. Tout se passe alors comme si la fente 2
\'etait bouch\'ee. Cela semble contredire l'exp\'erience pr\'ec\'edente. En
fait, pour faire cette deuxi\`eme exp\'erience, on doit s'y prendre
diff\'eremment et cela modifie le r\'esultat. Plus pr\'ecis\'ement, {\bf faire
une mesure influe sur le r\'esultat des exp\'eriences}.
\end{itemize} 

\noindent Nous avons d\'ecrit ci-dessus de mani\`ere tr\`es impr\'ecise \`a
la fois le d\'eroulement de l'exp\'erience et son r\'esultat. L'id\'ee
\'etait seulement de montrer que, comme en relativit\'e g\'en\'erale, il
faut remettre s\'erieusement en question sa mani\`ere de voir les choses
pour b\^atir une th\'eorie coh\'erente.

\chapter{M\'ecanique ondulatoire} \label{mecanique_ondulatoire}
Plusieurs tentatives ont \'et\'e n\'ecessaires avant d'aboutir \`a la formulation 
actuelle de la m\'ecanique quantique. Plus pr\'ecis\'ement, au milieu des ann\'ees 1920, il y avait deux approches 
concurrentes pour mod\'eliser les ph\'enom\`enes quantiques: celle de  Heisenberg,
Born, Jordan et Dirac, appel\'ee {\em m\'ecanique des matrices}, et 
celle de Schr\"odinger, appel\'ee {\em m\'ecanique ondulatoire} mise au point dans une s\'erie de huit articles 
datant 
de 1926 et dont nous pr\'esentons
la d\'emarche dans ce chapitre. Cela permettra de comprendre de mani\`ere relativement intuitive comment
 ces deux th\'eories seront unifi\'ees pour donner la m\'ecanique quantique telle que formalis\'ee actuellement et 
dont les principes seront pr\'esent\'es \`a partir du Chapitre  \ref{mecanique_quantique}. 

\section{Mod\'elisation d'une onde} \label{onde_modele}
En raison de la nature ondulatoire de la mati\`ere d\'ecrite plus haut, nous avons besoin de nous int\'eresser de plus 
pr\`es \`a ce qu'est une onde et \`a la mani\`ere dont on peut la mod\'eliser math\'ematiquement. 
Nous avons expliqu\'e qu'une onde se d\'efinissait comme la propagation d'une certaine pertubation 
"ondulatoire" \`a travers un milieu. Nous avons \'egalement vu que ces ondes pouvaient \^etre de nature 
bien diff\'erentes selon qu'il s'agisse d'un son, d'une vague ou encore d'une onde 
gravitationnelle comme rencontr\'ee en relativit\'e g\'en\'erale. D'un point de vue math\'ematique, 
ces ph\'enom\`enes ont tous en commun de faire intervenir les solutions de l'\'equation des ondes 
\begin{eqnarray} \label{eq_ondes_bas}
\square \phi = 0
\end{eqnarray}
(o\`u $\square:= -\partial_{tt} + \Delta = - \partial_{tt} +\sum_i \partial_{ii}$ est l'op\'erateur D'Alembertien)
ou tout au moins des solutions d'\'equations dont la forme est tr\`es proches et qui peuvent avoir des termes suppl\'ementaires en fonction du milieu observ\'e. C'est de cette constatation que nous allons partir: une onde sera mod\'elis\'ee
 par une fonction $\phi:M \to \mC$, solution d'une certaine \'equation et  o\`u $M$ est l'espace-temps que l'on consid\`ere. Le plus souvent on travaille avec celui de la m\'ecanique classique puisque de nombreux points de la 
physique quantique s'accordent mal avec les th\'eories relativistes.  
La fonction $\phi$ doit contenir toute l'information n\'ecessaire \`a la description de l'onde. \\

\noindent Une solution  \'evidente de l'\'equation \eref{eq_ondes_bas} est la fonction 

\begin{eqnarray} \label{onde_monoc} 
\phi(x,t) = \phi_0 e^{i(k\cdot x - \omega t)}
\end{eqnarray}
o\`u $x$ repr\'esente le vecteur position, $t$ le temps, $k$ le vecteur d'onde (c'est-\`a-dire le vecteur de propagation d'onde),  $\omega$ est la fr\'equence de l'onde et $x \cdot k$ est le produit scalaire. On suppose ici que l'on travaille sur $\mR^n \times \mR$ avec le produit scalaire canonique sur $\mR^n$, espace des vecteurs positions.  $\phi_0$ est une constante que nous laissons de c\^ot\'e pour l'instant. 
La raison pour laquelle le vecteur $k$ est appel\'e vecteur d'onde est la suivante:  on a 
$$\phi(x,t) = \phi(x + w k/||k||,t+1),$$
ce qui signifie que le ph\'enom\`ene observ\'e, quel qu'il soit, a les m\^emes  propri\'et\'es en $(x,t)$ et en $(x + w k/||k||,t+1)$ et a donc un vecteur de propagation \'egal \`a $w k/||k||$.
 Par ailleurs, on observe que, pour que $\phi$ soit solution de \eref{eq_ondes_bas}, il faut et il suffit que $w^2= ||k||^2$. Autrement dit, le vecteur de propagation est bien le vecteur $k$. 

\begin{definition}
Une onde de la forme \eref{onde_monoc} est appel\'ee {\em onde monochromatique}.
\end{definition}

\noindent On remarquera que cette d\'efinition est purement math\'ematique et n'est pour l'instant reli\'ee \`a aucun
ph\'enom\`ene physique.

\section{Mod\'elisation d'une particule en m\'ecanique ondulatoire et \'equation de Schr\"odinger} 
Compte-tenu de l'exp\'erience des fentes de Young d\'ecrite plus haut, l'id\'ee (audacieuse)  va \^etre de mod\'eliser les particules de la 
m\^eme mani\`ere que des ondes, \`a savoir par une fonction $\psi: M  \to \mC$ contenant toute l'information 
sur la particule. Ici, $M$ est l'espace-temps de la m\'ecanique classique et s'\'ecrit $M = M' \times \mR$  o\`u 
$M'$, est le "domaine d'existence" de la particule. Dans ce qui suit, on prendra $M'= \mR^3$.  Pour tenir compte 
de l'aspect "mat\'eriel" de la particule, on va imposer: 

\noindent {\bf Principe 1: } 
{\em la probabilit\'e de trouver la particule \`a l'instant $t$ dans une portion $\Om \subset \mR^3$ est \'egale \`a 

\begin{eqnarray} \label{principe1} 
\int_\Om |\psi(x,t)|^2 dx.
\end{eqnarray} }

\noindent Toute fonction suffisament d\'erivable $\psi:M \to \mC$ n'est pas n\'ecessairement admissible. Par exemple, le principe 1 impose que l'on ait 
\begin{eqnarray} \label{normalisation}
\int_{\mR^3} |\psi(x,t)|^2 dx = 1. 
\end{eqnarray} 

\noindent L'exp\'erience des fentes de Young sugg\`ere \'egalement que les ondes puissent se superposer ce que conduit \`a supposer 

\noindent {\bf Principe de superposition: } {\em Toute combinaison lin\'eaire de fonctions d'onde est \'egalement une fonction d'onde admissible}

\noindent Avant m\^eme de l'expliquer,  \'enon\c{c}ons  le principe fondamental de la m\'ecanique ondulatoire: 

\noindent {\bf Principe 2a: } {\em La fonction d'onde $\psi: M \to \mC$ d'une particule de masse $m$ \'evoluant dans le vide et 
soumise \`a aucune interaction est solution de {\em l'\'equation de Schr\"odinger} 
\begin{eqnarray} \label{eq_schr}
i \hbar \partial_t \psi = - \frac{\hbar^2}{2m} \Delta \psi
\end{eqnarray} 
o\`u $\hbar$ est une constante universelle appel\'ee {\em constante de Planck} et o\`u le laplacien $\Delta$ est 
un laplacien spatial, avec la convention de signe suivante: 
$$\Delta = \partial_{11}+\partial_{22} + \partial_{33}.$$}
 
{\bf On prendra garde au fait que le laplacien ici est celui habituellement utilis\'e en m\'ecanique quantique. Il est \'egal l'oppos\'e de 
celui utilis\'e en g\'eom\'etrie diff\'erentiel (et en particulier dans \cite{humbert:10})}. \\

\noindent La constante de Planck $\hbar$ poss\`ede les dimensions d’une \'energie multipli\'ee par le temps, ou de mani\`ere \'equivalente d'une quantit\'e de mouvement par une longueur.  Voici sa valeur exprim\'ee en Joule.Secondes: 
$$\hbar = 1,054571628×10^{-34}  \hbox{ J.s}$$
avec une incertitude standard de $\pm 0.000 000 053×10^{-34}$  J.s. 

\noindent Dans le cas d'une particule plac\'ee dans un potentiel, on a

\noindent {\bf Principe 2b: } {\em La fonction d'onde $\psi: M \to \mC$ d'une particule plac\'ee dans un potentiel $V(x,t)$  est  solution de: 
\begin{eqnarray} \label{eq_schr2}
i \hbar \partial_t \psi = - \frac{\hbar^2}{2m} \Delta \psi +V \psi.
\end{eqnarray} }

\section{D'o\`u vient l'\'equation de Schr\"odinger ? } \label{douvient}
Dans \cite{humbert:10}, nous avons expliqu\'e 
comment des consid\'erations physiques habiles permettaient d'arriver \`a l'\'equation d'Einstein. 
Les justifications pour postuler les principes 2a et 2b ci-dessus sont beaucoup plus audacieuses. 
D'ailleurs, Schr\"odinger lui-m\^eme n'\'etait pas satisfait de son \'equation. 
La raison principale  de l'\'eriger en tant que principe est
 qu'elle donne des r\'esultats remarquablement conformes par rapport aux 
mesures exp\'erimentales. 
Par contre, les arguments th\'eoriques qui aboutissent \`a 
l'\'equation ont \'et\'e mis au point apr\`es coup et sont  relativement flous. 
Nous n'allons pas les r\'ep\'eter en d\'etail ici. L'une des id\'ees est de 
minimiser une fonctionnelle d'\'energie concernant la fonction d'onde et 
d'enlever les termes non-lin\'eaires de l'\'equation qui en d\'ecoule. 
Cette id\'ee mime le cas classique qui consiste pour d\'ecrire le mouvement des particules \`a 
minimiser une fonctionnelle d'\'energie sur les trajectoires possibles.
 Les raisons pour lesquelles on ignore purement et simplement les termes non-lin\'eaires sont plus 
qu'approximatives. 

\noindent Une deuxi\`eme raison est la suivante: commen\c{c}ons par imaginer ce que pourrait 
\^etre une particule qui se d\'eplace \`a vitesse constante $v$ dans une direction donn\'ee $k$. 
L'id\'ee naturelle serait de prendre pour de telles particules la fonction d'onde monochromatique 
donn\'ee par \eref{onde_monoc}: 
$$  \psi(x,t) = \psi_0 e^{i(k\cdot x - \omega t)}.$$
Il faut d'abord remarquer que cette fonction particuli\`ere n'est pas admissible puisqu'elle ne v\'erifie pas
\eref{normalisation}. 
Rappelons que la quantit\'e de mouvement et l'\'energie totales d'un syst\`eme jouent un r\^ole particulier 
en physique du fait qu'elles sont conserv\'ees avec le temps, du moins dans le cas de syst\`emes isol\'es sans 
potentiel. Beaucoup de th\'eories physiques sont construites en s'appuyant sur ces grandeurs, en particulier en relativit\'e g\'en\'erale (voir \cite{humbert:10}). C'est pourquoi la fonction d'onde ci-dessus sera plut\^ot \'ecrite en termes de quantit\'e de mouvement $p$ et d'\'energie $E$ avec 
\begin{eqnarray} \label{condnec}
E= \hbar \omega \hbox{  et } p^2 = 2mE.
\end{eqnarray} 
Dans le paragraphe \ref{vitesse_p}, nous expliquerons pourquoi le $p$ que nous avons ainsi d\'efini  
s'interpr\`ete bien comme une quantit\'e 
de mouvement. 
Ces conditions sont donn\'ees pour que les ondes d\'efinies ci-dessous par \eref{onde_broglie} et qui sont 
celles que nous consid\`ererons d\'esormais soient solutions de l'\'equation \eref{eq_schr}:  
\begin{eqnarray} \label{onde_broglie}
  \psi(x,t) = \psi_0 e^{\frac{i}{\hbar} (p\cdot x - Et)}.
\end{eqnarray} 
Bien que ne satisfaisant pas \eref{normalisation}, ces fonctions d'onde seront consid\'er\'ees comme "de base" 
pour mod\'eliser au moins localement des particules libres 
de quantit\'e de mouvement $p$ et d'\'energie $E$ et sont appel\'ees 
{\em ondes de de Broglie}.  Avec cette vision, nous pouvons revoir l'exp\'erience des fentes de Young. 
On r\'esout l'\'equation 
\eref{eq_schr} avec les conditions aux limites suivantes: 
\begin{itemize} 
\item $\psi \equiv 0$ en tout point de la plaque opaque, hormis sur les deux trous.
\item Notons $z$ la coordonn\'ee de l'axe perpendiculaire \`a l'\'ecran et \`a la plaque opaque. On veut 
que si $z\to -\infty$ et $t \to -\infty$, l'onde tende vers une onde de de Broglie dans la direction de $z$ (on n\'eglige l'onde r\'efl\'echie sur la plaque opaque).
\item Pour $z \to +\infty$, $\psi \to 0$. 
\end{itemize}
On peut alors montrer que le probl\`eme admet bien une et une seule solution qui conduit \`a des r\'esultats conformes \`a l'exp\'erience.

\section{Compatibilit\'e de l'\'equation de Schr\"odinger avec la d\'efinition de la fonction d'onde.}

Une premi\`ere remarque est la suivante: on peut se demander si la normalisation \eref{normalisation}, 
qui doit \^etre v\'erifi\'ee pour tout $t$, est compatible avec l'\'equation \eref{eq_schr}. 
En fait, non seulement on a compatibilit\'e entre ces relations mais on a m\^eme plus: l'\'equation \eref{eq_schr} 
implique que si la relation \eref{normalisation} est vraie pour un $t \in \mR$, alors elle l'est pour tout $t$. 
En effet, avec \eref{eq_schr} et en utilisant le fait que $\overline{\psi}$ satisfait:
$$i \hbar \partial_t \overline{\psi} =  \frac{\hbar^2}{2m} \Delta \overline{\psi}$$
on voit que 
$$
\begin{aligned} 
\partial_t \int_{\mR^3} |\psi|^2 dx & =  \int_{\mR^3} (\partial_t  \overline{\psi} )\psi dx +  \int_{\mR^3} (\partial_t  \psi) \overline{\psi} dx \\
& =\frac{i \hbar}{2m} \left( \int_{\mR^3} \overline{\psi } \Delta \psi dx -  \int_{\mR^3} ( \Delta\overline{\psi }) \psi dx \right) \\
& = 0.
\end{aligned}
$$
Ici, on a suppos\'e tout de m\^eme qu'une fonction admissible 
permettait de faire l'int\'egration par partie de la derni\`ere ligne, c'est-\`a-dire par exemple que $\psi \in C^2(M) \cap H^{2,2}(\mR^2)$, $H^{2,2}$ d\'esignant l'espace de Sobolev des fonctions $L^2$ dont toutes les d\'eriv\'ees  d'ordre inf\'erieur ou \'egal \`a 2 sont aussi dans $L^2$. Nous ne nous arr\^eterons pas sur ces consid\'erations. Remarquons seulement que de telles hypoth\`eses sont physiquement plausibles: l'action d'une particule \`a l'infini est tr\`es faible, voir inexistante. On peut donc supposer que sa fonction d'onde ainsi que ses d\'eriv\'ees d\'ecroissent suffisamment pour permettre ce genre de calcul.

\section{Paquet d'ondes} 
Malgr\'e ce qui a \'et\'e dit dans le paragraphe \ref{douvient}, on a envie de penser que les ondes de de Broglie, m\^eme si elles ne peuvent mod\'eliser de mani\`ere globale une particule, doivent pouvoir le faire de mani\`ere locale au moins. 
Puisque par ailleurs, une combinaison lin\'eaire de fonctions d'onde admissibles est aussi une fonction d'onde admissible, on a envie de consid\'erer une combinaison linaire de fonctions d'ondes de de Broglie. Malheureusement, cela ne suffit pas pour que la relation \eref{normalisation} soit satisfaite. Pour cela, il faut consid\'erer une "somme infinie non d\'enombrable d'ondes de de Broglie", autrement dit "l'int\'egrale d'une famille d'ondes de de Broglie" que nous pouvons indexer par le vecteur $k$ de propagation de chacune de ces ondes de de Broglie. Plus pr\'ecis\'ement, on dira que 

\begin{definition} 
Un {\em paquet d'onde } est une fonction d'onde de la forme 
$$\psi(x,t) = \int_{\mR^3} \phi(k) e^{\frac{i}{\hbar}(k\cdot x -E(k) t)} dk$$
o\`u l'on doit se rappeler que  $E(k) = \frac{k^2}{2m}$ (voir \eref{condnec}).
\end{definition} 

\noindent En fait, ces consid\'erations sont loin d'\^etre gratuites: le r\'esultat ci-dessous montre que toute 
fonction d'onde admissible est un paquet d'ondes. 
Plus pr\'ecis\'ement, nous montrons le r\'esultat suivant: 
\begin{theorem} \label{paquetdonde}
Soit $\psi: M \to \mC$ une fonction d'onde admissible solution de l'\'equation de Schr\"odinger \eref{eq_schr}. Alors, il existe une fonction $\phi: \mR^3 \to \mC$ (en particulier, qui ne d\'epend pas de $t$) telle que 
$$\psi(x,t) = \frac{1}{(2 \pi \hbar)^{3/2}} \int_{\mR^3}\phi(p) e^{\frac{i}{\hbar}(p\cdot x - E t)} dp$$
o\`u l'on rappelle que $E=|p|^2/(2m)$. De plus, on a la normalisation  
$$\int_{\mR^3} |\phi(p)|^2 dp = \int_{\mR^3} |\psi(x,t)|^2 dx =1.$$ 
\end{theorem}

\noindent {\bf D\'emonstration:}
Soit $f(x,t)$ une fonction admissible solution de \eref{eq_schr}. On d\'efinit pour tous $(p,t)$,  
$$g(p,t) = \frac{1}{(2\pi \hbar)^{3/2}} \int_{\mR^3} f(x,t) e^{-\frac{i}{\hbar} p\cdot x} dx.$$
On calcule alors, en utilisant l'\'equation \eref{eq_schr} puis une int\'egration par parties: 

$$
\begin{aligned} 
\partial_t g(p,t)  & =  \frac{1}{(2\pi \hbar)^{3/2}} \int_{\mR^3} \partial_t f(x,t) e^{-\frac{i}{\hbar} p\cdot x} dx \\
& =  \frac{1}{(2\pi \hbar)^{3/2}} \int_{\mR^3} \frac{i \hbar}{2m} (\Delta_x f(x,t)) e^{-\frac{i}{\hbar} p\cdot x} dx \\
& =   \frac{1}{(2\pi \hbar)^{3/2}} \int_{\mR^3} \frac{i \hbar}{2m}  f(x,t) \Delta_x e^{-\frac{i}{\hbar} p\cdot x} dx \\
& =   \frac{1}{(2\pi \hbar)^{3/2}} \int_{\mR^3} \frac{i \hbar}{2m}  f(x,t)(-\frac{|p|^2}{\hbar^2} ) e^{-\frac{i}{\hbar} p\cdot x} dx \\
& = \frac{- i |p|^2}{2m \hbar} g(x,t).
\end{aligned} 
$$
Ainsi, en posant 
$$\phi(p,t) = g(p,t) e^{\frac{i |p|^2}{2m \hbar} t} = 
g(p,t)e^{\frac{i}{\hbar} Et},$$
on trouve que $\partial_t \phi(p,t)=0$ et donc que $\phi(p,t)=\phi(p)$ ne d\'epend pas de $t$. D'apr\`es la relation \eref{inverse_tf} 
de l'appendice \ref{fourier},  
et puisque que $g$ est la transform\'ee de Fourier inverse  de $f$, on a
$$ \begin{aligned} 
f(x,t)& =  \frac{1}{(2\pi \hbar)^{3/2}} \int_{\mR^3} g(p,t) e^{\frac{i}{\hbar} p\cdot x} dp \\
& =  \frac{1}{(2\pi \hbar)^{3/2}} \int_{\mR^3} \phi(p) e^{-\frac{i}{\hbar}Et}  e^{\frac{i}{\hbar} p\cdot x} dp.
\end{aligned} 
$$
qui est exactement ce que nous cherchions \`a obtenir.  $\square$

\section{Vitesse d'une particule dans le vide} \label{vitesse_p}
\subsection{D\'efinition} 
Revenons maintenant \`a la mod\'elisation d'une particule. Les d\'efinitions donn\'ees plus haut doivent permettre 
de caract\'eriser compl\`etement une particule, en particulier de retrouver sa vitesse. Nous avons d\'efini 
la vitesse de propagation d'une 
onde monochromatique dans le paragraphe \eref{onde_modele} mais comme nous l'avons vu, ces fonctions d'onde 
(celles que nous avons appel\'ees aussi ondes de de Broglie) ne peuvent en aucun cas 
mod\'eliser une particule puisqu'elles ne satisfont pas la relation \eref{normalisation}. Par ailleurs, 
la notion de vitesse qui appara\^{\i}t pour de telles ondes n'est pas g\'en\'eralisable \`a une fonction 
d'onde admissible quelconque. 
D'ailleurs, il n'y a aucune raison de dire qu'il existe un vecteur $k$ tel que $\psi(x+k,t+1)=\psi(x,t)$ comme 
c'est le cas pour les 
ondes de de Broglie. Les d\'efinitions que nous allons utilisons sont donc les suivantes et 
semblent effectivement bien plus naturelles: d'abord, pour parler de vitesse, 
il faut parler de position de la particule. Mais par d\'efinition, une particule n'a justement pas de 
position  pr\'ecise puisqu'elle ne poss\`ede que des probabilit\'es de pr\'esence. 
Ce qui est naturel, c'est de d\'efinir l'esp\'erance de sa position, autrement dit l'esp\'erance de la 
densit\'e de probabilit\'e $\psi$: 
nous parlerons de {\em centre} de la particule que nous d\'efinirons comme 
\begin{eqnarray} \label{centre} 
<x>:= \int_{\mR^3} x |\psi|^2 dx.
\end{eqnarray} 
Ce centre repr\'esente donc la "position la plus probable de la particule". Bien \'evidemment, 
on peut avoir un grand \'ecart type ce qui va donner les "incertitudes de mesure". 
Il est alors naturel de d\'efinir la vitesse de la particule comme la vitesse de $<x>$ qui cette fois a 
bien un sens pr\'ecis. 

\begin{definition}
On consid\`ere une particule de fonction d'onde $\psi$. Sa vitesse sera alors d\'efinie par 
$$<v> := \partial_t <x>= \int_{\mR^3} x \partial_t |\psi|^2 dx.$$
\end{definition}

\noindent M\^eme si la d\'efinition semble naturelle, on peut se poser la question de savoir si la vitesse ainsi d\'efinie
a un sens raisonnable. En particulier, il semblerait justifi\'e de demander qu'une particule \'evoluant dans le vide sans aucune 
contrainte soit anim\'ee d'un mouvement rectiligne uniforme. Le r\'esultat suivant montre que c'est bien le cas: 
\begin{prop} \label{mru}
Le vecteur vitesse d'une particule \'evoluant dans le vide et soumise \`a aucun potentiel est constant. 
\end{prop}
La d\'emonstration sera faite dans le paragraphe suivant. \\

\noindent Nous avons vu dans le paragraphe pr\'ec\'edent que si $\psi$ est la fonction d'onde d'une particule dans le vide, alors
on peut l'\'ecrire sous la forme d'un paquet d'ondes: 
\begin{eqnarray} \label{paquetdonde2} 
 \psi(x,t) = \frac{1}{(2 \pi \hbar)^{3/2}} \int_{\mR^3}\phi(p) e^{\frac{i}{\hbar}(p\cdot x - E t)} dp
\end{eqnarray}
o\`u  $E=|p|^2/(2m)$ et o\`u   
$$\int_{\mR^3} |\phi(p)|^2 dp = \int_{\mR^3} |\psi(x,t)|^2 dx =1.$$ 

\noindent La notation $p$ pour la variable de $\phi$ \'evoque une quantit\'e de mouvement. En fait, 
pour de telle particules, nous d\'efinirons m\^eme: 
\begin{definition}
 La {\em quantit\'e de mouvement moyenne} du paquet d'onde est d\'efinie par 
$$< p >= \int_{\mR^3} p |\phi(p)|^2 dp.$$
\end{definition}

\noindent Cette d\'efinition est justifi\'ee 
par la proposition suivante. 

\begin{prop} \label{qdm}
On a $<p>= mv$ o\`u $v$ est le vecteur vitesse (constant en vertu de la proposition \ref{mru}).  
\end{prop}

\noindent {\bf Cons\'equence importante: } {\em La fonction $|\phi|^2$ s'interpr\`ete donc comme une densit\'e de probabilit\'e pour la 
mesure de la quantit\'e de mouvement. Il existe d'ailleurs d'autres moyens, plus physiques, de le voir. 
En pr\'esence d'un potentiel $V$, on peut montrer que le r\'esultat reste le m\^eme: 
si $\psi$ est la fonction d'onde d'une particule, la densit\'e de probabilit\'e de la quantit\'e de mouvement sera donn\'ee
par $|\phi|^2$ o\`u $\phi$ est la transform\'ee de Fourier de $\psi$. Notons qu'avec un potentiel, la fonction $\phi$ 
peut d\'ependre de $t$. }

\subsection{D\'emonstration de la proposition \ref{mru}}
Puisque la particule \'evolue dans le vide et sans contrainte, sa fonction d'onde $\psi$ est solution de 
l'\'equation \eref{eq_schr}. En prenant le conjugu\'e de cette \'equation, on voit que $\overline{\psi}$ v\'erifie 
$$i \hbar \partial_t \overline{\psi} =  \frac{\hbar^2}{2m} \Delta \overline{\psi}.$$
On en d\'eduit que  
$$
\begin{aligned} 
<v> & =  \int_{\mR^3} x  \partial_t |\psi|^2 dx \\
& =  \int_{\mR^3} x  (\partial_t  \overline{\psi} )\psi dx +  \int_{\mR^3} (\partial_t  \psi) \overline{\psi} dx \\
& =\frac{i \hbar}{2m} \int_{\mR^3} x \left( \overline{\psi } \Delta \psi dx - 
( \Delta\overline{\psi }) \psi  \right) dx \\
& = \frac{i \hbar}{2m}  \int_{\mR^3}  \Delta(x  \overline{\psi }) \psi  -  \Delta (x {\psi }) \overline{\psi} dx. \\.
\end{aligned}
$$
On observe maintenant que, en \'ecrivant $x= x_1 e_1 + x_2 e_2+ x_3 e_3$ et en remarquant que $\Delta x=0$, 
$$
   \Delta (x \psi)  =   \sum_{i=1}^{3} \partial_i \psi + x \Delta \psi.$$  
  La m\^eme formule est valable avec $\overline{\psi}$ si bien que: 
$$  <v> = \frac{i \hbar}{m} \left( \sum_{i=1}^3 (\int_{\mR^3} 
\left( \psi \partial_i \overline{\psi} - \overline{\psi}  \partial_i \psi \right) dx \right) - 
\frac{i \hbar}{2m} \int_{\mR^3} x \left( \overline{\psi } \Delta \psi dx - 
( \Delta\overline{\psi }) \psi  \right) dx.$$
Or la derni\`ere int\'egrale du membre de droite n'est autre que $<v>$ tel que calcul\'e ci-dessus. Ainsi,
\begin{eqnarray} \label{v_t=} 
2<v> = \frac{i \hbar}{m}  \sum_{i=1}^3 (\int_{\mR^3} \left( \psi \partial_i \overline{\psi} - 
\overline{\psi}  \partial_i \psi \right) dx.
 \end{eqnarray}
On \'ecrit maintenant que 
$$\psi \partial_i \overline{\psi} = \partial_i |\psi|^2 -  \overline{\psi}  \partial_i \psi $$
L'int\'egrale du premier terme du membre de droite est nulle. On obtient donc que 
\begin{eqnarray} \label{p=}
2 <v> = -\frac{2 i \hbar}{m}  \sum_{i=1}^3 \int_{\mR^3} \overline{\psi}  \partial_i\psi   dx.
\end{eqnarray}
On d\'erive maintenant cette expression en fonction de $t$ et on utilise de nouveau \eref{eq_schr}: 
$$
\begin{aligned} 
\partial_t <v> & = -\frac{i \hbar}{m}  \int_{\mR^3} (\partial_i \partial_t  \psi )\overline{\psi}+ \partial_i \psi
\partial_t \overline{\psi} dx \\
& =-\left( \frac{i \hbar}{2m} \right)^2 \int_{\mR^3}\left(   \Delta (\partial_i\psi  ) \overline{\psi}  - 
\partial_i \psi  \Delta \overline{\psi}  \right) dx  \\
& = 0.
\end{aligned}
$$
Celta termine la preuve de la proposition \ref{mru}. 

\subsection{D\'emonstration de la proposition \ref{qdm}}
On rappelle que le th\'eor\`eme de Plancherel-Parseval (voir Appendice \ref{fourier}) dit que si $f$, $h$ sont des fonctions $L^2$ 
et si $\hat{f}$, $\hat{h}$ sont leurs transform\'ees de Fourier respectives alors 
$$\int_{\mR^3} f(x) g(x) dx = \int_{\mR^3} \hat{f}(p) \hat{g}(p) dp.$$

\noindent On applique ce r\'esultat avec $f=\overline{\psi}$ et $h= - i \hbar \partial_i \psi$ ($i:=1,2,3$). 
On a donc  $\hat{f} = \overline{\phi}$ et, en vertu de \eref{paquetdonde2},  
$\hat{h}=p_i \phi$ o\`u $p_i$ est la $i$-\`eme coordonn\'ee de $p$. 
On obtient 
$$- i \hbar  \int_{\mR^3}   \overline{\psi}(x) \partial_i \psi(x) dx = <p_i>:= \int_{\mR^3} p_i |\phi(p)|^2 dp. $$
Le r\'esultat se d\'eduit alors de \eref{v_t=}. 

\section{Principe d'incertitude d'Heisenberg} \label{incertitude}
Le principe d'incertitude d'Heisenberg est l'un des points remarquables et contre-intuitifs de la 
m\'ecanique ondulatoire et de la m\'ecanique quantique en g\'en\'eral: plus on conna\^{\i}t 
pr\'ecis\'ement la position d'une particule, 
moins on conna\^{\i}t son vecteur vitesse et vice-versa. 
Avant d'\'enoncer pr\'ecis\'ement ce principe, il convient de remarquer le fait suivant:  
supposons que $g$ est une fonction
telle que $|g|^2$ est la densit\'e de probabilit\'e d'une variable al\'eatoire qui repr\'esente une 
grandeur physique, la position
ou la vitesse par exemple. 
Alors, l'esp\'erance de $|g|^2$, c'est-\`a-dire 
$$<y>:= \int_{\mR^3} y |g(y)|^2 dy$$
repr\'esente la moyenne des mesures que nous ferons de cette grandeur physique. 
Par ailleurs, plus l'\'ecart-type sera grand, plus nous aurons des fluctuations dans ces 
mesures m\^eme si en moyenne 
nous trouvons la valeur $<y>$. 
Autrement dit, 
la quantit\'e 
$$\Delta y^2:= <y^2>  - <y>^2$$
repr\'esente l'incertitude dans la mesure de $y$.

\noindent Prenons maintenant une particule dans le vide repr\'esent\'ee par un paquet d'ondes (voir Th\'eor\`eme 
\ref{paquetdonde})

\begin{eqnarray} \label{paquetdonde3} 
 \psi(x,t) = \frac{1}{(2 \pi \hbar)^{3/2}} \int_{\mR^3}\phi(p) e^{\frac{i}{\hbar}(p\cdot x - E t)} dp
\end{eqnarray}
o\`u  $E=|p|^2/(2m)$ et o\`u   
$$\int_{\mR^3} |\phi(p)|^2 dp = \int_{\mR^3} |\psi(x,t)|^2 dx =1.$$ 

Alors, nous avons vu qu'il \'etait naturel de d\'efinir la position de la particule par 
$$<x>:= \int_{\mR^3} x |\psi(x)|^2 dx.$$
D'apr\`es la discussion ci-dessus, 
la quantit\'e 
$$\Delta x_i^2:= <x_i^2> - <x_i>^2$$
(o\`u $x_i$ est la $i$-\`eme coordonn\'ees de $x$) repr\'esentera l'incertitude  
que nous aurons en mesurant la coordonn\'ee
$x_i$ de la particule. 
De m\^eme, nous avons vu (voir Proposition \ref{qdm}) 
que la fonction $|\phi|^2$ s'interpr\'etait comme une densit\'e de probabilit\'e dans la mesure 
de la quantit\'e de mouvement (et donc dans la vitesse) de la particule. Ainsi, en notant $p_i$ la $i$-\`eme coordonn\'ee
de $p$, la quantit\'e 
$$\Delta p_i^2:= <p_i^2>-<p_i>^2$$
repr\'esentera l'incertitude de mesure de la quantit\'e de mouvement de la particule sur l'axe des $x_i$. \\

\noindent Le principe d'incertitude d'Heisenberg dit alors que 

\begin{prop} \label{incert}
 On a la relation 
$$\Delta x_i \Delta p_i \geq \hbar/2.$$
\end{prop}

\noindent Ce r\'esultat dit que si $\Delta x_i$ est petit, c'est-\`a-dire si on peut mesurer avec pr\'ecision la vitesse de 
la particule, alors $\Delta p_i$ ne peut pas \^etre trop petit, c'est-\`a-dire que la vitesse de la particule (sur l'axe $x_i$)
ne peut pas \^etre mesur\'ee avec autant de pr\'ecision qu'on le souhaite. 

\begin{remark} 
Nous \'etudions ici le cas d'une particule dans le vide mais le r\'esultat resterait le m\^eme avec un contexte physique 
plus g\'en\'eral. 
\end{remark}

\noindent La Proposition \ref{incert} est en fait un r\'esultat classique des 
transform\'ees de Fourier dont nous donnons tout de m\^eme la preuve en raison de l'importance qu'elle a 
pour la th\'eorie.  \\

\noindent {\bf D\'emonstration: } 
De mani\`ere plus g\'en\'erale, si $f$ est une fonction lisse qui d\'ecro\^{\i}t suffisamment \`a l'infini et si $\hat{f}$ est sa 
transform\'ee de Fourier, nous montrons que 
\begin{eqnarray} \label{pih}
 \left(\int_{\mR^3} |f|^2  \right)^2 \leq \frac{4}{\hbar^2}  \int_{\mR^3} (x_i- <x_i>)^2 |f|(x)^2 dx \int_{\mR^3} (p_i- <p_i>)^2 
|\hat{f}|(p)^2 dp 
\end{eqnarray}
Cette in\'egalit\'e implique de mani\`ere triviale la Proposition \ref{incert} en prenant $f= \psi$ 
car en d\'eveloppant
$(x_i- <x_i>)^2$ et  $(p_i- <p_i>)^2$, on voit que 
$$\Delta x_i^2= \int_{\mR^3} (x_i- <x_i>)^2 |f|(x)^2 dx \; \hbox{ et } \; 
\Delta p_i^2 =  \int_{\mR^3} (p_i- <p_i>)^2 
|\hat{f}|(p)^2 dp. $$ 
Soit $g$ une fonction lisse qui tend suffisamment vite vers $0$ \`a l'infini. 
En int\'egrant par partie par rapport \`a la variable $x_i$ et en appliquant l'in\'egalit\'e de H\"older,  on a 
$$\begin{aligned} 
  \left( \int_{\mR^3} g^2 \right)^2 & =  \left(-2  \int_{\mR^3} x_i g(x) \partial_ig (x) dx \right)^2 \\
& \leq 4 \int_{\mR^3}  x_i^2 g(x)^2 dx \int_{\mR^3} (\partial_i g)(x)^2 dx \\
& \leq 4 \int_{\mR^3}  x_i^2 g(x)^2 dx \int_{\mR^3} (\hat{\partial_i g} )(p)^2 dp
 \end{aligned}
$$
o\`u pour la derni\`ere ligne, nous avons utilis\'e le th\'eor\`eme de Plancherel-Parseval (voir Appendice \ref{fourier}) et o\`u
$(\hat{\partial_i g}$ d\'esigne la transform\'ee de Fourier de la fonction $\partial_i g$. 
On remarque maintenant que 

$$\hat{\partial_i g} = \frac{i p_i}{\hbar} \hat{g}$$ ce qui montre que 
$$
  \left( \int_{\mR^3} g^2 \right)^2  \leq \frac{4}{\hbar^2} 
\int_{\mR^3}  x_i^2 g(x)^2 dx \int_{\mR^3} p_i^2 \hat{g}^2 dp.$$

\noindent Soit maintenant $f$ une fonction lisse qui tend suffisamment vite vers $0$ \`a l'infini. 
En appliquant cette in\'egalit\'e avec $g(x)= e^{a x_i} f(x +b x_i)$, $a,b$ \'etant des r\'eels quelconques, et en faisant ensuite
les changements de variables $x_i'= x_i + b$ et $p_i'=p_i+a$, on montre que pour tous r\'eels  $a,b \in \mR$, on a 
 $$
  \left( \int_{\mR^3} f^2 \right)^2  \leq \frac{4}{\hbar^2} 
\int_{\mR^3}  (x_i-b)^2 f(x)^2 dx \int_{\mR^3} (p_i-a)^2 \hat{f}^2 dp$$
ce qui donne l'in\'egalit\'e cherch\'ee. $\square$

\section{Mesures physiques} 
Nous avons expliqu\'e comment donner une d\'efinition coh\'erente de position ou de quantit\'e de mouvement d'une 
particule. Cependant, la fonction d'onde est cens\'ee contenir toutes les propri\'et\'es de la particule \'etudi\'ee. 
Il faut aussi pouvoir mesurer ces grandeurs physique d'une mani\`ere ou d'une autre. Comme nous l'avons 
\'evoqu\'e dans le chapitre \ref{naissance}, ces mesures influent sur le r\'esultat des exp\'eriences. Il faut 
donc r\'eussir \`a inclure ces param\`etres dans la th\'eorie. 
Rappelons ce qu'on a trouv\'e jusqu'\`a pr\'esent: 
\begin{itemize} 
 \item la position (moyenne) d'une particule \`a l'instant $t$ est donn\'ee par 
\begin{eqnarray} \label{forme1} 
<x>:= \int_{\mR^3} \overline{\psi(x,t)} \left(x \psi(x,t) \right) dx;
\end{eqnarray}
\item en vertu de \eref{p=}, 
la quantit\'e de mouvement (moyenne) d'une particule \`a l'instant $t$ et projet\'ee sur l'axe $x_i$ est 
donn\'ee  par 
$$<p_i>=\int_{\mR^3} \overline{\psi(x,t)} \left(\frac{\hbar}{i} \partial_i \psi \right) dx$$ 
et donc la quantit\'e de mouvement (moyenne) vectorielle  $p= mv$ est donn\'ee par 
\begin{eqnarray} \label{forme2} 
 <p>:=\int_{\mR^3} \overline{\psi(x,t)} \left(\frac{\hbar}{i} \nabla \psi \right) dx.
\end{eqnarray}
\end{itemize}

\noindent On remarque que les observables ``position`` et ``quantit\'e de mouvement`` s'\'ecrivent sous la forme 
\begin{eqnarray} \label{formefin} 
\int_{\mR^3} \overline{\psi} \left(\hat{A} \psi \right) dx
\end{eqnarray}
o\`u $\hat{A}$ est un op\'erateur \'evalu\'e en la fonction d'onde. Lorsque l'on \'etudie la position de la
particule, $\hat{A} \psi = x \psi$ et lorsque l'on \'etudie la quantit\'e de mouvement, 
on prend $\hat{A} = \frac{\hbar}{i} \nabla$. On va voir aussi que l'\'energie 
totale peut elle-aussi s'\'ecrire sous une forme similaire. Pour cela,  
supposons que la particule soit soumise \`a un potentiel $V$ qui ne d\'epend pas de $t$ (c'est donc un 
syst\`eme isol\'e)
si bien que sa fonction d'onde 
est solution de 
$$i \hbar \partial_t \psi = - \frac{\hbar^2}{2m} \Delta \psi +V \psi.$$
Notons $\hat{H}$ l'op\'erateur $- \frac{\hbar^2}{2m} \Delta  +V$. On a donc 
\begin{eqnarray} \label{eq_schr2'} 
 i\hbar \partial_t \psi = \hat{H} \psi,  
\end{eqnarray}
et donc aussi 
\begin{eqnarray} \label{eq_schr2''} 
 -i\hbar \partial_t \overline{\psi} = \hat{H} \overline{\psi}.  
\end{eqnarray}
Notons $E$ la quantit\'e 
$$E:= \int_{\mR^3} \overline{\psi} \left(\hat{H} \psi \right) dx.$$
Cette quantit\'e est bien sous la forme \eref{formefin} et nous allons voir de deux mani\`eres diff\'erentes qu'elle 
s'interpr\`ete comme l'\'energie totale du syst\`eme. D'abord,
 on a la relation suivante: 
$$\partial_t E = 0.$$
En effet, en utilisant le fait que $H$ est auto-adjoint puis les relations \eref{eq_schr2'} et \eref{eq_schr2''}
$$\begin{aligned} 
   \partial_t E  &  =  \int_{\mR^3} \partial_t \overline{\psi} \left(\hat{H} \psi \right) dx
+\int_{\mR^3} \left( \hat{H} \overline{\psi} \right)  \partial_t \psi  dx \\
& =-i \hbar \int_{\mR^3} \left(\hat{H} \overline{\psi} \right)  \left(\hat{H} \psi \right) dx
+i \hbar \int_{\mR^3}\left(\hat{H} \overline{\psi} \right)  \left(\hat{H} \psi \right)  dx\\
& = 0. 
  \end{aligned}
 $$
\noindent Or, c'est un postulat physique tr\`es classique que d'interpr\'eter toute int\'egrale premi\`ere scalaire 
(autrement dit une quantit\'e conserv\'ee avec le temps)
d'un syst\` eme isol\'e comme l'\'energie totale du syst\`eme. Avant de montrer que cette interpr\'etation peut se voir 
d'une autre mani\`ere (voir \eref{energie2} ci-dessous), il convient de pr\'eciser un peu les choses. \\

\noindent Ce qui pr\'ec\`ede conduit \`a proposer le \\

\noindent {\bf Principe 3:} {\em \`A chaque grandeur physique $a$, on associe un op\'erateur $\hat{A}$ hermitien qui agit sur
l'espace des fonctions d'onde et tel que la grandeur physique observ\'ee soit en moyenne donn\'ee par 
$$<a>:= \int_{\mR^3} \overline{\psi} \left(\hat{A} \psi \right) dx.$$}

\begin{remark}
On utilisera la plupart du temps les notations ci-dessus: si $a$ est une grandeur physique, $\hat{A}$ d\'esignera son op\'erateur associ\'e et $<a>$ d\'esignera la moyenne de $a$. Par extension, les r\'esultats des  mesures exp\'erimentales  de $a$ seront parfois not\'es eux aussi $<a>$.
\end{remark}
\noindent Il faut alors faire quelques observations importantes: 
soit $a$ une grandeur physique d'op\'erateur associ\'e $\hat{A}$. 
Imaginons une exp\'erience qui mesure  la valeur de $a$ \`a l'instant $t$. 
Notons $\alpha$ la valeur effectivement mesur\'ee. On va faire une hypoth\`ese importante: \`a l'instant $t$, 
on n'a pas d'incertitude sur la valeur de $a$ puisque, justement, on vient d'en faire  la mesure et 
puisqu'on en conna\^it le r\'esultat. 
Le fait de faire cette mesure a donc chang\'e la situation: on sait qu'\`a l'instant $t$,  $<a>= \alpha$ et 
$\Delta a^2= 0$, bien que pour l'instant, nous n'ayons pas donn\'e de moyen de calculer cette valeur. 
Pour tenir compte de ce principe, 
on va faire l'hypoth\`ese (peut-\^etre assez audacieuse) qu'\`a l'instant $t$, on a $\hat{A} \psi = \alpha  \psi$, 
ce qui donne 
bien que:
$$<a> = \int_{\mR^3} \overline{\psi} \hat{A}\psi dx= \alpha \int_{\mR^3} |\psi|^2 dx= \al.$$  
Cette hypoth\`ese signifie que  nous supposons que l'appareil de mesure a donn\'e la valeur exacte de $a$ \`a l'instant $t$. Il  ob\'eit donc  aux lois de la physique
classique  et n'a pas de comportement quantique. Cependant, faire cette hypoth\`ese  a de nombreux avantages: 
\begin{enumerate} 
\item d'abord, elle implique 
que les valeurs mesur\'ees exp\'erimentalement sont n\'ecessairement des valeurs propres de l'op\'erateur $\hat{A}$, qui puisque $\hat{A}$ est auto-adjoint, 
sont r\'eelles. Notons que ces valeurs propres peuvent \^etre vectorielles par exemple dans le cas de la
quantit\'e de mouvement o\`u le vecteur $p$ mesur\'e est une valeur propre qui v\'erifie  $\hat{A}\psi = p\psi$.  Lorsque le spectre de $\hat{A}$ est discret, comme c'est le cas pour l'op\'erateur $\hat{H}$ d\'efini ci-dessus associ\'e \`a l'\'energie totale du syst\`eme, il en d\'ecoule que la quantit\'e mesur\'ee ne peut prendre qu'un ensemble discret, quantifi\'e de valeurs. C'est la raison pour laquelle on parle de m\'ecanique quantique.  Cela nous conforte 
dans l'id\'ee que cette hypoth\`ese est bonne puisque, exp\'erimentalement, ce fait avait d\'ej\`a \'et\'e observ\'e pour l'\'energie (voir le chapitre \ref{naissance}).
\item Soit maintenant une fonction d'onde que nous d\'ecomposons dans l'espace des fonctions propres de l'op\'erateur 
$\hat{A}$. Simplifions la situation en admettant que c'est possible.   
Chaque fonction propre de $\hat{A}$ associ\'ee \`a la valeur propre $\al$  est fonction propre  
de $\hat{A}^2$  de valeur propre $\al^2$. R\'eciproquement, si $\al^2$ est valeur propre de $\hat{A}^2$ alors $\al$ est valeur propre de $\hat{A}$. Il suffit de voir que
si $\psi$ est fonction propre de $\hat{A}^2$ alors $\psi +\frac{1}{\alpha} \hat{A}(\psi)$ est fonction propre de $\hat{A}$.
 Cela signifie que la grandeur physique $a^2$ sera associ\'ee \`a l'op\'erateur 
$\hat{A}^2$. On pourrait faire le m\^eme raisonnement et constater que  l'op\'erateur 
associ\'e \`a $a^n$ est $\hat{A}^n$ pour $n \in \mN$. Par extension, si $f$ est  une fonction analytique, on voit que l'op\'erateur associ\'e \`a $f(a)$   est $f(\hat{A})$. De la m\^eme mani\`ere, si $a$ et $b$ sont des grandeurs 
physiques dont les op\'erateurs associ\'ees sont $\hat{A}$ et $\hat{B}$ alors on pourra se
dire que la grandeur $f(a,b)$ sera 
associ\'ee \`a l'op\'erateur $f(\hat{A},\hat{B})$. En fait, ce n'est pas tout \`a fait vrai parce qu'il est possible
que les op\'erateurs $\hat{A}$ et $\hat{B}$ ne commutent pas. C'est un probl\`eme dont nous ne nous occupons pas maintenant
mais qui aura une v\'eritable importance en m\'ecanique quantique.  

\item En vertu de la remarque pr\'ec\'edente, on peut maintenant calculer facilement l'incertitude de la mesure
de $<a>$ en posant 
$$\Delta a^2:= <a^2> - <a>^2$$
o\`u $<a^2>$ est calcul\'ee avec \eref{formefin} et avec l'op\'erateur $\hat{A}^2$. 
Apr\`es coup, nous en d\'eduisons que si on a mesur\'e $<a>= \alpha$ \`a l'instant $t$, l'hypoth\`ese que l'on a faite
en disant que $\psi(\cdot,t)$ doit \^etre une fonction propre de $\hat{A}$ avec valeur propre $\alpha$ conduit 
\`a avoir \`a la fois $<a>=\alpha$ et $\Delta a^2=0$, qui est exactement ce qu'on voulait. 
\item Encore une fois, consid\'erons un op\'erateur $\hat{A}$ associ\'e \`a une grandeur physi--que $a$ et \'ecrivons 
(si c'est possible: voir le th\'eor\`eme \ref{decomp} de l'appendice \ref{hilbert} ainsi que les corollaires qui en d\'ecoulent)
$\psi= \sum_{\al} a_\al \psi_\al$ o\`u $(a_\al)_\al$ est la famille des valeurs propres de $\hat{A}$ et o\`u $\psi_\al$ est la projection orthogonale dans $L^2$ sur le sous-espace propre associ\'e \`a $a_\al$. Les sous-espaces propres 
\'etant orthogonaux deux \`a deux, on a
$$1= \int_{\mR^3} |\psi|^2 dx = \sum_\al a_\al^2 \int_{\mR^3} | \psi_\al|^2 dx.$$
Cela conduit \`a interpr\'eter $\int_{\mR^3} |\psi_\al|^2 dx$  comme la probabilit\'e de trouver $a_\al$ en faisant une mesure de la grandeur physique $a$.
\end{enumerate}

En r\'esum\'e, nous \'enon\c{c}ons \\
 
\noindent  {\bf Principes 3a:} 
{\em
\begin{itemize} 
 \item Si l'on fait une mesure d'une grandeur physique $a$ associ\'ee \`a un op\'erateur 
$\hat{A}$, nous supposerons que la fonction d'onde $\psi(\cdot,t)$ est une fonction propre de $\hat{A}$. 
\item Les seules mesures exp\'erimentales possibles de $<a>$ sont des valeurs propres de $\hat{A}$. 
\item l'op\'erateur associ\'e \`a la grandeur $f(a)$ ($f$ \'etant une fonction raisonnable) est l'op\'erateur 
$f(\hat{A})$. 
\end{itemize}
}
 
\noindent Terminons le paragraphe par la remarque suivante:  si on consid\`ere une particule 
isol\'ee (i.e.
le potentiel $V$ ne d\'epend pas de $t$), l'\'energie cin\'etique sera donn\'ee
par $1/2 m |v|^2 = \frac{|p|^2}{2m}$ o\`u $p=mv$ est la quantit\'e de mouvement de la particule. L'op\'erateur associ\'e
\`a l'\'energie cin\'etique est donc en vertu de ce qui pr\'ec\`ede \'egal \`a 
$$\hat{E_c}= \frac{1}{2m} |\hat{A}|^2$$
o\`u $\hat{A}= \frac{\hbar}{i}\nabla$ est l'op\'erateur associ\'e \`a $p$. On obtient donc que
$$\begin{aligned} 
   \hat{E_c} & = -\frac{\hbar^2}{2m} \left| (\partial_1,\partial_2,\partial_3) \right|^2 \\
& =  -\frac{\hbar^2}{2m}( \partial_1^2+ \partial_2^2+\partial_3^2)\\
&= -\frac{\hbar^2}{2m} \Delta.
  \end{aligned}
$$
De plus, par d\'efinition, l'\'energie potentielle de la particule est \'egale \`a $V$ et il n'y a aucune incertitude sur cette valeur puisqu'elle ne d\'epend que de $V$. Autrement 
dit, l'op\'erateur associ\'e est $\hat{E_p}\psi:= V \psi$ ce qui donne bien 
$$<V> = V \; \;
\hbox{ et } \; \; \Delta V^2= 0.$$
Puisque l'\'energie totale du syst\`eme est \'egale \`a la somme de l'\'energie cin\'etique et de l'\'energie 
potentielle du syst\`eme, l'op\'erateur associ\'e \`a l'\'energie totale du syst\`eme doit \^etre 
 \begin{eqnarray} \label{energie2}
  \hat{H} = \hat{E_c} + \hat{E_p}= - \frac{\hbar}{2m} \Delta + V
 \end{eqnarray}
qui est exactement ce que nous avions d\'ej\`a postul\'e au d\'ebut du paragraphe en nous basant sur d'autres consid\'erations.

\chapter{Principes de la m\'ecanique quantique} \label{mecanique_quantique}
Dans le chapitre pr\'ec\'edent, nous avons expliqu\'e les intuitions physiques qui ont conduit \`a d\'efinir la m\'ecanique ondulatoire. 
Nous allons maintenant oublier la physique et nous concentrer sur les math\'ematiques. Nous verrons qu'il est naturel de simplifier la th\'eorie. C'est ce qui conduira \`a poser les principes de m\'ecanique quantique qui seront donc plus simples \`a manier d'un point de vue math\'ematique  mais plus \'eloign\'es de la r\'ealit\'e physique dans leur formulation.   C'est pourquoi il est n\'ecessaire pour comprendre les aspects physiques de la m\'ecanique quantique d'avoir auparavant compris la m\'ecanique ondulatoire. 

\section{Principes de base et premi\`eres cons\'equences} 
\subsection{Les principes} 
R\'esumons ce que nous avons fait jusqu'\`a pr\'esent  en m\'ecanique ondulatoire: un 
syst\`eme physique est d\'ecrit par une
 fonction d'onde $\psi$ d\'efinie sur $\mR^3 \times \mR$, \`a valeurs complexes  
et solution de l'\'equation de Schr\"odinger. 
La fonction $\psi(\cdot, t)$ d\'ecrit l'\'etat du syst\`eme \`a l'instant $t$. 
Chaque grandeur physique  $a$ est associ\'ee
\`a  un op\'erateur $\hat{A}$ et la mesure moyenne de $a$ est donn\'ee par 
$$<a> = \int_{\mR^3} \overline{\psi} (\hat{A} \psi) dx.$$
Cette int\'egrale n'est rien d'autre que le produit hermitien $L^2$ $(\psi, \hat{A}\psi)$. 
Par ailleurs, nous avons vu que pour que la grandeur physique $a^2$ soit associ\'ee \`a l'op\'erateur $\hat{A}^2$,
 il fallait d\'ecomposer la fonction d'onde dans la "base" des fonctions propres de $\hat{A}$. 
Pour finir, nous avons vu le r\^ole jou\'e par les transform\'ees de Fourier  qui sont 
en un sens la d\'ecomposition des fonctions sur la "base" non d\'enombrable des $e^{ix}$. 
Ces consid\'erations  nous am\`enent \`a consid\'erer la fonction d'onde $\psi$ 
comme \'etant un vecteur d'un espace de Hilbert hermitien  o\`u nous pourrons alors parler de 
d\'ecomposition dans une base hilbertienne, de produit hermitien et d'op\'erateurs, 
qui sont tous les ingr\'edients que nous avons utilis\'es jusqu'\`a pr\'esent. Plus pr\'ecis\'ement, nous allons poser \\

\noindent {\bf PRINCIPE I : principe de superposition}\\
 {\em \`A chaque syst\`eme physique est associ\'e un espace de Hilbert hermitien $\epsilon_H$ et 
\`a une fonction $\psi : \mR \to \epsilon_H$ (la fonction d'onde) qui v\'erifie pour chaque $t$, $\| \psi \|=1$. 
Le terme ``fonction d'onde'' n'est plus adapt\'e \`a la situation: on dira plut\^ot que $\psi$ 
est le {\bf vecteur d'\'etat} du syst\` eme.\em} \\

\begin{remark}
Dans cette partie, nous utilisons les notations math\'ematiques usuelles qui sont un peu diff\'erentes de celles g\'en\'eralement utilis\'ees par les physiciens. 
Nous reviendrons ult\'erieurement sur ce point.
\end{remark} 

\noindent Pour que la situation soit bien claire, en m\'ecanique ondulatoire, l'espace $\epsilon_H$ \'etait l'espace des fonctions $L^2$ sur $\mR^3$. 
Nous choisissons donc simplement d'une part de ne plus nous limiter \`a 
$\epsilon_H= L^2$ et d'autre part \`a utiliser le langage des  espaces de Hilbert. Ce point de vue a de nombreux avantages que nous exposons d\`es maintenant: 
\begin{enumerate}
\item d'abord, les notations sont beaucoup plus simples. Par exemple, il est beaucoup plus pratique 
d'\'ecrire $(\psi,\hat{A}\psi)$ plut\^ot que  $\int_{\mR^3} \overline{\psi} (\hat{A} \psi) dx$. 
Les calculs gagneront en limpidit\'e. Nous verrons aussi que les grandeurs physiques seront mod\'elis\'ees par 
des op\'erateurs de $\epsilon_H$ sans doute plus facile \`a appr\'e--hender et \`a manier que les 
op\'erateurs diff\'erentiels utilis\'es en m\'ecanique ondulatoire 
(par exemple pour mod\'eliser la quantit\'e de mouvement). 
\item Le fait de travailler sur $\mR^3$ fixait implicitement un syst\`eme de coordonn\'ees. Or dans la plupart des situations physiques, il n'y a pas de 
syst\`eme de coordonn\'ees canonique. Dans le cadre plus g\'en\'eral de la m\'ecanique quantique, nous ne sommes pas contraints \`a de tels choix. 
\item Dans de nombreux cas de figure, nous n'avons pas besoin de manipuler des objets 
aussi compliqu\'es que des fonctions $L^2$. Dans certaines situations, nous pourrons m\^eme travailler 
dans un espace de Hilbert \`a deux dimensions. Cela simplifie 
consid\'erablement les probl\`emes et nous permet de nous concentrer sur 
les grandeurs physiques qui nous int\'eressent.
\item A posteriori, nous verrons que c'est cette formulation plus g\'en\'erale qui permettra de mettre en 
\'evidence le "spin", grandeur intrins\`eque \`a la particule et 
purement quantique, c'est-\`a-dire qui ne poss\`ede aucun analogue en physique classique ni m\^eme en m\'ecanique ondulatoire.
\end{enumerate} 

\noindent Une fois ces bases \'etablies, il s'agit de r\'e\'ecrire les principes de m\'ecanique ondulatoire avec ce nouveau langage.\\

\noindent {\bf PRINCIPE II: mesure des grandeurs physiques} \\
{\em \begin{enumerate}
\item[a)] \`A toute grandeur physique $a$ est associ\'e un op\'erateur auto-adjoint $\hat{A}:\ep_H \to \ep_H$. On dit que $\hat{A}$ est 
l'{\em observable} qui repr\'esente la grandeur $a$. 
\item[b)] {\em (principe de quantification)} On mesure \`a l'instant $t$ une grandeur physique $a$ d'observable $\hat{A}$. 
Quel que soit l'\'etat du syst\`eme, c'est-\`a-dire quel que soit le vecteur d'\'etat $\psi$, 
le r\'esultat de la mesure obtenue sera une valeur propre de $\hat{A}$ (r\'eelle d'apr\`es la proposition \ref{vpreelle} 
de l'appendice \ref{hilbert}).
\item[c)] {\em (principe de d\'ecomposition spectrale) } Notons $(a_\al)_{\al}$ 
les valeurs propres de l'observable $\hat{A}$. Soit $\psi$  le vecteur d'\'etat du syst\`eme. 
La probabilit\'e de trouver $a_\al$ en mesurant la grandeur physique $a_\al$ \`a l'instant $t$ est \'egale \`a 
$\| \psi_\al \|^2$ o\`u $\psi_\al$ est la projection orthogonale de $\psi(t)$ sur le 
sous-espace propre de $\hat{A}$ associ\'e \`a la valeur propre $a_\al$.
\item[d)]  Imm\'ediatement apr\`es avoir mesur\'e la valeur $a_\al$ \`a l'instant $t$, le nouvel \'etat du
 syst\`eme   est $\frac{\psi_\al}{\| \psi_\al \|}$. Plus pr\'ecis\'ement, 
si $P_\al$ est le projecteur sur le sous-espace propre associ\'e \`a $a_\al$, on a 
$\psi(t) = P_\al (\lim_{s\to t^-} P_\al(\psi(s))$ 
\end{enumerate}}

\begin{remark}
Il y a quelques remarques \`a faire sur le principe $IIc)$. De mani\`ere sous-jacente, il faut que probabilit\'e 
totale des mesures possibles 
soit \'egale \`a $1$ et donc qu'on puisse trouver une base hilbertienne d\'enombrable de vecteurs propres. 
C'est possible pour 
un grand nombre d'op\'erateurs (par exemple les op\'erateurs compacts: voir le th\'eor\`eme \ref{decomp} dans 
l'appendice \ref{hilbert} ainsi que
les corollaires qui en d\'ecoulent)  mais pas pour tous. En particuliers, les op\'erateurs qui repr\'esentent une 
grandeur physique quelconque peuvent ne pas \^etre born\'es: c'est le cas en m\'ecanique ondulatoire 
avec les op\'erateurs position, quantit\'e de mouvement et \'energie. 
Malgr\'e tout, nous gardons ce principe tel quel d'abord parce qu'il finira par fournir un mod\`ele tr\`es efficace
qui pr\'edit avec une grande exactitude ce qui se passe, ensuite parce qu'il y aura moyen de construire la th\'eorie 
en \'etant plus rigoureux: il faut pour cela consid\'erer des vecteurs propres "ext\'erieurs \`a  $\ep_H$". 
On trouvera dans l'appendice \ref{hilbert} (paragraphe \ref{discussion}) une discussion sur le sujet. 
Une deuxi\`eme remarque justement li\'ee \`a cette discussion est la suivante: 
le principe $IIc)$ est \'enonc\'e pour un op\'erateur $\hat{A}$ 
dont le spectre est discret mais peut \^etre 
extrapol\'e \`a des  op\'erateurs \`a spectre continu. Dans le cas discret, l'id\'ee est d'\'ecrire 
$\psi= \sum_{\al} a_\al \psi_\al$ o\`u $(a_\al)_\al$ est la famille des valeurs propres de $\hat{A}$ et 
o\`u $\psi_\al$ est la projection orthogonale sur le sous-espace propre associ\'e \`a $a_\al$. 
Les sous-espaces propres \'etant orthogonaux deux \`a deux, on a
$$1= \| \psi\|^2 dx = \sum_\al a_\al^2 \|\psi_\al\ \|^2.$$
C'est pourquoi on pose comme principe d'interpr\'eter $\|\psi_\al\|^2$  comme la probabilit\'e de trouver $a_\al$ 
en faisant une mesure de la grandeur physique $a$. Si maintenant $\hat{A}$ est \`a spectre continu 
 (comme c'est le cas pour l'observable position), on interpr\'etera l\`a aussi $\psi_\al$ comme 
une densit\'e de probabilit\'e. 
Par exemple, si $[x,y]$ est un intervalle de valeurs propres, la probabilit\'e de faire une 
mesure dont la valeur est comprise dans $[x,y]$ sera 
 $$\int_x^y |\psi_z|^2 dz$$
o\`u $\psi_z$ est la projection de $\psi$ sur l'espace propre associ\'e \`a $z$.  
Notons aussi que pour de tels op\'erateurs, les mesures obtenues ne sont pas quantifi\'ees. 
Par exemple, en se pla\c{c}ant dans le cadre de la m\'ecanique ondulatoire, 
l'observable $i$-\` eme coordonn\'ee de la position est donn\'ee par $\hat{A} \psi (x)= x_i \psi(x)$. Autrement dit, 
tout vecteur $x$ est valeur propre de $\hat{A}$ et le vecteur propre associ\'e est la fonction 
$\psi_{x}$ qui vaut $1$ en $x$ et $0$ ailleurs. Le projet\'e de $\psi$ sur l'espace propre associ\'e 
 \`a $x$ est $\psi(x)$. Une autre mani\`ere de le dire est que le projecteur associ\'e est la distribution de Dirac 
en $x$ (lire \`a ce propos la discussion \ref{discussion} dans l'appendice \ref{hilbert}). Si $A \subset\mR^3$, la probabilit\'e de mesurer une position dans $A$ est 
 $\int_{A} |\psi(x)|^2 dx$ et on retrouve bien l'interpr\'etation densit\'e de probabilit\'e de la 
 fonction $|\psi|^2$. Notons que math\'ematiquement, ce que l'on vient de faire  n'est absolument pas rigoureux puisque dans 
 $L^2$, une fonction nulle partout sauf  en un point est \'egale \`a la fonction nulle.  
\end{remark}

\noindent Il faut maintenant donner l'\'equation d'\'evolution du vecteur d'\'etat $\psi$.  En s'inspirant de 
la m\'ecanique
ondulatoire, on posera le principe suivant:

\noindent {\bf PRINCIPE III: \'evolution du syst\`eme}
{\em Soit $\psi(t)$ le vecteur d'\'etat \`a l'instant $t$. Tant que le syst\`eme n'est soumis 
\`a aucune observation,
son \'evolution au cours du temps est r\'egie par l'\'equation de Schr\"odinger 
\begin{eqnarray} \label{schrnew}
 i \hbar \partial_t \psi = \hat{H} \psi
\end{eqnarray}
o\`u $H$ est l'observable \'energie.} 

\noindent Remarquons qu'il n'y a plus de laplacien ici mais que l'on suppose qu'\`a chaque situation physique est 
associ\'ee une ``observable \'energie $\hat{H}$'' qui n'est pas donn\'ee a priori. 
Pour finir avec les principes 
fondamentaux, on va essayer de donner une r\`egle de base en ce qui concerne le choix de l'espace de Hilbert $\ep_H$. 
On a vu en m\'ecanique ondulatoire que le bon cadre pour d\'ecrire le mouvement d'une particule dans $\mR^3$ \'etait de 
travailler dans $L^2(\mR^3)$ (o\`u les fonctions sont \`a valeurs complexes). On aurait pu aussi travailler dans $L^2(\mR)$ ou dans $L^2(\mR^2)$ si on avait \'etudi\'e
le mouvement de la particule le long d'un axe ou dans un plan. Dans un sens, 
chaque coordonn\'ee (on parle plut\^ot de 
{\em degr\'e de libert\'e}) s'\'etudie s\'epar\'ement. De m\^eme si on avait travaill\'e avec 
deux particules \'evoluant dans 
$\mR^3$ (donc avec six degr\'es de libert\'e), on aurait travaill\'e dans $L^2(\mR^6)$. 
Rappelons que $L^2(\mR)$ est s\'eparable,  une  base hilbertienne \'etant 
donn\'ee par les polyn\^omes de Hermite: 
$$H_n(x) = (-1)^n e^{x^2} \frac{d^n}{dx^n} e^{-x^2}.$$
De m\^eme, $L^2(\mR^r)$ est s\'eparable et une base hilbertrienne est donn\'ee par 
$$(H_{k_1} \otimes \cdots \otimes H_{k_r})_{k_1,\cdots,k_r \in \mN}$$
o\`u l'on rappelle (voir le paragraphe \ref{tensoriel} de l'appendice \ref{hilbert}) que pour $(x_1,\cdots, x_r) \in \mC^k$, on a
$$H_{k_1} \otimes \cdots \otimes  H_{k_r}(x_1,\cdots, x_r) := H_{k_1}(x_1)H_{k_2}(x_2)\cdots H_{k_r}(x_r).$$
Autrement dit, $L^2(\mR^r) = L^2(\mR) \otimes \cdots \otimes L^2(\mR)$ ($r$ fois). 
Il est donc naturel de poser le principe suivant: \\

\noindent {\bf PRINCIPE IV: degr\'es de libert\'e}
{\em Chaque degr\'e de libert\'e est d\'ecrit dans un espace de Hilbert. Lorsqu'un syst\`eme physique poss\`ede
$N$ degr\'es de libert\'e, on consid\`erera que l'espace de Hilbert du syst\`eme sera le produit tensoriel des
 $N$ espaces de Hilbert associ\'es \`a chacun des degr\'es de libert\'e.}\\

\subsection{Notations et vocabulaire de Dirac} 
Dans ce texte, nous utilisons le langage et les notations usuelles des espaces de Hilbert. En particulier, tout 
\'el\'ement $\psi$ de l'espace de Hilbert $\ep_H$ peut-\^etre vu ou bien comme un vecteur ou bien une 
forme lin\'eaire continue en posant $\psi^* = (\psi,\cdot)$. Remarquons d'ailleurs que toute forme lin\'eaire continue
de $\ep_H$ peut s'\'ecrire de cette mani\`ere: c'est le th\'eor\`eme de repr\'esentation de Riesz \ref{riesz} pr\'esent\'e
en appendice. En physique, si $\psi \in \ep_H$, on parlera 
\begin{itemize}
 \item de {\em ket} lorsque $\psi$ sera vu comme un vecteur et on a coutume de noter $|\psi>$ 
\`a la place de $\psi$;
\item de {\em bra} lorsque $\psi$ est vu comme une forme lin\'eaire continue. La notation habituelle est alors
$<\psi|$. 
\end{itemize}
En g\'en\'eral, le produit hermitien de $x,y \in \ep_H$ est not\'e $<x|y>$ et les notations pr\'ec\'edentes indiquent
qu'un ket est moralement ``\`a droite'' du produit hermitien tandis qu'un bra est ``\`a gauche''. 

\subsection{Conservation de la norme} 
Nous avons d\'efini un vecteur d'\'etat comme un vecteur unitaire d\'ependant du temps. Comme c'\'etait le cas en m\'ecanique
ondulatoire, l'\'equation de Schr\"odinger permet de montrer la conservation de la norme avec le temps. 
On a en effet, en notant $\psi$ le vecteur d'\'etat: 
$$\begin{aligned}
   \partial_t \| \psi \|^2 & = (\partial_t \psi, \psi) + (\psi , \partial_t\psi)\\
& = (-\frac{i}{\hbar} \hat{H} \psi, \psi) + (\psi ,  -\frac{i}{\hbar} \hat{H}) \\
& =    \frac{i}{\hbar} (\hat{H} \psi, \psi)-  \frac{i}{\hbar} (\psi, \hat{H} \psi).
  \end{aligned}
$$
L'op\'erateur $\hat{H}$ \'etant auto-adjoint, on trouve que  $\partial_t \| \psi \|^2 =0$.

\section{Pr\'edictibilit\'e de l'\'evolution d'un syst\`eme: le chat de Schr\"odinger}
Il est naturel de se poser la question suivante: est-ce qu'il est possible de pr\'edire 
l'\'evolution d'un syst\`eme physique en l'observant \`a un instant donn\'e ? Des principes de base, 
on peut tirer deux conclusions pour r\'epondre \`a
cette question: \\

\noindent {\bf 1)} la premi\`ere consiste \`a remarquer que tant que l'on ne fait pas de mesure, l'\'evolution
du syst\`eme est r\'egie par l'\'equation de Schr\"odinger qui est compl\`etement d\'eterministe. Afin de s'en 
convaincre, on va le montrer lorsqu'on est en pr\'esence d'un syst\`eme isol\'e, c'est-\`a-dire un syst\`eme 
dont l'observable \'energie $\hat{H}$ ne d\'epend pas du temps. 
Dans ce cas, notons $E_\al$ les valeurs propres de $\hat{H}$ et $\psi_\al$ les vecteurs associ\'es. Nous avons 
vu que l'un 
des principes de base de la m\'ecanique quantique consistait \`a supposer que les vecteurs propres des observables 
formaient une base hilbertienne de $\ep_H$. Remarquons que c'est de toute mani\`ere tr\`es justifi\'e en ce qui 
concerne l'observable \'energie puisqu'en m\'ecanique ondulatoire, on avait  
$\hat{H}= -\frac{\hbar^2}{2m} \Delta +V$ et que cet op\'erateur poss\`ede bien une base hilbertienne de vecteurs propres
dans $L^2$. 
On note $\la_\al(t)$ les coefficients du vecteur d'\'etat $\psi(t)$  dans son \'ecriture dans 
la base hilbertienne 
$(\psi_\al)_\al$ autrement dit: $\psi(t)= \sum \la_\al(t) \psi_\al$. 
On a alors d'apr\`es l'\'equation de Schr\"odinger et puisque par d\'efinition $\hat{H} = E_\al \psi_\al$: 
$$\begin{aligned} 
   \partial_t \psi(t) &= -\frac{i}{\hbar } \hat{H} \psi(t) \\
& = -\frac{i}{\hbar }\sum \la_\al(t) \hat{H} \psi_\al \\
&= -\frac{i}{\hbar }\sum \la_\al(t) E_\al \psi_\al.
  \end{aligned}
$$
Puisque par ailleurs
$$\partial_t \psi(t) = \sum \la_\al'(t) \psi_\al,$$
on obtient en identifiant que 
$$\la_\al'(t) =  -\frac{i}{\hbar } \la_\al(t) E_\al$$
c'est-\`a-dire $\la_\al(t) = \la_\al(0) e^{\frac{i}{\hbar}E_\al t}$
et donc 
\begin{eqnarray} \label{decompsi} 
 \psi(t) = \sum  \la_\al(0) e^{\frac{i}{\hbar}E_\al t} \psi_\al.
\end{eqnarray}
\noindent L'\'etat du syst\`eme \`a l'instant $t$ ne d\'epend donc que de l'\'etat du syst\`eme \`a $t=0$. \\

\noindent {\bf 2)} la seconde conclusion que nous pouvons d\'eduire des principes de base est qu'une mesure brise
la pr\'edictibilit\'e ``pass\'ee`` du syst\`eme. Plus pr\'ecis\'ement, en connaissant l'\'etat du syst\`eme juste apr\`es 
une mesure, on ne peut pas conna\^{\i}tre l'\'etat du syst\`eme juste avant cette mesure. En effet,
en notant encore une fois $\la_\al$ les valeurs propres de $\hat{H}$ et 
$\psi_\al$ les vecteurs propres associ\'es (qui peuvent d\'ependre du temps si le syst\`eme n'est pas isol\'e), 
$\psi(t)$ se d\'ecompose \`a l'instant $t$ 
 $$\psi(t)= \sum \la_\al \psi_\al.$$
Si maintenant on effectue une mesure en $t'>t$, les principes de base disent qu'on aura 
$$\psi(t')=  \frac{\psi_\al}{\| \psi_\al \|}$$
\`a l'instant $t'$. En particulier, on n'a absolument aucune indication pour retrouver les coefficients 
$\la_\al$ 
apparaissant \`a l'instant $t$. \\

\noindent Remarquons maintenant que les deux conclusions que nous venons de pr\'esenter semblent conduire \`a
un paradoxe. Consid\'erons un syst\`eme physique $S$ pour lequel on poss\`ede un d\'etecteur $D$ qui 
mesure une certaine quantit\'e physique $a$. Il y deux fa\c{c}ons possibles de voir les choses: 

\begin{enumerate}
 \item on peut d'abord consid\'erer le syst\`eme $S$ comme nous l'avons fait ci-dessus. Le fait de faire une 
mesure brise la pr\'edictibilit\'e du syst\`eme. 
\item On peut aussi, consid\'erer le syst\`eme form\'e de $S$ et de tous les atomes du d\'etecteur $D$. Dans ce cas, 
le d\'etecteur n'influe plus sur le syst\`eme qui est donc r\'egi par l'\'equation de Schr\"odinger. 
\end{enumerate}
Ces deux visions des choses semblent aboutir \`a une contradiction. En fait, dans la deuxi\`eme version, c'est la 
lecture du d\'etecteur qui joue le r\^ole de la mesure. Autrement dit, il faut consid\'erer que l'on 
travaille dans l'espace de Hilbert $\ep_S \otimes \ep_D$ o\`u $\ep_S$ et $\ep_D$ sont les espaces de Hilbert respectivement 
associ\'es \`a $S$ et $D$. Juste avant la mesure, le vecteur d'\'etat du syst\`eme s'\'ecrira
$$\psi(t)= \sum a_{\alpha } S_\al \otimes D_\beta$$
o\`u l'on a d\'ecompos\'e le vecteur $\psi(t)$ dans une base hilbertienne de vecteurs propres de l'observable 
$\hat{A} \otimes \hat{D}$ o\`u  $\hat{A}$ est associ\'e \`a la grandeur physique $a$ et o\`u  $\hat{D}$ est associ\'e 
\`a grandeur physique ''affichage du d\'etecteur``. Le fait de lire ce qui est marqu\'e sur le d\'etecteur donnera 
\`a la fois une mesure de la grandeur physique $a$ ainsi que l'affichage du d\'etecteur. C'est une remarque int\'eressante
parce qu'elle montre que les ph\'enom\`enes macrospiques doivent aussi \^etre soumis aux lois de la m\'ecanique quantique
si l'on veut \'eviter les paradoxes. On r\'esume  ce ph\'enom\`ene en disant qu'avant de lire le d\'etecteur, l'\'etat
du syst\`eme est une superposition de plusieurs \'etats o\`u le d\'etecteur affiche des r\'esultats diff\'erents.  
Dans les livres de vulgarisation, on pr\'esente en g\'en\'eral les choses en disant qu'un chat est enferm\'e
dans une bo\^{\i}te dans laquelle se trouve un syst\`eme l\'etal dont le d\'eclenchement est al\'eatoire. 
Dans la perception classique, on aura tendance \`a dire qu'avant d'ouvrir la bo\^{\i}te, il existe une r\'eponse \`a
la question: ''le chat est-il mort ou vivant?'', c'est-\`a-dire qu'\`a tout instant donn\'e, le chat est soit mort
soit vivant. En m\'ecanique quantique, on supposera que l'\'etat du syst\`eme est une superposition de l'\'etat 
''chat mort`` et l'\'etat ''chat vivant``. Ce n'est qu'au moment o\`u la bo\^{\i}te est ouverte qu'un des 
deux \'etats se d\'ecide.

\subsection{Conservation de l'\'energie}
Par d\'efinition, l'\'energie totale du syst\`eme est
$E(t)=(\psi, \hat{H} \psi)$. Si le syst\`eme est isol\'e, c'est-\`a-dire si $\hat{H}$ ne d\'epend pas du temps, on a:
$$\partial_t E(t) = (\partial_t \psi,\hat{H} \psi)+ (\psi, \hat{H} \partial_t \psi).$$
En utilisant le fait que $\hat{H}$ est auto-adjoint et l'\'equation de Schr\"odinger, on a 
$$\begin{aligned} 
\partial_t E(t) & =(-\frac{i}{\hbar} \hat{H} \psi,\hat{H} \psi) + (\hat{H} \psi, -\frac{i}{\hbar} \hat{H} \psi,\hat{H} \psi)\\
& =\frac{i}{\hbar}(\hat{H} \psi,\hat{H} \psi) - \frac{i}{\hbar}(\hat{H} \psi,\hat{H} \psi)\\
& = 0.
\end{aligned}$$
On retrouve le principe de conservation de  l'\'energie totale d'un  syst\`eme isol\'e.\\

\section{Commutateurs des observables}
\subsection{Postulats} \label{postulat}
Soient $\hat{A}, \hat{B}$ deux observables. A priori, il n'y a aucune raison pour que ces deux  op\'erateurs commutent. C'est la raison pour laquelle on introduit leur commutateur 
$$[\hat{A},\hat{B}] : = \hat{A} \hat{B} -\hat{B}\hat{A}.$$
Par  ailleurs, il n'y a  aucun  moyen de d\'eduire ces  commutateurs des principes de base d\'ecrits plus haut. Il faut donc les postuler. Pour certains d'entre eux, on peut se baser sur ce qui a \'et\'e fait en m\'ecanique ondulatoire. 
Par exemple, notons $\hat{x}_i$ les observables associ\'ees aux composantes du vecteur position, $\hat{p}_i$ ceux associ\'es aux composantes du vecteur quantit\'e de mouvement. 
On rappelle qu'en m\'ecanique ondulatoire, $\hat{x}_i$ est simplement la multiplication par $x_i$, que 
$\hat{p}_i = \frac{\hbar}{i} \partial_{x_i}$ et que $\hat{H}= -\frac{\hbar^2}{2m} \Delta +V$.  Ainsi, les $\hat{x_i}$ commutent entre eux ainsi que les $\hat{p}_i$ et aussi  $\hat{x}_i$ et $\hat{p}_j$ si $i \not=j$. Par contre, 
$$\hat{x}_i \hat{p}_i  - \hat{p}_i \hat{x}_i = \frac{\hbar}{i} x_i \partial_{x_i} (\cdot) - \frac{\hbar}{i} \partial_{x_i} (x_i \cdot) = -\frac{\hbar}{i} \cdot$$
et aussi 
$$\hat{x}_i  \hat{H} - \hat{H} \hat{x}_i= -\frac{\hbar^2r}{2m} x_i \Delta(\cdot) +    \frac{\hbar^2}{2m} 
\Delta( x_i \cdot) = \frac{\hbar^2}{2m}  \partial_{x_i} (\cdot).$$
En r\'esum\'e, on supposera en m\'ecanique  quantique:\\

\begin{postulat} \label{postulat1}  
Pour $i,j$,  on fait les hypoth\`eses:

$$[\hat{x}_i ,\hat{x}_j ]= [\hat{p}_i, \hat{p}_j] =0,$$
\[ [\hat{x}_i,\hat{p}_j]= \left| \begin{array}{cc}
0 & \hbox{ si  } i \not= j \\
- i \hbar & \hbox{ si }   i=j,  \end{array} \right. \]
et enfin
$$[\hat{x}_i, \hat{H} ] = \frac{\hbar^2}{2m}  \partial_{x_i}.$$ 
\end{postulat}

\noindent Les relations de commutation jouent un r\^ole fondamental dans la th\'eorie: il est 
m\^eme tr\`es fr\'equent qu'on ne d\'efinisse
un op\'erateur que par ce type de relations.  Nous allons donner quelques exemples d'applications. 
\subsection{Principe d'incertitude}
Soit $\hat{A}$ une observable associ\'ee \`a une grandeur physique $a$. Par principe, la moyenne $<a>$  de $a$ sera
$$<a> = ( \psi,\hat{A} \psi).$$
La variance $\Delta a^2= <a^2> - <a>^2$ sera
donc $$(\psi, (\hat{A}^2 -<a>^2) \psi)= (\psi, \hat{A}_0^2\psi)$$
o\`u l'on d\'efinit pour toute observable $\hat{A}$,  $\hat{A}_0 := \hat{A} - <a>$. 
De cette observation et du postulat ci-dessus, on d\'eduit ais\'ement le principe  d'incertitude. 
On utilise pour cela l'in\'egalit\'e de Cauchy-Schwarz 
qui dit que pour tous $u,v \in \ep_H$, on a $|(u,v) |\leq \| u\|  \|v \|.$ 
On remarque aussi que pour  deux observables $\hat{A}, \hat{B}$, on a $[\hat{A}, \hat{B}] = [\hat{A}_0, \hat{B}_0]$.
On a donc, en utilisant le postulat ci-dessus:

$$\begin{aligned} 
\hbar{h} & = \left|i \hbar \| \psi\|^2 \right| \\
& = \left|(\psi, [\hat{x}_i,\hat{p}_i] \psi)\right| \\
& = \left|(\psi, [(\hat{x}_i)_0,(\hat{p}_i)_0] \psi)\right| \\ 
& =\left|((\hat{x}_i)_0, (\hat{p}_i)_0 \psi) + ((\hat{p}_i)_0, (\hat{x}_i)_0 \psi) \right| \\
& \leq 2 \| (\hat{x}_i)_0 \psi\|  \| (\hat{x}_i)_0 \psi\| \\
& \leq 2 ( \psi,  (\hat{x}_i)_0^2 \psi) (\psi,(\hat{p}_i)_0^2 \psi) \\
& \leq 2 {\Delta x_i}^2 {\Delta p_i}^2
\end{aligned}
$$
On retrouve donc le principe d'incertitude de m\'ecanique ondulatoire:
\begin{eqnarray} \label{princ_incertitude}
 {\Delta x_i}^2 {\Delta p_i}^2 \geq \hbar/2.
\end{eqnarray}

\subsection{Th\'eor\`eme d'Ehrenfest: conservation de la quantit\'e de mouvement} 
Soit $\hat{A}$ une observable associ\'ee \`a une grandeur physique $a$ (qui peut d\'ependre du temps). En d\'erivant termes \`a termes et en utilisant l'\'equation 
de Schr\"odinger puis le fait que les op\'erateurs sont auto-adjoints, on voit imm\'ediatement que 
\begin{eqnarray} \label{ehrenfest}
 \partial_t <a>= \partial_t (\psi, \hat{A} \psi) = \frac{1}{i \hbar} (\psi, [\hat{A}, \hat{H}] \psi) + (\psi, 
\partial_t\frac{A} \psi).
\end{eqnarray}
 Cette relation, bien qu'\'evidente, est connue sous le nom de th\'eor\`eme d'Ehrenfest. Supposons 
que l'observable $\hat{A}$
ne d\'epende pas du temps et commute avec $\hat{H}$. Alors la quantit\'e physique $a$ (en fait sa moyenne $<a>$) est conserv\'ee
dans le temps. C'est par exemple le cas lorsque $a=p_i$, projection de la quantit\'e de mouvement sur les axes et 
lorsqu'on est dans le cas d'un syst\`eme isol\'e ($\partial_t H = 0$) et d'une particule libre. 
 On rappelle qu'une particule libre est une particule qui n'est soumise \`a aucun potentiel $V$. L'\'energie
totale est alors \'egale \`a l'\'energie cin\'etique 
$$E= 1/m v^2= \frac{p_1^2+p_2^2+p_3^2}{2m}.$$ 
En m\'ecanique quantique, on choisit  donc de lier $\hat{H}$ et $\hat{p}_i$ par  
$$\hat{H} :=\frac{p_1^2+p_2^2+p_3^2}{2m}.$$
On en d\'eduit du postulat \ref{postulat1} que les  $ \hat{p}_i$ commutent avec $\hat{H}$ ce qui implique la conservation de la 
quantit\'e de mouvement.    


\chapter{Du moment cin\'etique \`a la d\'efinition du  spin}  \label{momeci}
 Nous essayerons de donner une d\'efinition la plus simple possible du moment cin\'etique interne des atomes  dans le cadre de la m\'ecanique quantique. En \'etudiant le spectre
de l'observable ainsi d\'efinie, nous nous apercevrons que certaines valeurs propres (donc des valeurs potentiellement mesurables)  ne peuvent pas provenir du moment cin\'etique orbital des particules constituantes de l'atome. Il y a donc deux possibilit\'es: ou bien ces valeurs propres ne sont la traduction d'aucune r\'ealit\'e physique ou au contraire, elles correspondent \`a ce qui doit \^etre une sorte de moment cin\'etique intrins\`eque \`a chaque particule de l'atome. L'exp\'erience  de Stern et Gerlach montrera que c'est la deuxi\`eme solution qui est la bonne. Cette nouvelle grandeur physique des particules sera appel\'ee {\em spin}.  

\section{Moment cin\'etique d'une particule en physique classique} \label{momentci}
Consid\'erons une particule massive en mouvement dans $\mR^3$ mod\'elis\'ee par un point $M$. On d\'efinit son {\em moment cin\'etique} autour d'un point $O$ fixe par 
$$\ol{L} = \ol{OM} \wedge \ol{p}$$
o\`u $\wedge$ est le produit vectoriel de $\mR^3$ et o\`u $\ol{p}$ est la quantit\'e de mouvement de la particule. Chaque composante de $\ol{L}$ mesure la "quantit\'e de mouvement de rotation" autour de l'axe correspondant. 
Par exemple, supposons que la particule tourne dans le plan $z=0$ o\`u le point $O$ a \'et\'e pris 
comme origine et o\`u $(x,y,z)$ sont les coordonn\'ees sur $\mR^3$, alors les coordonn\'ees  
$\ol{L}_x$ et $\ol{L}_y$ de $\ol{L}$ sont nulles puisque par d\'efinition $\ol{L}$ est orthogonal 
\`a $\ol{OM}$ et \`a $\ol{p}$ et donc au plan $z=0$. Notons aussi que si dans ce plan, $M$ ne tourne pas 
autour de $O$, autrement dit, si $\ol{OM}$ et  $\ol{p}$  sont colin\'eaires, alors $\ol{L} = \ol{0}$. 
Dans le plan $z=0$, on peut utiliser les coordonn\'ees polaires. 
On a donc $\ol{OM}= r (\cos(\Theta),\sin(\Theta),0)$ et $\ol{p}= m  \ol{v} = 
m (r'\cos(\Theta), r'\sin(\Theta),0)  + m r \Theta'(-\sin(\Theta), \cos(\Theta))$. 
On trouve ainsi que 
\[ \ol{L} = \left( 
\begin{array}{c}
0 \\ 
0  \\
 m r \Theta' 
\end{array} \right) = \hbox{masse} \cdot  \left( 
\begin{array}{c} 
0 \\ 
0  \\ 
\hbox{vitesse de rotation}  
\end{array} \right) \]
et c'est pourquoi $\ol{L}$ s'interpr\`ete comme une quantit\'e de mouvement de rotation. Il faut 
remarquer que la d\'efinition d\'epend du point $O$.\\

\noindent De par sa d\'efinition, le moment cin\'etique joue un r\^ole important en m\'ecanique du solide. 
Mais il intervient de mani\`ere encore plus fondamentale lorsqu'on 
cherche \`a \'etudier la structure interne de la mati\`ere. 
En effet, de mani\`ere grossi\`ere, un atome est constitu\'e d\'electrons gravitant autour d'un noyau. 
Leur mouvement de rotation, mesur\'e par leur moment cin\'etique que l'on qualifie dans cette situation de
 {\em mouvement cin\'etique  orbital}, induit un champ magn\'etique qui interagit avec le champ magn\'etique 
environnemental. Le moment cin\'etique orbital est donc une composante intrins\`eque de l'atome qui avec 
sa charge lui conf\`ere ses propri\'et\'es \'electromagn\'etiques.   

\section{D\'efinition du mouvement cin\'etique en m\'ecanique quantique}
Une mani\`ere de construire les observables de grandeurs physiques de  la m\'ecanique quantique est de 
s'appuyer sur leur d\'efinition classique. Par exemple, le moment cin\'etique peut se d\'efinir par la 
formule $\ol{L}= \ol{x} \wedge \ol{p}$. On pourra donc poser 
\[ \left\{ \begin{array}{ccc} 
    \hat{L}_1 & = \hat{x}_2 \hat{p}_3 - \hat{p}_2      \hat{x}_3 &= [\hat{x}_2,\hat{p}_3]; \\
\hat{L}_2& = \hat{x}_3 \hat{p}_1 - \hat{p}_3      \hat{x}_1 & = [\hat{x}_3,\hat{p}_1]; \\
\hat{L}_3 &= \hat{x}_1 \hat{p}_2 - \hat{p}_1     \hat{x}_2 &= [\hat{x}_1,\hat{p}_2]
   \end{array}
\right. \]
ce qui se notera  
\begin{eqnarray} \label{defL}
 \hat{L}= \hat{x} \wedge \hat{p}.
\end{eqnarray}
Cette fa\c{c}on de proc\'eder a l'inconv\'enient d'\^etre rigide puisqu'elle oblige \`a d\'efinir les $\hat{x}_i$ 
et les $\hat{p}_i$ et impose aussi l'espace de Hilbert dans lequel on travaille. 
L'autre fa\c{c}on de proc\'eder, et qui sera le plus souvent f\'econde, sera de d\'efinir 
les observables uniquement par leurs commutateurs. 
Par exemple, on peut d\'efinir les op\'erateurs $\hat{x}_i$ et $\hat{p}_i$ uniquement avec le 
postulat \ref{postulat1} \'enonc\'e au chapitre pr\'ec\'edent. Nous avons vu que ces seules 
informations sont suffisantes pour  retrouver des propri\'et\'es physiques importantes 
comme la conservation de la quantit\'e de mouvement d'une particule libre et le principe d'incertitude. 
On d\'ecide donc de suivre cette d\'emarche.
On d\'eduit imm\'ediatement du postulat \ref{postulat1} que  si l'on d\'efinit $\vec{L}$ avec \eref{defL}, on a 
$$ \hat{L} \wedge \hat{L} = i \hbar \hat{L}.$$
Les coordonn\'ees de cette relation donnent les $[\hat{L}_i,\hat{L}_j]$. 
C'est sur cette d\'efinition que nous alllons nous appuyer pour d\'efinir les moments cin\'etiques en m\'ecanique quantique. Il est usuel de les noter $\hat{J}$. 
\begin{definition}
En m\'ecanique, une observable de moment cin\'etique est une observable  vectorielle $\hat{J}$ v\'erifiant 
\begin{eqnarray} \label{defJ}
 \hat{J} \wedge \hat{J} = i \hbar \hat{J}.
\end{eqnarray} 
\end{definition}

\section{Spectre d'un moment cin\'etique} 
Dans ce paragraphe, nous \'etudions le spectre d'une observable de moment cin\'etique. Nous discuterons de ses cons\'equences physiques dans les paragraphes suivants.  
\begin{prop} \label{valeurpropreJ}
Soit 
$$\hat{J}:= \left( \begin{array}{c}
\hat{J}_x \\ \hat{J}_y \\ \hat{J}_z 
\end{array} \right) $$
 une observable de moment cin\'etique  o\`u l'on a not\'e $(x,y,z)$ les coordonn\'ees de $\mR^3$.
 On d\'efinit  $\hat{j}^2:= \| \hat{J} \|^2:= J_x^2 + J_y^2 +J_z^2$.  Alors, 
 \begin{enumerate} 
\item Les valeurs propres de $\hat{j}^2$ sont de la forme $j(j+1) \hbar^2$ o\`u $j \in \frac{1}{2} \mN$. 
\item Les valeurs propres de $\hat{J}_x$, $\hat{J}_y$ et $\hat{J}_z$ sont de la forme $m \hbar$ o\`u $m \in \frac{1}{2} \mZ$. 
\item Soient $m\hbar$, $m  \in \frac{1}{2} \mZ$ une valeur propre de $\hat{J}_z$ (on pourrait prendre  $\hat{J}_x$ ou $\hat{J}_y$) et $j(j+1) \hbar^2$, $j \in  \in \frac{1}{2} \mZ$ une valeur propre de $\hat{j}^2$.  
Alors 
$$m \in \{ -j, -j +1, \cdots,j-1, j \}.$$
\end{enumerate}
\end{prop} 

\noindent {\bf D\'emonstration} \\
On va travailler avec $\hat{J}_z$, les cas de $\hat{J}_x$ et $\hat{J}_y$ se traitant de mani\`ere analogue.
D'abord, on calcule \`a partir des relations \eref{defJ} que $[\hat{j}^2, \hat{J}_z]=0$. Cela implique, 
d'apr\`es le th\'eor\`eme 
\ref{diagsimul} de l'appendice \ref{hilbert}, qu'il existe une base hilbertienne de vecteurs propres communs \`a $\hat{j}^2$ et $\hat{J}_z$. 
Les valeurs propres de $\hat{j}^2$ \'etant positives (puisque $\hat{j}^2$ est associ\'ee \`a une grandeur 
physique positive), on peut les \'ecrire sous la forme $j(j+1)  \hbar^2$ o\`u $j \in \mR$. 
De m\^eme celles de $\hat{J}_z$ s'\'ecrivent sous la forme $m \hbar$. On notera $\psi_{j,m}$ un vecteur propre associ\'e simultan\'ement aux valeurs propres  $j(j+1) \hbar^2$ de $\hat{j}$ et $ m \hbar$ de $\hat{J}_z$. Le but est de montrer que $j,m \in 1/2 \mZ$. On introduit maintenant les op\'erateurs 
$$\hat{J}_+ := \hat{J}_x + i \hat{J}_y \; \hbox{ et } \; \hat{J}_- := \hat{J}_x - i \hat{J}_y.$$
Ces op\'erateurs ne sont pas auto-adjoints: on a $\hat{J}_{\pm}^* = \hat{J}_{\mp}$. Des relations de 
commutations \eref{defJ}, on tire 
\begin{eqnarray} \label{commuj+-}
[\hat{j}^2, \hat{J}_\pm]= 0  \; \hbox{ et } \; [\hat{J}_z, \hat{J}_\pm]= \pm \hbar \hat{J}_\pm.
\end{eqnarray}
On en d\'eduit que, pour $j,m$ fix\'es tels que $\psi_{j,m}$ existe, on a: 
$$\hat{j}^2 \hat{J}_\pm \psi_{j,m}= \hat{J}_\pm \hat{j}^2 \psi_{j,m} = j(j+1) \hbar^2 \hat{J}_\pm \psi_{j,m}.$$
et 
$$\hat{J}_z \hat{j}_\pm\psi_{j,m}= (\hat{J}_\pm \hat{J}_z \pm   \hbar \hat{J}_\pm ) \psi_{j,m}= (m \pm 1 ) \hbar \hat{J}_{\pm} \psi_{j,m}$$
ce qui fait que $\hat{J}_\pm \psi_{j,m}$, s'il est non nul, est un vecteur propre commun aux op\'erateurs $\hat{j}^2$ et $\hat{J}_z$ associ\'e aux valeurs propres $j(j+1)  \hbar^2$ et $(m \pm 1) \hbar$. Consid\'erons maintenant l'ensemble
$$\Om:= \left\{ \hat{J}_\pm^k \psi_{j,m}| k \in \mN \right\}.$$
Cet ensemble est constitu\'e de vecteurs qui sont tous dans le m\^eme espace propre de $\hat{j}^2$ et qui sont propres pour $\hat{J}_z$ associ\'es \`a des valeurs propres de la forme $\al \hbar$ o\`u $\al \in m+ \mZ$. Physiquement, il est raisonnable de penser que seul un nombre fini de vecteurs de $\Om$ sont non nuls puisque que $\hat{J}_z$ est une observable mesurant une coordonn\'ees de moment cin\'etique tandis que $\hat{j}^2$ est une observable mesurant la norme du m\^eme moment cin\'etique. 
On va montrer rigoureusement ce r\'esultat. On \'ecrit 
$$\| \hat{J}_\pm \psi_{j,m}\|^2 = (\psi_{j,m}, \hat{J}_\pm^* \hat{J}_\pm \psi_{j,m}) = (\psi_{j,m},  \hat{J}_\mp \hat{J}_\pm \psi_{j,m}).$$
En uilisant que 
$\hat{J}_\mp \hat{J}_\pm = \hat{j}^2 - \hat{J}_z^2 \mp \hbar \hat{J}_z$, on obtient 
\begin{eqnarray} \label{vpnot0}
\| \hat{J}_\pm \psi_{j;m}\|^2= (j(j+1) - m(m \pm1)) \hbar^2 \| \psi_{j,m} \|^2
\end{eqnarray}
ce qui prouve que, puisque  $\| \hat{J}_\pm \psi_{j;m}\|^2 \geq 0$,
\begin{eqnarray}
- j \leq m \leq j
\end{eqnarray}
On remarque aussi que si $m \not= \pm j$, $\hat{J}_\pm \psi_{j,m} \not= 0$.
On peut r\'e\'ecrire \eref{vpnot0} en rempla\c{c}ant $\psi_{j,m}$ par  $\hat{J}_\pm^k \psi_{j,m}$ et faire le m\^eme raisonnement. 
Une r\'ecurrence \'evidente dit que, puisque $\Om$ est fini, il existe des entiers $k, k'$ 
tels que $m+k=j$ et $m-k'= -j$ (sinon, $\hat{J}_\pm(\psi_{j,m})$ n'est jamais nul). En soustrayant ces in\'egalit\'es, 
on obtient que $2j = k + k'$ et donc $j \in \frac{1}{2} \mZ$. 
On en d\'eduit aussi que $m\in \frac{1}{2} \mZ$ ce qui d\'emontre la proposition.
$\square$

\section{Mouvement cin\'etique orbital}
 Dans le paragraphe pr\'ec\'edent, nous avons travaill\'e avec un  moment cin\'etique d\'efini par  
\eref{defJ}. Cette d\'efinition est purement math\'ematique: nous avons uniquement utilis\'e les relations de 
commutation d'un moment cin\'etique classique. Il est naturel de se demander quel est le sens physique de cette 
d\'efinition. Pour cela, 
repla\c{c}ons-nous dans le cadre de la m\'ecanique ondulatoire avec les obervateurs $\hat{x}$, $\hat{p}$ et 
on d\'efinit $L$ par \eref{defL}. Cette d\'efinition plus physique doit correspondre \`a la situation qui nous a 
inspir\'es: celle d'une particule, d'un \'electron en particulier, en orbite autour d'un noyau. 
Un tel moment cin\'etique sera appel\'e {\em moment cin\'etique orbital}. Nous montrons que 

\begin{prop} \label{valeurpropreL}
Soit $\hat{L}$ un observable de moment cin\'etique orbital. Notons $\hat{L}_x$, $\hat{L}_y$ et $\hat{L}_z$ les observables "coordonn\'ees" associ\'ees. Leurs valeurs propres sont de la forme $m \hbar$ o\`u $m \in \mZ$.
\end{prop}
\noindent {\bf D\'emonstration}\\
Comme dans le paragraphe pr\'ec\'edent, nous travaillons seulement avec $\hat{L}_z$.  La d\'efinition \eref{defL} ainsi que les observables de position et de quantit\'e de mouvement de la m\'ecanique ondulatoire nous disent que 
$$\hat{L}_z= i \hbar( x \partial_y - y \partial_x).$$
On prendra garde au fait qu'ici, $(x,y,z)$ sont les coordonn\'ees sur $\mR^3$ tandis que dans \eref{defL}, 
$\hat{x}$ repr\'esente l'observable vectoriel position. En utilisant les coordonn\'ees sph\'eriques 
$(r,\phi,\Theta)$ sur $\mR^3$, on voit que 
\begin{eqnarray} \label{Lz} 
\hat{L}_z= -i \hbar \partial_\phi.
\end{eqnarray}
Soit $\psi_m(r,\phi,\Theta)$ - nous n'\'ecrivons pas la d\'ependance en temps - une fonction d'onde propre pour $\hat{L}_z$, c'est-\`a-dire v\'erifiant 
$$\hat{L}_z \psi_m = m \hbar \psi_m$$
o\`u $m \in \mR$.  L'\'ecriture \eref{Lz} permet de dire que $\psi_m$ est de la forme
$$\psi_m(r,\phi,\Theta) = \Psi_m(r,\Theta) e^{i m \phi}.$$
Puisqu'en coordonn\'ees sph\'eriques, $(r,\phi,\Theta) = (r, \phi + 2\pi, \Theta)$, il faut que 
$\psi_m(r,\phi,\Theta) = \psi (r, \phi + 2\pi, \Theta)$ ce qui montre que $m$ doit \^etre entier.  $\square$

\section{Interpr\'etation physique des propositions \ref{valeurpropreJ} et \ref{valeurpropreL}: spin d'une particule}
La proposition \eref{valeurpropreJ} dit que les valeurs propres des coordonn\'ees d'un moment cin\'etique d\'efini par \eref{defJ} sont de la forme $m \hbar$ o\`u $m \in 1/2 \mZ$. Autrement dit, ce sont les seules valeurs que l'on peut mesurer exp\'erimentalement. D'apr\`es la proposition \ref{valeurpropreL}, celles de ces valeurs qui ne sont pas enti\`eres ne peuvent pas caract\'eriser des mouvements cin\'etiques orbitaux. Il y a donc plusieurs solutions:
\begin{enumerate}
\item ou bien ces valeurs  non enti\`eres ne sont pas valeurs propres:
 nous avons en effet montr\'e que toute valeur propre \'etait de la forme  $m \hbar$ o\`u $m \in 1/2 \mZ$ mais nous n'avons pas montr\'e que toutes ces valeurs \'etaient effectivement des valeurs propres;
\item ou bien ces valeurs sont des valeurs propres mais ne caract\'erisent aucune r\'ealit\'e physique: il est en effet raisonnable d'imaginer que la d\'efinition \eref{defJ}, bas\'ee uniquement 
sur des relations de commutation, n'est pas la bonne. N'oublions pas que le passage de 
la m\'ecanique ondulatoire \`a la m\'ecanique quantique s'est fait en oubliant un peu la physique et en math\'ematisant la th\'eorie;
\item ou bien ces valeurs sont des valeurs propres et caract\'erisent une grandeur physique 
qui n'a donc aucun analogue classique. 
\end{enumerate}
C'est l'exp\'erience de Stern et Gerlach qui va trancher la question: c'est la troisi\`eme version qui est la bonne.
 Les valeurs demi-enti\`eres mesurent une quantit\'e physique que l'on appelle {\em spin} et qui est une propri\'et\'e 
intrins\`eque des particules. On l'\'etudiera plus en d\'etail dans le chapitre suivant.
Son interpr\'etation physique est encore discut\'ee \`a l'heure actuelle. 
Certains l'imaginent comme une rotation d'une particule sur elle-m\^eme mais cette vision des choses n'est pas satisfaisante \`a bien des \'egards. 
On peut d'ailleurs se poser la question: est-elle r\'eellement le reflet d'une quantit\'e humainement 
imaginable comme une rotation ? Ce n'est pas du tout \'evident et on pourrait penser qu'il s'agit simplement 
d'une grandeur purement math\'ematique 
qui doit \^etre prise en compte lorsque l'on \'etudie le comportement d'une particule. 

\section{Exp\'erience de Stern et Gerlach} 
La construction est la suivante: on projette des atomes horizontalement dans un champ magn\'etique 
vertical et on mesure
la d\'eviation de la trajectoire sur un \'ecran plac\'e un peu plus loin. On peut montrer que 
l'interaction de l'atome avec le champ magn\'etique provient 
uniquement de la composante $J_z$ du moment
cin\'etique interne $J$.
Nous ne rentrerons pas dans les d\'etails pratiques de l'exp\'erience 
ni m\^eme dans son interpr\'etation physique exacte. 
N\'eanmoins, ce qu'on observe est que les atomes, lorsqu'il s'agit d'atomes
d'argent,  vont frapper l'\'ecran 
\`a deux endroits seulement de l'\'ecran. On peut raisonnablement penser que l'atome n'a, 
en termes de moment cin\'etique
que deux \'etats possibles, c'est-\`a dire que l'op\'erateur $\hat{J}_z$ 
n'a que deux valeurs propres. Cela sugg\`ere que le param\`etre $j$ de la proposition  
\ref{valeurpropreJ} prend la valeur $1/2$ avec probabilit\'e $1$ et que les valeurs propres de $\hat{J}_z$ sont 
$\pm 1/2 \hbar$. On pourra aussi \`a ce propos, lire la remarque \ref{remarquefs}.  Plus g\'en\'eralement, en supposant que le param\`etre $j$ soit fix\'ee, 
on voit que $\hat{J}_z$ 
a $2j+1$ valeurs propres si $j \in \mZ$ et $2j$ valeurs propres si $j \in 1/2 \mZ 
\setminus \mZ$. On observera d'autres exemples d'atomes o\`u apparaissent un nombre pair de taches, impliquant des 
valeurs
demi-enti\`eres pour $j$ et $m$. Autrement dit, les conclusions que nous avons tir\'ees sont bonnes, 
ce qui 
sera confirm\'e par toutes les pr\'edictions exp\'erimentales que nous pouvons en tirer et le moment 
cin\'etique interne de l'atome
d'argent n'est pas uniquement d\^u \`a des moments cin\'etiques orbitaux. 

\begin{remark}
 Il est int\'eressant 
de voir la d\'emarche adopt\'ee par les physiciens face \`a ce type de probl\`eme. Les r\'esultats de l'exp\'erience 
disent que l'atome peut \^etre dans deux \'etats distincts, correspondant chacun \`a une valeur propre de 
l'observable associ\'ee. Il y a donc au moins deux vecteurs propres. On essaye donc de voir ce qui se passe en 
supposant que l'espace de Hilbert associ\'e au degr\'e de libert\'e $J_z$ est exactement de dimension $2$, engendr\'e 
par les deux \'etats propres mis en \'evidence par l'exp\'erience. C'est cette d\'emarche qui conduira 
au formalisme du spin $1/2$ d\'ecrit dans la suite. Les r\'esultats obtenus exp\'erimentalement montrent 
que cette hypoth\`ese est bonne. On est donc dans une situation o\`u l'on travaille
avec des espaces de Hilbert beaucoup plus simples que $L^2$ rencontr\'e en m\'ecanique ondulatoire
 et qui justifie encore une fois la formalisme de la m\'ecanique 
quantique. 
\end{remark}


\chapter{Le spin et son formalisme}
Le chapitre pr\'ec\'edent montrait comment exp\'erimentalement, on pouvait mettre en \'evidence le {\em spin} : 
grandeur purement quantique s'apparentant \`a un moment cin\'etique interne. Dans ce chapitre, nous allons 
montrer comment d'une part mod\'eliser cette grandeur physique et d'autre part, \'etudier les cons\'equences de son existence. Nous \'etudierons ici seulement le {\em spin $1/2$}.   
\section{D\'efinition} 
R\'esumons d'abord les conclusions du chapitre pr\'ec\'edent. Le but \'etait de d\'efinir un moment 
cin\'etique en m\'ecanique quantique. Au lieu de nous appuyer  sur la d\'efinition classique de cette 
grandeur vectorielle, 
nous avons d\'ecid\'e de d\'efinir l'observable associ\'ee en nous basant uniquant sur les 
relations de commutation entre ses coordonn\'ees. Nous d\'efinissions donc un 
observable  $\hat{J}= (\hat{J}_x,\hat{J}_y,\hat{J}_z)$ v\'erifiant les relations \eref{defJ}. On observait alors que 
\begin{enumerate} 
\item les valeurs propres de l'observable $\hat{j}= \hat{J}^2_x +\hat{J}^2_y +\hat{J}_z^2$ \'etaient de la forme 
$j(j+1) \hbar^2$ avec $j \in 1/2 \mZ$; 
\item les valeurs propres de l'observable $\hat{J}_z$ 
\'etaient de la forme 
$m \hbar$ avec $m \in 1/2 \mZ$; 
\item les valeurs non-enti\`eres de $j$ et $m$, dont l'existence physique 
\'etait mise en \'evidence par l'exp\'erience de Stern et Gerlach, ne pouvaient pas correspondre au
 moment cin\'etique "classique" d'une particule, bien qu'induisant un comportant similaire
de l'atome en pr\'esence d'un champ magn\'etique. 
\end{enumerate}

\noindent Nous avons avons aussi montr\'e exp\'erimentalement que dans le cas de l'atome d'argent, 
l'observable $\hat{j}$ n'avait en fait qu'une seule valeur propre correspondant \`a $j=1/2$ et que 
l'observable $\hat{J}_z$ n'avait que deux valeurs propres correspondant \`a $m= \pm 1/2$. 
Nous dirons que cet atome a un {\em spin $1/2$}. 
Exp\'erimentalement, on peut montrer que les particules elles-m\^emes, et pas seulement les atomes,
 poss\`edent ce "moment magn\'etique propre". Plus pr\'ecis\'ement, 
les \'electrons, protons, neutrons, neutrinos, quarks poss\`edent un spin $1/2$. 
On peut aussi montrer que d'autres particules poss\`edent un spin plus
\'elev\'e correspondant \`a d'autres valeurs  plus grandes de $j$ et $m$. On classe 
d'ailleurs les particules selon qu'elles aient ou non un spin entier: 
\begin{definition}
 Les particules de spin entier sont appel\'ees des {\em bosons} et les particules de spin demi-entier sont appel\'ees 
des {\em fermions}.
\end{definition}

\noindent Malgr\'e tout, nous n'avons toujours donn\'e aucune d\'efinition math\'ematique pr\'ecise du spin. 
Nous le faisons maintenant, en nous limitant au spin $1/2$. 

\begin{definition} 
Les particules pour lesquelles on ne peut mesurer que deux \'etats distincts  pour chaque 
coordonn\'ee du moment cin\'etique propre
(par une exp\'eri--ence de type Stern et Gerlach par exemple) seront dites 
\`a {\em spin $1/2$}. Cette grandeur physique $S= (S_x,S_y,S_Z)$ est un degr\'e de libert\'e 
\`a elle toute seule: elle ne peut pas se lire 
\`a travers les observables position et vitesse par exemple (contrairement \`a un moment cin\'etique classique).
Ce fait n'est pas du tout \'evident. On lui attribue donc   
un espace $\ep_{spin}$ de dimension $2$ (une dimension pour chaque valeur propre de $J_z$) 
et un observable vectoriel 
$\hat{S} = (\hat{S}_x, \hat{S}_y,\hat{S}_z)$ v\'erifiant les relations de commutation 
\begin{eqnarray} \label{defS}
\hat{S} \wedge \hat{S} = i \hbar \hat{S}
\end{eqnarray}
On pose comme principe que chacune des observables $\hat{S}_x, \hat{S}_y,\hat{S}_z$ poss\`ede 
exactement les valeurs propres $\pm 1/2 \hbar$. 
\end{definition} 

\begin{remark} \label{remarquefs}

\noindent 

 \begin{enumerate} 
  \item Puisque les mesures des coordonn\'ees $S_x,S_y,S_z$ ne donnent que les valeurs possibles
$\pm 1/2 \hbar$, la mesure de $\| S \| := S_x^2+S_y^2+S_z^2$ donne la valeur 
$\frac{3}{4} \hbar^2 = j(j+1)\hbar^2$ 
o\`u  $j= 1/2$. En particulier, cela implique que $\hat{s}:= \hat{S}_x^2+ \hat{S}_y^2+ \hat{S}_z^2$ est \'egal \`a 
$ \frac{3}{4} \hbar^2 Id$. 
\item On pourrait penser que le spin $S$ poss\`ede trois degr\'es de libert\'e (un pour chaque coordonn\'ee): c'est en fait 
proscrit par le fait que les observables $\hat{S}_x, \hat{S}_y,\hat{S}_z$ ne commutent
pas et donc que ces quantit\'es ne sont pas ind\'ependantes les unes des autres. Cela dit, le choix de mod\'eliser 
le spin $1/2$ par un espace de Hilbert \`a deux dimensions ne s'impose pas. Il se trouve que, puisqu'il y a au moins 
deux \'etats, cet espace est au moins de dimension $2$. On fait donc ce choix dans un premier temps pour voir ce que l'on peut en d\'eduire. C'est dans un deuxi\`eme temps que  
l'exp\'erience viendra confirmer que le mod\`ele est bon.  
 \end{enumerate}
\end{remark}

\noindent Il faut maintenant se poser la question de savoir si la d\'efinition ci-dessus peut\^etre r\'ealis\'ee. 
Autrement dit, a-t-on, dans un espace de dimension $2$ 
des op\'erateurs $\hat{S}_x$, $\hat{S}_y$ et $\hat{S}_z$ v\'erifiant les relations \ref{defS} ? 
La r\'eponse est oui et il est facile de le voir. En effet, soit $\ep_{spin}$ un espace de Hilbert de dimension $2$. 
Si la solution existe, choisissons une base orthonorm\'ee $(\phi^z_+,\phi^z_-)$ 
correspondant aux \'etats propres de $\hat{S}_z$. On a donc $\hat{S}_z (\phi^z_\pm) = \pm \frac{1}{2}\hbar \phi^z_\pm$. 
En \'ecrivant $\hat{S}_{xy}(\phi_\pm)$ dans la base $(\phi^z_+,\phi^z_-)$  et en utilisant les relations 
\eref{defS}, on voit facilement que 
\begin{eqnarray} \label{sanal} 
 \hat{S}_x(\phi^z_\pm)= \frac{\hbar}{2} \phi^z_{\mp} \; \hbox{ et } \, \hat{S}_y(\phi^z_\pm)= \pm \frac{i \hbar}{2} 
\phi^z_{\mp},
\end{eqnarray}
autrement dit, sous forme matricielle: 
\begin{eqnarray} \label{smat} 
 \hat{S}_x = \left( \begin{array}{cc} 0 & \hbar/2 \\ \hbar/2 & 0 \end{array} \right) \; \; ,  \; \;
\hat{S}_y = \left( \begin{array}{cc} 0 & -i\hbar/2 \\ i\hbar/2 & 0 \end{array} \right) \; \; ,  \; \;
\hat{S}_z = \left( \begin{array}{cc}  \hbar/2& 0\\ 0 & - \hbar/2 \end{array} \right).
\end{eqnarray}

\section{Mod\'elisation de l'\'etat spatial d'une particule de spin $1/2$}
D'apr\`es ce qui pr\'ec\`ede, pour mod\'eliser l'\'etat spatial d'une particule, nous choisirons de travailler 
dans l'espace de Hilbert
$$\ep_H = \ep_{externe} \otimes \ep_{spin}$$
o\`u $\ep_{externe}= L^2(\mR^3)$ (l'espace des fonctions \`a valeurs complexes et de carr\'e somma--ble sur $\mR^3$) et o\`u 
$\ep_{spin}$ est un espace de Hilbert de dimension $2$ poss\'edant une base canonique (au signe pr\`es) 
ind\'ependante du temps $(\phi^z_+, \phi^z_-)$ 
correspondant aux \'etats propres de $\hat{S}_z$. Tout vecteur $\psi$ de l'espace de Hilbert 
$\ep_H$ s'\'ecrit alors de  mani\`ere  unique comme 
$$\psi = \psi_+  \otimes \phi^z_+ + \psi_-\otimes \phi^z_-.$$
Autrement dit, en fixant une fois pour toutes $(\phi^z_+,\phi^z_-)$, le vecteur d'\'etat de la particule est 
une fonction d'onde \`a deux composantes 
$$\psi(t)= \left( \begin{array}{c} \psi_+(x,t) \\ \psi_-(x,t)  \end{array} \right). $$
En revenant \`a l'interpr\'etation physique du vecteur d'\'etat, la probabilit\'e de pr\'esence de la particule
sur $\Om \subset \mR^3$ \`a l'instant $t$ est 
$$\int_{\Om} |\psi_+|^2(x,t) + |\psi_-|^2(x,t) dx. $$
Les int\'egrales 
$$\int_{\Om} |\psi_\pm|^2(x,t) dx$$
repr\'esentent la probabilit\'e de d\'etecter la particule dans $\Om$ et de mesurer que la projection de son spin sur l'axe $z$
est $\pm \hbar/2$. 

\begin{remark}
 Le spin est compl\`etement ind\'ependant des grandeurs $x,y, z$. Autrement dit, une observable d\'ependant uniquement
de ces grandeurs n'agira que sur $\ep_{externe}$: elle sera de la forme $\hat{A} \otimes Id_{\ep_{spin} } $ et 
commutera avec toute observable ne d\'ependant que du spin, c'est-\`a-dire de la forme 
$Id_{\ep_{externe}} \otimes \hat{B}$.   
\end{remark}

\section{Moment magn\'etique d'une particule} 
Nous allons bri\`evement parler du moment magn\'etique d'une particule en m\'ecanique quantique pour mettre en 
\'evidence le ''comportement spinoriel'' du spin (voir le paragraphe suivant). 
En physique classique, on d\'efinit le {\em moment magn\'etique} cr\'e\'e par une distribution de courant $j(x)$
par 
$$\vec{\mu}= \frac{1}{2} \int_{\mR^3} x \wedge j(x) dx.$$ 
On peut v\'erifier que cette d\'efinition ne d\'epend pas de l'origine lorsque le syst\`eme est constitu\'e de boucles 
ferm\'ees.
Son interpr\'etation physique est simple: si le syst\`eme est plac\'e dans un champ magn\'etique $\vec{B}$, les boucles de courant
vont avoir tendance \`a tourner pour se placer dans une position perpendiculaire \`a ce champ magn\'etique. On peut
montrer exp\'erimentalement que la force qui 
s'applique sur la boucle de courant est  $\nabla (\vec{\mu}, \vec{B})$. L'\'energie potentielle associ\'ee
est le produit scalaire 
\begin{eqnarray} \label{energiepot}
 E = - (\vec{\mu} , \vec{B}).
\end{eqnarray}
 
\noindent Si maintenant on consid\`ere un \'electron en orbite autour d'un noyau, nous pouvons v\'erifier que 
\begin{eqnarray} \label{momentmagcl}
 \vec{\mu} = 1/2 q x \wedge \vec{v}
\end{eqnarray}
o\`u $\vec{v}$ est sa vitesse et o\`u $q$ est sa charge. Par ailleurs, nous avons expliqu\'e dans le paragraphe 
\ref{momentci} du chapitre \ref{momeci} que son moment cin\'etique \'etait donn\'e par 
\begin{eqnarray} \label{momentcicl}
 \vec{L} = m x \wedge \vec{v}.
\end{eqnarray}
Autrement dit, les moments magn\'etique et cin\'etique sont proportionnels. C'est une remarque qui joue 
un r\^ole important en physique mais nous ne nous attarderons pas sur le sujet ici. 
Cependant, en m\'ecanique quantique, nous avons vu que le spin a une influence sur le comportement d'une particule qui 
s'apparente au mouvement orbital d'un \'electron. Nous choississons donc de lui attribuer un moment magn\'etique 
$\hat{\mu}$  proportionnel \`a $\hat{S}$ mais dont l'int\'er\^et est, comme dans le cas classique, de tenir compte de la
charge de la particule consid\'er\'ee alors que le moment cin\'etique tient compte de sa masse.  
Nous d\'efinissons donc 
\begin{eqnarray}
 \hat{\mu} = \gamma \hat{S}
\end{eqnarray}
o\`u $\gamma= 2 q/m $ o\`u $q$ est la charge de la particule consid\'er\'ee et o\`u $m$ est sa masse, 
ce qui peut \^etre compar\'e avec les formules 
\eref{momentmagcl} et \eref{momentcicl}.

\section{''Comportement spinoriel``} 
En m\'ecanique classique, une rotation de $2 \pi$ agit comme l'identit\'e. Nous allons voir que ce n'est pas le cas en 
m\'ecanique quantique.
Pour cela consid\'erons un champ magn\'etique fixe $\vec{B}_0= B_0 \vec{z}$ parall\`ele \`a l'axe $z$. 
Nous ne nous pr\'eoccupons pas des 
variables spatiales et travaillons donc dans $\ep:= \ep_{spin}$ l'espace de Hilbert 
\`a deux dimensions qui permet de d\'ecrire l'\'etat de spin. Notons $(\phi^z_+,\phi^z_-)$ la base canonique associ\'ee. 
La formule \eref{energiepot} nous sugg\`ere d'\'ecrire 
l'hamiltonien li\'e au syst\`eme comme  
$\hat{H} = - (\hat{\mu}, \hat{B}_0 )$. Il faut donner un sens \`a cette expression: puisque $\vec{B}_0$ est parall\`ele 
\`a l'axe des $z$, prendre le produit scalaire avec $\vec{B}_0$ revient \`a multiplier $\mu_z$ par le nombre $B_0$.
Puisque $\hat{\mu} = \gamma \hat{S}$, on va travailler avec le hamiltonien
\begin{eqnarray}
 \hat{H} = - \gamma B_0 \hat{S}_z= \om_0 \hat{S}_z
\end{eqnarray}
o\`u l'on a not\'e $\om_0= - \gamma B_0$. Les niveaux d'\'energie (donc les valeurs propres de $\hat{H}$) 
de la particule sont donc $\pm \om_0 \hbar / 2$ et les vecteurs propres associ\'es sont $\phi^z_\pm$. 
Supposons qu'\`a $t=0$, le spin de la particule soit dans l'\'etat $\psi(0)= \al \phi^z_+ + \beta \phi^z_-$ avec 
$|\al|^2 + |\be|^2 =1$.  Notons $\psi(t)= \al(t) \phi^z_+ + \beta (t) \phi^z_-$
L'\'equation de Schr\"odinger dit que 
$i \hbar \partial_t \psi =  \om_0 \hat{S}_z \psi(t)$, c'est-\`a-dire
$$i \hbar (  \al'(t) \phi^z_+ + \beta' (t) \phi^z_-) = \om_0 \hbar /2 ( \phi^z_+ - \phi^z_-)$$
ce qui conduit au syst\`eme: 
\[ \left\{ \begin{array}{ccc}
            \al'(t) & = & - i \om_0/2 \al(t), \\
\be'(t) & = & i \om_0/2 \al(t).
           \end{array}
\right. \]
Avec les conditions initiales $\al(0)= \al$ et $\be(0)=\be$, nous obtenons 
\begin{eqnarray} \label{vecteurdetat}
 \psi(t) = \al e^{-i\om_0t/2} \phi^z_+ + \be e^{i \om_0 t/2} \phi^z_-.
\end{eqnarray}

\noindent Nous voulons maintenant calculer les moyennes $<\mu_x>$, $<\mu_y>$ et $<\mu_z>$. En revenant \`a la 
d\'efinition de la mesure d'une grandeur physique en m\'ecanique 
quantique, on a 

$$  \begin{aligned} 
    < \mu_z > & =(\psi, \hat{\mu}_z \psi) \\
& = \gamma \hbar/2 ( \al e^{-i\om_0t/2}  \phi^z_+ + \be e^{i \om_0 t/2} \phi^z_-,   
\al e^{-i\om_0t/2}  \phi^z_+ - \be e^{i \om_0 t/2} \phi^z_-),\\
  \end{aligned}
$$
c'est-\`a-dire
\begin{eqnarray} \label{muzmoy}
 < \mu_z > = \gamma \hbar/2 (|\al|^2 - |\beta|^2 ).
\end{eqnarray}

\noindent 
De m\^eme, en utilisant les relations \eref{sanal}, on calcule  que 
\begin{eqnarray} \label{muxmoy}
 < \mu_x > = C \cos(\om_0 t + \phi)
\end{eqnarray}
et 
\begin{eqnarray} \label{muymoy}
 < \mu_y > = C \cos(\om_0 t + \phi)
\end{eqnarray}
o\`u l'on a \'ecrit $\gamma \hbar \overline{\al} \beta = C e^{i \phi}$. 

\noindent On peut tirer deux conclusions des relations \eref{muzmoy}, \eref{muxmoy} et \eref{muymoy}:
\begin{enumerate}
 \item la premi\`ere est que la moyenne de la composante du 
moment magn\'etique (et donc du spin) dans l'axe parall\`ele au champ magn\'etique ne d\'epend pas du temps alors que
sur des axes orthogonaux, cette moyenne est sinuso\"{\i}dale de p\'eriode $\om_0$. On appelle ce ph\'enom\`ene la 
{\em pr\'ecession de Larmor}. C'est un fait directement   observable par l'exp\'erience et qui permet ainsi 
de confirmer 
le mod\`ele th\'eorique que l'on vient de construire. 
\item D'apr\`es les relations   \eref{muzmoy}, \eref{muxmoy} et \eref{muymoy}, les moyennes des coordonn\'ees du moment
magn\'etique sont $2 \pi/\om_0$-p\'eriodiques. En particulier, l'influence sur le comportement de la particule 
du moment magn\'etique est $2 \pi/\om_0$-p\'eriodique. Or, la construction m\^eme du spin 
sugg\`ere de l'interpr\'eter comme une sorte
de rotation interne et la $2 \pi/\om_0$-p\'eriodicit\'e correspond au moment o\`u cette rotation a  fait un tour complet. 
Cependant, en revenant \`a \eref{vecteurdetat}, on voit que l'\'etat $\psi$ n'est que $4\pi/\om_0$-p\'eriodique
et que $\psi(t+ 2 \pi /\om_0) = -\psi(t)$. Les lecteurs qui connaissent la construction du groupe $spin$ pourront faire
le rapprochement avec cette situation. 
\end{enumerate}



\chapter{Principe de Pauli} 
Dans ce chapitre, nous allons \'etudier l'interaction entre plusieurs particules et aboutirons encore une fois 
\`a des conclusions \'eloign\'ees de nos habitudes de physique classique.

\section{Syst\`emes de deux particules} 
Imaginons un syst\`eme compos\'e de deux (pour simplifier) particules d\'ecrites chacune dans des espaces de Hilbert $\ep_1$ 
et $\ep_2$. Les principes de la m\'ecanique quantique disent que l'on doit \'etudier le syst\`eme dans 
l'espace $\ep_1 \otimes \ep_2$. 
Il y a deux remarques importantes \`a faire: chacune sera l'objet d'un paragraphe. 

\subsection{\'Etats d'intrication} 
 Supposons dans ce paragraphe que les 
$\ep_i$ sont les espaces de dimension $2$ associ\'es au spin 
des deux particules, que l'on suppose donc de spin $1/2$. On a 
une base canonique $(\psi^i_+, \psi^i_-)$ dans chaque espace correspondant aux \'etats propres de la composante $z$ du
spin de chaque particule. Supposons qu'\`a un instant donn\'e l'\'etat du syst\`eme est 
$$\psi(t) = \al \psi^1_+ \otimes \psi^2_- + \beta \psi^1_-  \otimes \psi^2_+.$$
Il faut savoir qu'exp\'erimentalement, on est capable de produire ce 
genre de situation. Cet \'etat est tr\`es particulier parce que l'on peut faire le raisonnement suivant: si on mesure 
la composante en $z$ du spin de la premi\`ere particule,
on trouvera $\hbar/2$ avec probabilit\'e $|\al|^2$ et $-\hbar/2$ avec probabilit\'e $|\be|^2$. 
De m\^eme, pour la seconde particule, on trouvera 
 $\hbar/2$ avec probabilit\'e $|\be|^2$ et $-\hbar/2$ avec probabilit\'e $|\al|^2$. 
L\`a ou cela devient \'etonnant, c'est que la mesure {\bf simultan\'ee} 
de $\hbar/2$ pour chacune des deux particules a une probabilit\'e de $0$: plus pr\'ecis\'ement, on mesurera simultan\'ement 
$\hbar /2$ pour la premi\`ere particule et $-\hbar /2$ pour la seconde avec probabilit\'e $|\al|^2$ et 
 on mesurera
$-\hbar /2$ pour la premi\`ere particule et $\hbar /2$ pour la seconde avec probabilit\'e $|\be|^2$. 
Les deux particules ne sont plus dissoci\'ees: on dit que le syst\`eme est dans un {\em \'etat intriqu\'e}. 
On peut donner une d\'efinition math\'ematique plus pr\'ecise:
\begin{definition}
 On dit que le syst\`eme constitu\'e des deux particules est dans un {\em \'etat intriqu\'e} si le vecteur d'\'etat 
ne peut pas se mettre sous la forme $\psi (t) = \psi^1(t) \otimes \psi^2(t)$. 
\end{definition}
Dans ces situations, les particules d\'ependent l'une de l'autre. Einstein, Podolsky et Rosen se bas\`erent sur ces \'etats
pour formuler un paradoxe connu sous le nom de {\em paradoxe EPR} permettant d'infirmer les principes de la 
m\'ecanique quantique.  
Pour d\'ecrire tr\`es bri\`evement leur argument, ils disaient qu'en mesurant le spin de la premi\`ere particule mais sans 
toucher \`a la deuxi\`eme et en communiquant ensuite le r\'esultat trouv\'e, par exemple $\hbar/2$,  
on \'etait certain de mesurer $-\hbar/2$ pour la deuxi\`eme, contredisant en apparence le fait qu'on avait une probabilit\'e
$|\al|^2$ de trouver $-\hbar/2$. On aboutit \`a des pr\'edictions diff\'erentes en consid\'erant deux syst\`emes 
d'une seule particule
ou au contraire un seul syst\`eme de deux particules. 
L'erreur dans leur raisonnement consiste \`a dire qu'en faisant une mesure sur la premi\`ere particule, on 
ne change rien \`a l'\'etat de la seconde. Les \'etats intriqu\'es sont justement des \'etats o\`u l'on ne peut plus 
dissocier les deux particules et leur attribuer \`a chacun un \'etat bien distinct. 
Le paradoxe EPR sera d'ailleurs contredit par l'exp\'erience. 

\subsection{Indiscernabilit\'e des particules}
Un autre probl\`eme se pose: celui de savoir comment discerner deux particules. La plupart du temps, ce sont les 
crit\`eres physiques qui permettent de le faire. Mais si maintenant les particules sont de m\^eme nature, deux \'electrons
par exemple ? La situation est diff\'erente de celles que l'on rencontre au quotidien: deux objets identiques ont toujours
de petites diff\'erences qui permettent de les diff\'erencier. Au niveau des particules, il en va autrement: 
deux \'electrons sont identiques en tout point. Il reste alors une solution: celle d'attribuer un num\'ero \`a chacune
des particules et de suivre leur trajectoire. Mais l\`a encore, c'est une m\'ethode 
qui ne marche qu'au niveau macroscopique,
les trajectoires des particules \'etant mal d\'efinies. \\

\noindent Ce probl\`eme se lit dans les principes de la m\'ecanique quantique. En effet, prenons une situation physique 
quelconque o\`u  deux particules identiques sont soumises aux m\^emes conditions physiques. Les espaces
$\ep_1$ et $\ep_2$ sont alors \'egaux, de m\^eme que les hamiltoniens $\hat{H}_1$ et $\hat{H}_2$. Notons $\ep := \ep_1=\ep_2$   
et $\hat{H}:= \hat{H}_1=\hat{H}_2$. 
Le hamiltonien du syst\`eme form\'e par les deux particules, dont l'espace de Hilbert associ\'e est
$\ep\otimes \ep$ doit \^etre exactement
$$\hat{H}_{tot} = \hat{H} \otimes Id_{\ep} + Id_{\ep_1} \otimes \hat{H}_2.$$
En effet, les valeurs propres de $\hat{H}_{tot}$ doivent \^etre exactement les sommes des valeurs propres de 
$\hat{H}$ puisque l'\'energie totale du syst\`eme est \'egale \`a la somme de l'\'energie port\'ee par chaque
particule. 
Consid\'erons donc deux niveaux d'\'energie $E_1$ et 
$E_2$ distincts de $\hat{H}$ et notons  $\psi_1$ et $\psi_2$ les vecteurs propres correspondants.  
 On remarque que  la valeur propre $E_1 + E_2$ de $\hat{H}_{tot}$ poss\`ede un espace propre de dimension 
au moins $2$ contenant $\psi_1 \otimes \psi_2$ et $\psi_2 \otimes \psi_1$. La question est maintenant: supposons que l'on 
fasse une mesure de l'\'energie du syst\`eme \`a l'instant $t$ et que l'on trouve $E_1 + E_2$. Les principes de m\'ecanique
quantique disent qu'en $t$, le vecteur d'\'etat $\psi$ doit \^etre un vecteur propre normalis\'e de l'espace propre 
associ\'e \`a $E_1+E_2$. Le probl\`eme est maintenant qu'il n'y a aucune raison de priviligier un choix de vecteur propre
plut\^ot qu'un autre dans cet espace propre. On pourrait par exemple se dire qu'on a num\'erot\'e les particules et que l'on 
postule de prendre $\psi = \psi_1 \otimes \psi_2$ \`a l'instant $t$. Cette mani\`ere de proc\'eder ne tient pas puisque
si on avait juste choisi d'inverser les num\'eros $1$ et $2$ au d\'ebut de l'exp\'erience 
(ce choix n'influe en rien sur le syt\`eme), 
on aboutit \`a une situation diff\'erente \`a l'instant $t$ alors qu'on a fait exactement la m\^eme manipulation du syst\`eme, 
celle de mesurer son \'energie. On pourrait imaginer que ces consid\'erations n'ont aucune influence physique: c'est-\`a-dire
que dans les deux cas, on arrive \`a des vecteurs d'\'etat certes diff\'erents mais dont les diff\'erences ne
se lisent pas dans leurs effets physiques. Malheureusement, l'exp\'erience prouve que ce n'est pas le cas. \\

\noindent Pour r\'esoudre cette ambigu\"{\i}t\'e, il convient de formuler un nouveau principe appel\'e 
{\em principe de Pauli}  et d\'ecrit dans le paragraphe \ref{prpauli} 

\section{Spin d'un syst\`eme  de deux particules}
  Avant de parler du principe de Pauli, il convient d'expliquer comment mod\'eliser le spin de deux particules. D'abord,
en m\'ecanique classique, le moment cin\'etique total de deux syst\`emes est \'egal \`a la somme des moments cin\'etiques des
deux syst\`emes. On va garder cette m\^eme r\`egle pour le spin. Autrement dit, soient deux particules que nous allons supposer
de spin $1/2$ dont le spin est d\'ecrit dans les espaces $\ep^1_{spin}$ et $\ep^2_{spin}$. On va d\'ecrire le spin
du syst\`eme dans l'espace $\ep_s= \ep^1_{spin} \otimes \ep^2_{spin}$ qui est de dimension $4$. Si l'on note 
$\hat{S}^i$ ($i=1,2$) les observables vectorielles de spin associ\'ees, l'observable totale $\hat{S}$ sera
$$\hat{S}  = \hat{S}^1 \otimes Id_{\ep^1_{spin}} + Id_{\ep_{spin}^2} \otimes \hat{S}^2.$$
La raison est la m\^eme que dans le paragraphe pr\'ec\'edent \`a propos de l'\'energie: cela permet que les 
valeurs propres de $\hat{S}$, c'est-\`a-dire les valeurs mesurables possibles, soient les sommes des valeurs 
mesurables pour chacun des deux spins. 
Par abus de notation, nous \'ecrirons 
$$\hat{S}= \hat{S}^1+ \hat{S}^2.$$ 
Notons
 $(\si^{i}_+, \si^{i}_i)$ la base canonique de $\ep^i_{spin}$, 
donc \'etats propres de la composante en $z$ de $\hat{S}^i$. Notons enfin 
$\si_{\pm, \pm} := \si^1_{\pm} \otimes \si^2_\pm$. Les $\si_{\pm,\pm}$ forment une base de $\ep_s$.
De m\^eme que pour le cas d'une particule, nous allons poser 
\begin{eqnarray} \label{hatSt}
 \hat{s}= (\hat{S}_x)^2 +(\hat{S}_y)^2+ (\hat{S}_z)^2 = (\hat{S}^1_x + \hat{S}^2_x)^2 +(\hat{S}^1_y + \hat{S}^2_y)^2 + 
(\hat{S}^1_z + \hat{S}^2_z)^2 
\end{eqnarray}
dont nous notons $S(S+1) \hbar^2$ les valeurs propres. 
De m\^eme, nous param\'etrons les valeurs propres de $\hat{S}_z$ par $M \hbar$. 
Nous noterons $\si_{S,M}$ le vecteur propre commun associ\'e. 
Commen\c{c}ons par regarder les valeurs que peut prendre la norme de $S$ et donc les valeurs propres de $\hat{s}$. 
Puisque les composantes de $S^i$ ne prennent que les valeurs 
$\pm \hbar/2$, les composantes de $S$ ne prennent que les valeurs $0$ ou $\pm \hbar$. Par cons\'equent, 
on peut se dire qu'a priori,  $\| S \|^2$ prendra les valeurs $0$, $\hbar^2$, $2 \hbar^2$ ou $3\hbar^2$. En fait, 
ce n'est pas le cas. En effet, on calcule directement en utilisant les relations \eref{sanal} que 
$$\hat{s} (\si_{+,+}) = 2 \hbar^2 \si_{+,+}, \; \; \hat{s} (\si_{-,-}) = 2 \hbar^2 \si_{-,-}, $$
$$ \hat{s} \left( \frac{1}{\sqrt{2}} (\si_{+,-} + \si_{-,+})\right) =
 2 \hbar^2 \frac{1}{\sqrt{2}} (\si_{+,-} + \si_{-,+}) \; \; \hbox{ et } \; \; 
\hat{s} \left( \frac{1}{\sqrt{2}} (\si_{+,-} - \si_{-,+})\right )= 0.$$
Autrement dit, $\hat{s}$ ne poss\`ede que $2$ valeurs propres: 
\begin{itemize} 
 \item $2 \hbar^2$ correspondant \`a la valeur $S=1$ dont l'espace propre est de dimension $3$ engendr\'e par les vecteurs
$\Th_1:= \si_{+,+}$, $\Th_2:=\si_{-,-}$ et $\Th_3:=  \frac{1}{\sqrt{2}} (\si_{+,-} + \si_{-,+})$;
\item $0$  correspondant \`a la valeur $S=0$ dont l'espace propre est de dimension $1$ engendr\'e par 
$\Th_4:=  \frac{1}{\sqrt{2}} (\si_{+,-} - \si_{-,+})$.
\end{itemize}
Les \'etats de l'espace propre correspondant \`a $S=1$ sont appel\'es {\em \'etats triplets} et celui correspondant \`a 
$S=0$ est appel\'e {\em singulet}. 
Nous conservons les notations du paragraphe pr\'ec\'edent mais nous supposons cette fois que les particules 
sont les m\^emes. On a donc $\ep^1_{spin} = \ep^2_{spin}$ et on peut 
d\'efinir l'{\em op\'erateur d'\'echange} (pour le spin) $P_s$ par 
\[ P_s: \left| \begin{array} {ccc} 
              \ep_s & \to  & \ep_s  \\
u\otimes v & \mapsto & v \otimes u.
             \end{array} \right. \]
Les vecteurs de $\ep_s$ ne s'\'ecrivent  pas tous sous la forme $u \otimes v$ mais on \'etend la d\'efinition par 
lin\'earit\'e \`a $\ep_s$ tout entier. 
Alors, par d\'efinition des vecteurs $\Th_i$, on voit que:
\begin{eqnarray} \label{echange} 
 P_s(\Th_i)= \ep_i \Th_i 
\end{eqnarray}
avec $\ep_i= 1$ si $i=1,2,3$ et $\ep_4 = -1$. Autrement dit, les vecteurs propres
correspondant \`a $S=1$ sont ''totalement sym\'etriques'' tandis que ceux correspondant \`a $S=0$ sont 
``totalement anti-sym\'etrique``  Cette relation va jouer un r\^ole fondamental dans le paragraphe suivant.

\section{Principe de Pauli} \label{prpauli}
Nous avons dit qu'il fallait n\'ecessairement que deux particules de m\^eme nature soient indiscernables. 
Supposons que le syst\`eme compos\'e par ces deux particules soient associ\'e \` a l'espace de Hilbert
$\ep_H= \ep \otimes \ep$, l'espace de Hilbert $\ep$ \'etant lui-m\^eme associ\'e \`a chacune des deux particules. 
D\'efinissons l'op\'erateur d'\'echange entre les deux particules de la m\^eme mani\`ere que ci-dessus, c'est-\`a-dire 
par 

\[ P: \left| \begin{array} {ccc} 
              \ep & \to  & \ep \\
u\otimes v & \mapsto & v \otimes u.
             \end{array} \right. \]
 
\noindent Le principe que l'on veut formuler consiste donc \`a dire que, si $\psi$ est le vecteur d'\'etat du syst\`eme, 
$P(\psi)$ repr\'esente le m\^eme \'etat physique que $\psi$. Remarquons que deux vecteurs d'\'etat $\psi$ et $\psi'$ 
d\'ecrivent la m\^eme situation physique si et seulement si toutes les mesures que l'on peut faire sur le syst\`eme 
donnent les m\^emes r\'esultats, 
autrement dit si pour tout op\'erateur $\hat{A}$, on a 
$$(\psi, \hat{A} \psi) =(\psi',\hat{A} \psi').$$ 
Il est facile de voir que cette condition est r\'ealis\'ee si et seulement si $\psi$ s'\'ecrit $\psi'= e^{i\phi}
\psi$. 
En appliquant cette observation \`a notre situation, on a donc $P(\psi) = e^{i \phi} \psi$ et puisque par ailleurs
$P^2=Id$, $e^{i \phi} = \pm 1$. 
Cela signifie donc que l'on doit poser comme principe que tout vecteur d'\'etat d'un syst\`eme de deux particules 
est soit totalement sym\'etrique soit totalement antisym\'etriques. 
Reste \`a savoir quel est le bon choix \`a faire. D'ailleurs, ce bon choix, s'il existe, ne d\'epend-il que de la nature
des particules ? Autrement dit, un syst\`eme de deux particules identiques doit-il toujours avoir un vecteur d'\'etat 
sym\'etrique, anti-sym\'etrique ou bien les deux situations peuvent-elles arriver ? 
Le principe de Pauli postule la r\'eponse \`a ces questions et modifie  
la d\'efinition des {\em bosons} et {\em fermions} donn\'ee plus t\^ot. 

\noindent {\bf Principe de Pauli: } 
{\em Toutes les particules de la nature appartiennent \`a l'une des cat\'egories suivantes :
\begin{itemize}
 \item les {\em bosons} pour lesquels le vecteur d'\'etat de deux particules identiques est sym\'etri--que (i.e. 
$P\psi = \psi$); 
\item  les {\em fermions} pour lesquels le vecteur d'\'etat de deux particules identiques est anti-sym\'etrique (i.e. 
$P\psi = -\psi$). 
\end{itemize}
Toutes les particules de spin entier ou nul sont des bosons tandis que celles de spin demi-entier sont des fermions.
}

\noindent Il faut bien comprendre que ce principe est un postulat: rien ne nous a permis de 
d\'ecider qu'une particule de spin $1/2$ par exemple \'etait dans un \'etat totalement anti-sym\'etrique plut\^ot 
que sym\'etrique. Cela deviendra un th\'eor\`eme en m\'ecanique quantique relativiste (dont nous ne parlerons pas dans le 
chapitre correspondant) et par ailleurs, l'exp\'erience confirme que ce postulat est correct. \\

\noindent Ce principe r\'eduit l'espace dans lequel on travaille. Par exemple, consid\'erons deux particules de 
spin $1/2$ et regardons seulement leur spin. D'apr\`es ce principe, il doit \^etre dans l'\'etat 
singulet d\'ecrit dans le paragraphe pr\'ec\'edent. Le vecteur d'\'etat de spin doit donc \^etre \'egal \`a 
$\frac{1}{\sqrt{2}} (\si_{+,-} - \si_{-,+})$. Ce qui signifie que, quelle que soit la mesure effectu\'ee sur le spin, 
on trouve les particules dans un \'etat oppos\'e.  C'est pourquoi on parle 
du {\em principe d'exclusion de Pauli}. Il faut savoir que ce principe a un grand nombre de cons\'equences macroscopiques
dont nous ne parlerons pas ici.

\begin{remark}
 On prendra garde de ne pas confondre le fait d'\^etre dans un \'etat anti-sym\'etrique et le principe d'exclusion 
\'evoqu\'e ci-dessus. En effet, il existe des \'etats sym\'etriques qui ont la m\^eme cons\'equence:
 par exemple l'\'etat triplet $\frac{1}{\sqrt{2}} (\si_{+,-} + \si_{-,+})$.
\end{remark}

\chapter{M\'ecanique quantique relativiste}
De par sa construction qui consid\`ere le temps comme d\'ecoupl\'e des variables d'espace, la m\'ecanique quantique 
n'est pas compatible avec les principes de relativit\'e restreinte. On peut d'ailleurs v\'erifier que la th\'eorie n'est
pas invariante par le groupe de Lorentz ce qui signifie que deux observateurs relativistes 
ne sont pas physiquement \'equivalents. Exp\'erimentalement, on observe aussi que, comme on pouvait s'en douter, 
la m\'ecanique quantique n'est pr\'ecise que lorsque les ph\'enom\`enes observ\'es ne mettent en jeu que des particules \`a faible 
vitesse. Elle n'est par exemple pas un bon mod\`ele pour d\'ecrire toutes les exp\'eriences o\`u il y a 
interaction entre lumi\`ere et mati\`ere.  \\

\noindent Nous pr\'esentons dans ce chapitre les premi\`eres tentatives pour modifier la m\'ecanique quantique afin la 
rendre relativiste. C'est ainsi que nous aboutirons \`a l'\'equation de Dirac. C'est cependant une vision tr\`es 
incompl\`ete et qui dans de nombreuses situations o\`u elle ne s'applique pas, a \'et\'e abandonn\'ee au profit 
de la th\'eorie quantique des champs, 
incompl\`ete elle-aussi mais  permettant de r\'eparer certaines incoh\'erences \'evidentes de la 
m\'ecanique quantique relativiste. 
Pour comprendre comment on en est arriv\'es  \`a la th\'eorie quantique des champs, il est n\'ecessaire de 
comprendre cette th\'eorie. 


\section{\'Equation de Klein-Gordon} 
Nous allons d'abord chercher \`a trouver une \'equation relativiste en ne nous pr\'eo--ccupant que des 
variables d'espaces. 
En d'autres termes, nous commen\c{c}ons \`a travailler avec une particule de spin $0$. Dans ce cadre, pour 
b\^atir une th\'eorie 
relativiste,  il est naturel de  travailler dans l'espace de la relativit\'e restreinte, autrement dit, 
$\mR^4$ muni de la forme quadratique $\eta:= dx^2-dt^2 $ de signature $(3,1)$. Rappelons que dans ces coordonn\'ees,
la vitesse de la lumi\`ere vaut $1$. Nous n'allons pas chercher dans ce paragraphe \`a pousser la th\'eorie au maximum.
Nous nous contenterons de soulever quelques probl\`emes pos\'es par la m\'ecanique quantique relativiste. Nous serons
plus pr\'ecis dans le paragraphe suivant. 
Pour commencer, il faut essayer de 
d\'egager le coeur de la m\'ecanique quantique: il s'agit \'evidemment de l'\'equation de Schr\"odinger 
$i \hbar \partial_t \psi = \hat{H} \psi.$
En effet, tous les autres principes ne font que d\'efinir le cadre de travail et donnent le moyen de 
relier la th\'eorie \`a la situation physique. Par contre, c'est cette \'equation 
qui d\'ecrit l'\'evolution du syst\`eme. Nous allons donc commencer par supposer que 
la particule est d\'ecrite l\`a aussi par une fonction d'onde d\'efinie sur $\mR^4$ et essayer de trouver 
une ''\'equation de Schr\"odinger relativiste``.
Pla\c{c}ons-nous  dans la situation la plus simple o\`u la particule n'est pas charg\'ee, et \'evolue dans le vide sans 
contrainte. En m\'ecanique classique et dans ce contexte, l'\'energie du syst\`eme est compos\'e uniquement de son \'energie 
cin\'etique $E= 1/2 m \| \vec{v} \|^2= \| \vec{p} \|^2/2m$ o\`u $\vec{v}$ et $\vec{p}$ sont les vecteurs vitesse et
quantit\'e de mouvement. Cette expression nous a conduits \`a d\'efinir 
l'observable \'energie totale $\hat{H}$ par  (voir \eref{energie2}) 

$$  \hat{H} = \hat{E_c} + \hat{E_p}= - \frac{\hbar^2}{2m} \Delta. $$
 
\noindent Dans le m\^eme situation physique, que nous regardons cette fois d'un oeil relativiste, le calcul de 
l'\'energie totale tiendra compte de l'\'energie au repos de la 
particule donn\'ee par la c\'el\`ebre $E=mc^2$. Nous ne rappellerons pas ici les bases de la relativit\'e restreinte 
et donnons directement la valeur de l'\'energie en fonction du vecteur quantit\'e de mouvement $\vec{p}$ (qui appartient
\`a l'espace vu par l'observateur) et de la masse au repos: on a
\begin{eqnarray} \label{energierel}
 E= \sqrt{\eta(\vec{p},\vec{p}) - m^2}.
\end{eqnarray}
En proc\'edant de la m\^eme mani\`ere que pour obtenir \eref{energie2}, cela conduit \`a poser 
$$\hat{H} = \sqrt{- \hbar^2\Delta + m^2}.$$
 
\noindent Cette d\'efinition n'a pas de sens: pour lui en donner un, plut\^ot que sur \eref{energierel}, il faut 
se baser  sur l'expression 
\begin{eqnarray} \label{energierel2}
E^2= \eta(\vec{p},\vec{p}) - m^2
 \end{eqnarray}
et de ce fait, poser comme \'equation de Schr\"odinger relativiste
$$(i \hbar \partial_t)^2 \psi(x,t)= - \hbar^2 \Delta \psi(x,t)+ m^2\psi(x,t)$$
c'est-\`a-dire, 
\begin{eqnarray} \label{kleing}
 \square \psi + \left(\frac{m}{\hbar}\right)^2 \psi = 0
\end{eqnarray}
o\`u $\square := -\partial_t^2 + \Delta$ est l'op\'erateur d'Alembertien. 
Cette \'equation est appel\'ee {\em \'equation de Klein-Gordon}. On peut v\'erifier qu'elle est invariante par 
le groupe de Lorentz et dans un premier temps en tout cas, elle est un bon candidat pour \^etre une 
\'equation de Schr\"odinger relativiste bien que nous n'ayons pas expliqu\'e encore les principes de la th\'eorie.  

\section{Limites de l'\'equation de Klein-Gordon}
L'\'equation de Klein-Gordon est insatisfaisante \`a plusieurs \'egards:
\begin{enumerate} 
\item D'abord, les relations \eref{energierel} et \eref{energierel2} ne sont pas \'equivalentes. En effet, 
\eref{energierel2} est \'equivalente \`a 
$$E= \pm \sqrt{\eta(\vec{p},\vec{p}) - m^2}.$$
Autrement dit, d'un point de vue math\'ematique, cela conduit \`a des \'energies qui peuvent \^etre n\'egatives. Nous 
\'etudierons plus tard ce qu'il en est au niveau physique. 
\item Supposons qu'un observateur connaisse l'\'etat d'un syst\`eme
\`a l'instant $t=0$, c'est-\`a-dire qu'il connaisse $\psi$ \`a l'instant $t=0$, il ne peut pas en d\'eduire l'\'evolution
du syst\`eme: en effet, l'\'equation de Klein-Gordon \'etant d'ordre $2$ par rapport au temps, il lui faut en plus 
conna\^{\i}tre $\partial_t \psi$ en $t=0$. En particulier, il en sera de m\^eme pour des vitesses faibles. Cela contredit
la m\'ecanique quantique non relativiste pour laquelle la connaissance de l'\'etat du syst\`eme \`a un instant donn\'e
est suffisante pour conna\^{\i}tre son \'evolution.
\end{enumerate} 

\noindent Avant de regarder les cons\'equences physiques d'\'energies n\'egatives, il faut d'abord 
mettre au point la th\'eorie. Nous nous int\'eressons donc au point $(2)$. Faisons ce qui est habituel 
quand on regarde une EDO d'ordre $2$ que l'on voudrait ramener \`a l'ordre $1$ 
(en $t$ seulement): on 
d\'efinit le vecteur 
\[ \phi= \left( \begin{array}{c} 
                  \psi \\ \partial_t \psi
                 \end{array}
\right). \]
On est ramen\'es \`a l'\'equation d'ordre $1$ suivante:
\[ \partial_t \phi = \left( \begin{array}{ccc} 
                             0 & Id \\ \Delta + \frac{m^2}{\hbar^2} & 0  \\
                            \end{array} \right) \phi. \]
En fait, il sera plus commode de poser 
\begin{eqnarray} \label{Phichi} 
 \phi_1= \psi + \frac{i\hbar}{m} \partial_t \psi \; \hbox{ et }  \phi_2 = \psi - \frac{i\hbar}{m} \partial_t \psi
\end{eqnarray}
et de remarquer que la fonction d'onde d\'efinie par $\Phi= (\phi_1, \phi_2)$ v\'erifie l'\'equation
suivante 
\begin{eqnarray} \label{eqphi}
 \partial_t \Phi = \frac{1}{2} \left( 
\begin{array}{ccc} 
                            \frac{i \hbar }{m}  \Delta&  \frac{i \hbar }{m} \Delta + \frac{2im}{\hbar} \\ 
-\frac{i \hbar }{m} \Delta - \frac{2im }{\hbar}& - \frac{i \hbar }{m}  \Delta&\\
                            \end{array} \right)\Phi. 
 \end{eqnarray}
Si la vitesse de la particule est petite devant la vitesse de la lumi\`ere, nous pouvons n\'egliger son 
\'energie cin\'etique devant son \'energie interne et donc l'\'energie totale est \'egale \`a 
$E \simeq mc^2=m^2$. En raisonnant de la m\^eme mani\`ere que pour obtenir \eref{energie2} et \eref{energierel2}, cette 
\'egalit\'e se traduit en terme d'observables par 
$i \hbar \partial_t  \psi = m \psi$ et donc pour des situations non-relativistes, $\phi_2 \simeq 0$. 
En prenant $\phi_2= 0$ et en regardant la premi\`ere coordonn\'ee dans 
\eref{eqphi}, on obtient 
$$ \partial_t \phi_1 =   \frac{i \hbar}{2m}\Delta \phi_1.$$
Autrement dit, on retrouve l'\'equation de Schr\"odinger non-relativiste. \\

\noindent Il convient maintenant de faire une remarque importante: dans le cadre non-relativiste, l'interpr\'etation 
de la 
fonction d'ondes comme une densit\'e de probabilit\'e d\'epend de mani\`ere \'evidente de la forme produit 
(espace * temps) de l'espace-temps. En effet, il suffit d'imposer 
que sur les espaces $t=constante$, la fonction d'onde  soit de norme $1$. 
Le premier  probl\`eme que nous 
rencontrons ici est qu'en relativit\'e restreinte, le temps n'est pas universel et poser $t=constante$ avec la 
coordonn\'ee
$t$ de $\mR^4$ revient \`a choisir un observateur particulier. C'est en fait un faux probl\`eme: en effet, 
l'id\'ee sera de construire une densit\'e de probabilit\'e $P(x,t)$ telle que $\int_{\mR^3} P(x,0) dx= 1$ 
Ensuite, en se d\'ebrouillant bien (voir ci-dessous), on montre que l'\'equation de Klein-Gordon implique 
$\int_{\mR^3} P(x,T)^2 dx= 1$ pour tout $T$. L'\'equation de Klein-Gordon \'etant invariante par le groupe de Lorentz,
cette  conservation de norme sera vraie quel que soit la variable de temps $t$ choisie, pour peu qu'elle corresponde \`a
un observateur galil\'een. \\

\noindent Venons-en au deuxi\`eme probl\`eme: on veut obtenir la conservation de norme \'evoqu\'ee ci-dessus. 
Notons $(\cdot,\cdot)$ le produit scalaire de $L^2$. En m\'ecanique quantique non-relativiste, on posait 
$P(x,t) = |\psi|^2$ et l'\'equation de Schr\"odinger impliquait que 
$$\partial_t \int_{\mR^3} P(t,x) dx= \int_{\mR^3} (\partial_t \psi) \overline{\psi} +  
 \psi (\partial_t\overline{\psi}) dx= 0.$$
Si maintenant nous consid\'erons une solution de $\psi$ de l'\'equation de Klein-Gordon \eref{kleing}, 
et si nous posons 
$P(x,t) = |\psi|^2$, le m\^eme raisonnement ne fonctionne plus. Il ne fonctionne que lorsque la fonction 
d'onde est soumise \`a une \'equation du type $\partial_t \psi= \hat{A} \psi$ o\`u $\hat{A}$ est auto-adjoint. 
L'\'equation de Klein-Gordon n'est pas sous cette forme et il faut s'y prendre autrement. La premi\`ere possibilit\'e
consiste \`a poser:
\begin{eqnarray} \label{probas}
 P(x,t) = \frac{i \hbar }{m} (\partial_t \psi) \overline{\psi} +  
 \psi (\partial_t\overline{\psi})
\end{eqnarray}
o\`u la constante $\frac{i \hbar}{m}$, qui n'a aucune importance dans le raisonnement,  a \'et\'e choisie
 pour retrouver la situation non-relativiste en approximation. 
Avec cette expression et l'\'equation de Klein-Gordon, on obtient
$$ \begin{aligned} 
  \partial_t \int_{\mR^3} P(t,x) dx & = \int_{\mR^3} (\partial_tt \psi) \overline{\psi}  
 \psi (\partial_tt \overline{\psi}) dx \\
& =   \int_{\mR^3} (\Delta \psi + \frac{m^2}{\hbar^2} \psi)  \overline{\psi} - 
(\Delta  \overline{\psi} + \frac{m^2}{\hbar^2}  \overline{\psi}) \psi dx \
& = 0
\end{aligned} 
$$
o\`u la derni\`ere \'egalit\'e s'obtient en int\'egrant par partie. La formule \eref{probas} 
 pose de toute mani\`ere un probl\`eme 
important: on peut v\'erifier que la densit\'e de probabilit\'e $P(x,t)$ d\'efinie sous cette forme 
n'est pas positive partout. 
D'apr\`es Pauli et Weisskopf \cite{pw:34}, cela peut n\'eanmoins se r\'esoudre en modifiant 
l'interpr\'etation de $P(x,t)$ mais leur th\'eorie soul\`eve de nombreux probl\`emes. \\

\noindent La deuxi\`eme possibilit\'e
est de travailler en d\'efinissant $P(x,t)= \| \Phi \|^2$ o\`u $\Phi$ est  solution du syst\`eme 
\eref{Phichi}. Nous ne rentrerons pas dans les d\'etails dans ce paragraphe. L'\'equation de Dirac que nous 
\'etudierons  dans le 
prochain paragraphe se base sur le m\^eme principe: travailler avec une \'equation d'ordre $1$. Cela permet de 
d\'efinir la densit\'e de probabilit\'e comme la norme de la fonction d'onde. \\

\noindent Au final, quoi que l'on fasse, deux probl\`emes vont subsister:
\begin{enumerate} 
 \item l'existence de solutions d'\'energie n\'egative: c'est ce qui conduira Dirac \`a postuler l'existence 
d'un ''positron``, particule analogue 
\`a l'\'electron mais charg\'ee positivement (voir prochain paragraphe);
\item le probl\`eme de conservation du nombre de particules que nous n'avons pas encore \'evoqu\'e mais qui est pourtant
le plus fondamental:  exp\'erimentale--ment, on sait que le nombre de particules
d'un syst\`eme n'est pas une quantit\'e conserv\'ee. Or par construction m\^eme de la m\'ecanique quantique o\`u l'on 
construit l'espace de travail en fonction des degr\'es de libert\'e de chaque particule, la th\'eorie n'est pas 
adapt\'ee \`a la situation. C'est ce point particulier qui fera que la m\'ecanique quantique relativiste est 
abandonn\'ee dans de nombreuses situations au profit de la th\'eorie quantique des champs.  
\end{enumerate}

\section{L'\'equation de Dirac}
Nous allons maintenant essayer de construire une th\'eorie relativiste des particules de spin $1/2$, disons de l'\'electron
pour fixer les choses, c'est-\`a-dire qui tient compte 
du spin. Nous nous placerons dans une situation sans champ \'electromagn\'etique. De ce fait, nous ne verrons pas
les effets de la charge de la particule et n'en tiendrons donc pas compte. 
Les deux  probl\`emes \'evoqu\'es \`a la fin du paragraphe pr\'ec\'edent ne seront pas r\'esolus ici: 
comme nous l'avons dit, seule la th\'eorie quantique des champs sera \`a m\^eme de lever le probl\`eme de 
non-conservation du nombre de particules.  
Nous allons travailler  dans un espace de Hilbert 
$\ep_{exterieur} \otimes \ep_{spin}$. D'apr\`es le mod\`ele non-relativiste, il est naturel de consid\'erer que 
$\ep_{spin}$ est de dimension finie. La composante de spin \'etant vectorielle, elle comportera  
en relativit\'e restreinte quatre composantes, c'est-\`a-dire quatre observables l\`a o\`u il n'y en avait que 
trois dans la th\'eorie non-relativiste. C'est l'\'etude de ces quatre observables qui fixera la dimension
de $\ep_s$ mais cette discussion est repouss\'ee \`a plus tard. En tout cas, \`a ce point, 
en fixant une base $(\phi_1,\cdots,\phi_r)$ de l'espace $\ep_{spin}$, on pourra consid\'erer le vecteur d'\'etat $\psi=(\psi_1,\cdots,\psi_r)$ 
comme un vecteur \`a $r$ composantes, o\`u $r$ est la dimension de $\ep_{spin}$ et o\`u chaque composante est une 
fonction d'onde de l'espace $\ep:=\ep_{exterieur} = L^2(\mR^3)$. Le vecteur $\psi \in  \ep_{exterieur} 
\otimes \ep_{spin}$ 
original se retrouvera en \'ecrivant $\psi = \sum_{i=1}^{r} \psi_i \otimes \phi_i$. Dans la suite, nous jonglerons 
avec ces deux mani\`eres  de voir les choses. 
 Pour obtenir un mod\`ele satisfaisant, il faudra que le vecteur d'\'etat $\psi$ 
soit soumis \`a une \'equation g\'en\'eralisant celles de Schr\"odinger (qui ne tient pas compte 
des ph\'enom\`enes relativistes) et de Klein-Gordon (qui ne tient pas compte du spin) et qui devra avoir 
deux propri\'et\'es principales: 
\begin{enumerate}
\item elle devra \^etre invariante sous l'action du groupe de Lorentz;
 \item elle devra \^etre d'ordre $1$ en $t$ et plus pr\'ecis\'ement de la forme 
\begin{eqnarray} \label{formedirac}
  i\hbar \partial_t \psi = \hat{H}_D \psi
\end{eqnarray}
o\`u $\hat{H}_D$ est un op\'erateur auto-adjoint. En d\'efinissant la densit\'e de probabilit\'e de pr\'esence par 
$P(x,t) =\sum_{j= 1}^r \| \psi_j\|^2(x,t)$, la forme de l'\'equation montrera que 
$\partial_t \int_{\mR^3} P(x,t) dx=0$. La d\'emonstration est la m\^eme que celle du cadre non-relativiste. 
\end{enumerate}
La forme explicite de $\hat{H}_D$ d\'ependra bien \'evidemment du syst\`eme physique \'etudi\'e. Comme toujours, 
ce sont les situations les plus simples qui servent de guide pour construire des mod\`eles de syst\`emes plus 
\'elabor\'es. C'est pourquoi ici, 
 nous allons nous placer dans un syst\`eme physique compos\'e
d'un seul \'electron qui n'est soumis \`a aucun potentiel, qui pourrait d'ailleurs briser l'invariance sous l'action
du groupe de Lorentz.  
Il sera plus commode d'utiliser les notations $(x_1:=x,x_2:=y,x_3:=z,x_4:=t)$ pour simplifier les expressions. 
L'observable $\hat{H}_D$ repr\'esentera, comme dans le cas non-relativiste, l'\'energie totale du syst\`eme. 
Il y a plusieurs remarques \`a faire: en ce qui concerne la partie ''spatiale'', 
l'\'energie du syst\`eme d\'ependra \'evidemment de la quantit\'e 
de mouvement $\vec{p}=(p_1,p_2,p_3)$ (on se place du point de vue d'un observateur galil\'een qui ``voit'' les espaces 
$t=constante$) mais pas de la position $(x_1,x_2,x_3)$ puisque l'\'equation doit \^etre invariante par les translations. 
Autrement dit, $\hat{H}_D$ d\'ependra des 
observables $\hat{p}_j= -i \hbar \partial_j$ mais pas des $x_i$. De plus, 
le groupe de Lorentz, contrairement au groupe $O(3)$ agit
aussi sur la quatri\`eme composante. Dans l'\'equation \eref{formedirac}, le terme qui contient l'op\'erateur 
$\partial_t$ est lin\'eaire.  Il doit en \^etre de m\^eme pour les termes contenant les op\'erateurs $\partial_j$, 
et donc les op\'erateurs $\hat{p}_j$.   
On remarque enfin que l'\'energie au repos est une fonction lin\'eaire de la masse et c'est pourquoi nous supposerons que
le terme contenant $m$ doit aussi \^etre lin\'eaire. 
Pour finir, on va faire l'hypoth\`ese raisonnable que $\hat{H}_D$ se met sous la forme 
$$\hat{H}_D= \hat{H}_{exterieur} \otimes \hat{H}_{spin}$$ 
o\`u les observables $\hat{H}_{exterieur}$ et $\hat{H}_{spin}$ agissent respectivement sur $\ep_{exterieur}$  
et $\ep_{spin}$. 
Les remarques que nous avons faites ci-dessus conduisent \`a faire l'hypoth\`ese que 
$\hat{H}_D$ est de la forme 
\begin{eqnarray} \label{formehd}
 \hat{H}_D = \sum_{j=1}^3 \hat{\si}_j \hat{p} + m \si_4  
\end{eqnarray}
 o\`u l'op\'erateur $\hat{m}$ associ\'e \`a l'\'energie au repos est confondu avec $m$ (la multiplication par $m$) et o\`u
l'observable vectorielle 
\[ \hat{\si}= \left( \begin{array}{c} 
                   \hat{\si}_1 \\ \hat{\si}_2 \\ \hat{\si}_3 \\ \hat{\si}_4. 
                 \end{array} \right) \]
agit uniquement sur $\ep_{spin}$. 
Les arguments ci-dessus n'ont aucune rigueur mais comme toujours en physique, on propose des mod\`eles et l'exp\'erience
viendra confirmer leur efficacit\'e: c'est ce que se passe ici. \\
 
\noindent Remarquons tout de m\^eme que cela suppose l'existence d'une observable vectorielle de dimension $4$ agissant 
sur l'espace $\ep_{spin}$, ce qui n'est pas totalement d\'eraisonnable sachant que nous travaillons cette fois sur 
$\mR^4$. Il est naturel de supposer par ailleurs que l'\'equation de Klein-Gordon est v\'erifi\'ee lorsqu'on 
ne consid\`ere que les 
parties spatiales de $\hat{H}_D$. C'est-\`a-dire que l'on veut que 
\begin{eqnarray} \label{kgd}
 - \hbar^2 \partial_{tt} \psi = \left( \sum_{i=1}^2 \hat{p}_i^2 + \frac{m^2}{\hbar^2} \right) \psi.
\end{eqnarray}

\noindent Des relations \eref{formedirac} et \eref{formehd}, on tire que 
$$\left( i \hbar \partial_t \psi -   \sum_{j=1}^3 \hat{\si}_j \hat{p} - m \si_4 \right) \psi=0.$$
Et donc 
$$\left( i \hbar \partial_t \psi +   \sum_{j=1}^3 \hat{\si}_j \hat{p} + m\hbar \si_4 \right) 
\left( i \hbar \partial_t \psi -   \sum_{j=1}^3 \hat{\si}_j \hat{p} + m\hbar \si_4 \right) \psi=0,$$
c'est-\`a-dire, puisque les observables de spin commutent avec les observables externes: 

$$ \begin{aligned} 
 \Big[ -\hbar^2 \partial_{tt} \psi - \sum_{j=1}^3 \hat{\sigma}_j^2 \hat{p}_j^2 - m^2 \hat{\si}_4^2 - 
\sum_{1 \leq j< l \leq 3} & (\hat{\si}_j \hat{\si}_l +\hat{\si}_l  
\hat{\si}_j) \hat{p}_j \hat{p}_l \\  
& - \sum_{j= 1}^3    (\hat{\si}_j \hat{\si}_4 +\hat{\si}_4 
\hat{\si}_j) m \hat{\si}_4 \Big] \psi = 0.
   \end{aligned}$$

\noindent En comparant cette expression avec \eref{kgd}, nous obtenons que, pour tous 
$j,l \in \{1,\cdots4\}$, $j \not=l$,  
\begin{eqnarray} \label{relationspin}
 \hat{\al}_j^2 = 1 \; \hbox{ et } \;   \hat{\si}_j \hat{\si}_l +\hat{\si}_l 
\hat{\si}_j = 0. 
\end{eqnarray}
L'\'equation \eref{formedirac}, qui se r\'e\'ecrit sous la forme 
\begin{eqnarray} \label{eqdirac} 
 - i\hbar \partial_t \psi =   (\sum_{j=1}^3 \hat{\si}_j \hat{p} + m \si_4 ) \psi
\end{eqnarray}
o\`u les $\hat{\si}_j$ v\'erifient \eref{relationspin},  est appel\'ee {\em \'equation de Dirac}.

\section{\'Ecriture spinorielle de l'\'equation de Dirac} 
L'adjectif ``spinorielle`` qui appara\^{\i}t dans le titre de ce paragraphe ne signifie pas ``relatif au spin'' puisque 
de toutes mani\`eres, l'\'equation de Dirac est l\`a pour mod\'eliser le spin de l'\'electron. Il signifie que 
l'on  va exprimer l'\'equation de Dirac avec l'objet math\'ematique ``spineur`` que nous d\'ecrivons tr\`es 
rapidement dans l'appendice \ref{spineurs}. 
 Reprenons les notations du paragraphe pr\'ec\'edent et posons pour $i \in \{ 1,2,3 \}$,  
$$\hat{\gamma}_i:= \hat{\si}_4 \hat{\si}_i \; \; \hbox{ et } \hat{\gamma}_4 := \hat{\si}_4.$$
D'apr\`es les relations \eref{relationspin}, on voit que 
\[ \hat{\gamma}_i^2 = \left|  \begin{array}{ccc} 
 - 1 & \hbox{ si } & i \in   \{ 1,2,3 \} \\
 1 & \hbox{ si } & i=4 \end{array} \right. \]
et que, si $i \not= j$  
$$\hat{\gamma}_i \hat{\gamma}_j + \hat{\gamma}_j \hat{\gamma}_i = 0.$$
On reconna\^{\i}t ici les relations \eref{defclif}  de l'appendice \ref{spineurs} v\'erifi\'ees par les op\'erateurs 
$\rho(e_i)$ agissant sur $\Si_4$. Puisque $\Si_4$ est de dimension minimale 
parmi les espaces sur lesquels des op\'erateurs v\'erifiant de telles relations agissent, on peut raisonnablement
poser $\ep_{spin}:= \Si_4$ (rappelons que le spin $1/2$ correspond \`a la valeur minimale de spin non nulle et donc \`a 
l'espace $\ep_{spin}$ de plus petite dimension) 
et poser $\hat{\gamma}_i:= \rho(e_i)= e_i \cdot$ o\`u l'on a muni $\mR^4$ de la base 
canonique et o\`u $e_4$ correspond \`a la coordonn\'ee $t$. 
On travaille donc maintenant sur l'espace $\ep:= L^2(\mR^4) \otimes \Si_4$. Un vecteur d'\'etat s'\'ecrit 
$$\psi = \sum_{j=1}^4 \psi_j \otimes \th_j$$
o\`u l'on a choisi une base orthonorm\'ee $(\th_1,\cdots,\th_4)$ de $\Si_4$ et o\`u  les $\psi_j$ sont des fonctions 
de carr\'e int\'egrable et \`a valeurs complexes. Pour tous $(x,t) \in  \mR^4$ et $j$, $\psi_j(x,t) \in \mC$ et donc
$$\psi(x,t) = \sum_{j=1}^4 \psi_j(x,t) \otimes \th_j= 1 \otimes (\sum_{j=1}^4 \psi_j(x,t) \th_j)$$
ce qui fait que $\psi$ peut se voir comme  une fonction de $\mR^4$ \`a valeurs dans $\Si_4$. 
Avec ces notations, en multipliant \`a gauche par $\hat{\si}_4$, l'\'equation de Dirac se r\'eecrit 
$$i \hbar \sum_{j=1}^4 e_j \cdot \partial_j \psi - m\psi = 0$$
c'est-\`a-dire avec les notations de l'appendice \ref{spineurs}

\begin{eqnarray} \label{Diracspin}    
(i \hbar D(\psi) - m) \psi = 0. 
\end{eqnarray} 

\noindent L'\'equation mise sous cette forme est sans doute plus difficile \`a interpr\'eter d'un point de vue physique mais 
\`a l'avantage de faire intervenir des objets bien connus des math\'ematiciens. En particulier, l'op\'erateur $D$ se d\'efinit 
parfaitement de mani\`ere intrins\`eque sur les vari\'et\'es  lorentzienne. L'espace  
$(\mR^4,\eta)$ \'etant une vari\'et\'e lorentzienne, l'\'equation de Dirac ne d\'epend pas des coordonn\'ees choisies et est donc invariante par le groupe de Lorentz.   

\section{Une construction plus physique de l'espace $\ep_{spin}$}
L'\'ecriture spinorielle du paragraphe pr\'ec\'edent a l'avantage de donner une forme simple \`a l'\'equation
de Dirac mais se pr\^ete mal \`a son interpr\'etation physique. Nous allons essayer de pr\'esenter
la situation sous un angle diff\'erent. 
\`A la place de construire les op\'erateurs $\hat{\gamma}_i$, nous aurions pu poser 
$$\hat{\al}_1:= -i \hat{\si}_2 \hat{\si}_3, \; \hat{\al}_2:= -i \hat{\si}_3 \hat{\si}_1, \hat{\al}_3:= -i \hat{\si}_1 \hat{\si}_2$$
et 
aussi 
$$\hat{\be}_1 := - i\hat{\si}_1\hat{\si}_2 \hat{\si}_3, \; \hat{\be}_2 := - i \hat{\si}_4 \hat{\si}_1\hat{\si}_2 \hat{\si}_3 \;
\hbox{ et } \hat{\be}_3:= \hat{\si}_4.$$
De nouveau avec les relations \eref{relationspin}, on voit que 
\begin{itemize} 
\item les op\'erateurs $\hat{\al}_j$ commutent avec les op\'erateurs $\hat{\beta}_j$: on a donc envie de chercher $\ep_{spin}$ 
sous la forme 
$\ep_{spin} = \ep_{spin}^1 \otimes \ep_{spin}^2$ o\`u les op\'erateurs $\hat{\al}_j$ et $\hat{\be}_j$ 
agissent respectivement 
sur $\ep_{spin}^1$ et $\ep_{spin}^2$.
\item Les op\'erateurs vectoriels $(\hat{\al}_1,\hat{\al}_2,\hat{\al}_3)$ et   $(\hat{\be}_1,\hat{\be}_2,\hat{\be}_3)$ satisfont
aux m\^emes relations de commutation qu'un spin, c'est-\`a-dire qu'elles satisfont aux relations \eref{defS}.
\end{itemize}

\noindent  Le raisonnement men\'e dans le cadre non-relativiste s'applique alors sur chaque $\ep_{spin}^j$. On va donc d\'ecider 
de prendre $\ep_{spin}^1= \ep_{spin}^2= \ep_s$ o\`u $\ep_s$ est l'espace de Hilbert \`a deux dimensions avec lequel on 
a d\'ej\`a travaill\'e pour d\'efinir le spin $1/2$ en m\'ecanique quantique non-relativiste. 
Par ailleurs, les d\'efinitions des $\hat{\al}_j$ ne font intervenir que les coordonn\'ees spatiales $\hat{\si}_j$ 
avec $j \in \{1,2,3 \}$. On va lui faire jouer un r\^ole particulier ce qui permettra de retrouver 
l'interpr\'etation non-relativiste du spin. 
Pour cela, il suffit de voir que $\hat{\si}_4= \hat{\be}_3$ qui agit sur $\ep_{spin}^2 =\ep_2$ v\'erifie 
$\hat{\be}_3^2= 1$ et donc $\ep_{spin}^2= \hbox{Vect}(u_+,u_-)$ o\`u $(u_+,u_-)$ est une base orthonorm\'ee de 
$\ep_{spin}^2$ compos\'ee des vecteurs propres $u_\pm$ associ\'es aux valeurs propres $\pm 1$  de $\hat{\be}_3$. 
On peut donc \'ecrire 
$$\ep_{spin} = \ep^+ \oplus \ep_-$$
o\`u $\ep_\pm:= \ep_s \otimes \hbox{Vect}(u_\pm)$. 
On peut donc aussi \'ecrire 
$$\ep = (\ep_{exterieur} \otimes \ep_+) + (\ep_{exterieur} \otimes \ep_-).$$
Via cette d\'ecomposition, le vecteur d'\'etat $\psi$ se d\'ecompose en la somme 
$$\psi= \psi_+ + \psi_-$$
o\`u $\psi_\pm$ s'\'ecrivent $\tilde{\psi}_\pm \otimes u_\pm$ et sont donc des vecteurs propres associ\'es
\`a la valeur propre $\pm 1$ de l'observable  $\hat{\si}_4$. 
Puisque chacun des deux espaces $\ep_{exterieur} \otimes \ep_\pm$ est de dimension $2$, 
on peut convenir que les deux premi\`eres 
coordonn\'ees de $\psi$ correspondent \`a $\psi_+$ et les deux derni\`eres \`a $\psi_-$ et 
r\'e\'eecrire la notation initiale $\psi= (\psi_1,\psi_2,\psi_3,\psi_4)$ de la mani\`ere suivante 
\[\psi = \left( \begin{array}{c} \psi_+ \\ \psi_-  \end{array} \right) \] 
o\`u 
   \[\psi_+ = \left( \begin{array}{c} \psi_1 \\ \psi_2  \end{array} \right) \: \hbox{ et } 
\psi_- =\left( \begin{array}{c} \psi_3 \\ \psi_4  \end{array} \right).  \] 

\noindent Venons-en maintenant \`a l'interpr\'etation physique de $\psi_+$. Nous discuterons rapidement de celle 
de $\psi_-$ dans le paragraphe suivant  mais elle est beaucoup moins claire et les physiciens 
ne s'accordent pas sur ce point. \\ 

\noindent 
Reprenons l'\'equation de Dirac \eref{eqdirac}
$$- i \hbar \partial_t \psi = (\sum_{j=1}^3 \hat{\si}_j \hat{p}_j + \hat{\si}_4 m)$$
et pla\c{c}ons-nous \`a un niveau d'\'energie $E$. Autrement dit, consid\'erons un vecteur propre $\phi$ v\'erifiant 
\begin{eqnarray} \label{phipm}
(\sum_{j=1}^3 \hat{\si}_j \hat{p}_j + \hat{\si}_4 m)\phi = E \phi
\end{eqnarray}
que nous d\'ecomposons en $\phi= \phi_+ + \phi_- \simeq (\phi_+, \phi_-)$ associ\'e \`a la d\'ecomposition de $\ep$. 
Par construction,
$\hat{\si}_4 \phi_\pm = \pm \phi_\pm$. 
Par ailleurs, pla\c{c}ons-nous dans des conditions non-relativistes: 
la quantit\'e de mouvement est petite par rapport \`a la masse et donc $E\simeq m$. 
En reportant dans \eref{phipm} o\`u l'on n\'eglige les termes faisant appara\^i{\i}tre la 
quantit\'e de mouvement, on obtient 
$$\hat{\si}_4 m (\phi_+,\phi_-) \simeq m (\phi_+,\phi_-),$$ 
 c'est-\`a-dire 
 $$m (\phi_+, - \phi_-) \simeq m(\phi_+,\phi_-)$$ 
ce qui implique qu'\`a la limite non-relativiste $\phi_- \simeq 0$. C'est la raison pour laquelle on consid\`ere que 
le vecteur  $\psi_+$ s'interpr\`ete comme le vecteur
d'\'etat non-relativiste de la particule et contient tout l'information sur le spin. Plus  pr\'ecis\'ement, 
le spin vectoriel de la particule est donn\'e par 
l'observable vectorielle $(\hat{\al}_1, \hat{\al}_2, \hat{\al}_3)$.  \\
 
\section{De la th\'eorie du positron \`a la th\'eorie quantique des champs} 
Dans ce paragraphe, nous allons nous int\'eresser de plus pr\`es au probl\`eme des \'energies n\'egatives dont nous 
avons d\'ej\`a parl\'e plus haut. 
Pour commencer, pla\c{c}ons-nous dans la m\^eme situation qu'\`a la fin du paragraphe pr\'ec\'edent et consid\'erons
une situation non-relativiste. Nous avons vu alors qu'un \'etat d'\'energie \'etait d\'ecrit par un vecteur propre 
$(\psi_+,\psi_-)$ o\`u la composantes $\psi_-$ est tr\`es petite devant la composante $\psi_+$ (on parle d'ailleurs de
{\em grandes et petites composantes}). Rappelons que, avec les notations du paragraphe pr\'ec\'edent, les vecteurs 
$\psi_\pm$ s'\'ecrivent $\psi_\pm = \tilde{\psi}_\pm  \otimes u_\pm$. Consid\'erons maintenant le 
vecteur  $\tilde{\psi}:= (\tilde{\psi}_- \otimes u_+, \tilde{\psi}_+ \otimes u_-)$. D'apr\`es ce qui pr\'ec\`ede, la composante $\tilde{
\psi}_- \otimes u_+$ 
est petite par rapport \`a $  \tilde{\psi}_+ \otimes u_-$. En n\'egligeant encore une fois la quantit\'e de mouvement, 
ce vecteur est cette fois propre pour l'\'energie $\hat{H}_D \simeq m\hat{\si}_4$ mais avec valeur propre $-m$. 
Si l'on ne n\'eglige plus les termes de quantit\'e de mouvement, $\tilde{\psi}$  n'est pas n\'ecessairement un vecteur 
propre mais on obtient malgr\'e tout que 
$$( \tilde{\psi}, \hat{H}_D \tilde{\psi}) <0.$$
Cela sugg\`ere deux remarques:
\begin{itemize} 
 \item D'abord, il existe des \'etats d'\'energie n\'egative. Nous l'avons d\'ej\`a \'evoqu\'e mais nous en obtenons ici
une preuve. On peut en fait \^etre beaucoup plus pr\'ecis et d\'ecrire pr\'ecis\'ement les \'etats d'\'energie n\'egative,
ce qui ne nous sera pas utile ici. 
\item  La vecteur d'\'etat \`a deux composantes $\psi_-$ semble correspondre \`a une particule qui se comporte comme l'\'electron
mais dont la masse est n\'egative.  
\end{itemize}

\noindent La th\'eorie permet de pr\'edire assez pr\'ecis\'ement le comportement que pourrait avoir une 
particule de masse n\'egative. Mais ce comportement n'a jamais \'et\'e observ\'e exp\'erimentalement. 
Par contre, des ph\'enom\`enes similaires au niveau de la charge ont \'et\'e observ\'es. Pour simplifier les choses,
nous avons ignor\'e la charge mais avec des raisonnements similaires, nous pouvions mettre en \'evidence des 
\'etats correspondant \`a des ''\'electrons'' de masse $+m$ mais de charge $-e>0$ donnant lieu \`a des \'energies 
n\'egatives.    \\

\noindent Ces probl\`emes pourraient n'avoir aucune importance en postulant simplement que ces \'etats ne sont 
que math\'ematiques et que les vecteurs d'\'etat acceptables sont ceux qui appartiennent au 
sous-espace engendr\'e par les vecteurs propres  de l'\'energie correspondant \`a des valeurs 
propres positives. De cette mani\`ere, le vecteur $\tilde{\psi}$ d\'efini ci-dessus ne correspond simplement
\`a aucun \'etat physique. Malheureusement, en regardant ce qui se passe de plus pr\`es, on peut voir qu'un \'electron
plac\'e dans l'\'etat d'\'energie positif le plus faible a une probabilit\'e non nulle de basculer vers un 
\'etat d'\'energie
n\'egative. Pour donner une interpr\'etation physique \`a ce ph\'enom\`ene, Dirac 
imagine la th\'eorie suivante: dans le 
vide, tous les \'etats d'\'energie n\'egative sont satur\'es par un \'electron. Dans un syst\`eme physique, il peut 
arriver que l'un de ces \'electrons s'\'echappe, formant ainsi un ''trou``. Ce 
trou se comporte alors exactement comme une particule de masse $m$ et de charge $-e$ que l'on appelle {\em positron}. 
Cette mani\`ere de pr\'esenter les choses est un progr\`es \'evident qui permet d'expliquer bon nombre de faits 
exp\'erimentaux et qui permet de comprendre la cr\'eation/annihilation de paires de particules. Remarquons qu'il s'agit bien de paires de particules puisque quand un \'electron passe d'un \'etat d'\'energie n\'egative \`a un \'etat d'\'energie positive, l'\'electron lui-m\^eme appara\^{\i}t ainsi que le trou qu'il laisse, c'est-\`a-dire  son anti-particule.  Malheureusement, 
cette th\'eorie  poss\`ede aussi de gros inconv\'enients: lorsqu'un niveau d'\'energie n\'egative est occup\'e ou lib\'er\'e, le
nombre de particules que l'on observe n'est plus le m\^eme, ce qui est mal d\'ecrit par la m\'ecanique quantique 
relativiste puisque par d\'efinition, elle attribue un espace de Hilbert \`a chaque particule.
Il faut, pour rendre compte de nombreux ph\'enom\`enes, b\^atir une th\'eorie qui permette la cr\'eation/annihilation de 
particules dans sa formulation. C'est le but de la th\'eorie quantique des champs. Pour finir, disons que la th\'eorie quantique 
des champs elle-m\^eme est incompl\`ete. Quand on la regarde dans le cadre de la relativit\'e g\'en\'erale, elle 
dit que la notion de vide d\'epend des observateurs ce qui a conduit Hawking \`a \'emettre son hypoth\`ese de rayonnement
des trous noirs (inobserv\'ee \`a ce jour). Or ce fait, s'il est \'etabli, est mal compris dans le cadre de la th\'eorie. 
La th\'eorie des cordes est beaucoup plus efficace sur ce point.

\appendix

\chapter{Rappels sur les transform\'ees de Fourier} \label{fourier}
Le but de ce paragraphe est de rappeler les d\'efinitions et propri\'et\'es de base 
qui concernent la transformation de Fourier et qui nous seront utiles lors de l'\'etude de la m\'ecanique quantique. 
Soit $f:\mR^3 \to \mC$. On suppose que $f$ tend suffisamment vite vers $0$ \`a l'infini.
 Nous ne serons pas plus pr\'ecis dans ce texte. 

\begin{definition}
La fonction sur $\mR^3$ d\'efinie par 
$$g(p) = \frac{1}{(2 \pi \hbar)^{3/2}}  \int_{\mR^3} f(x) e^{\frac{i}{\hbar}  p\cdot x} dx$$
o\`u $\cdot$ repr\'esente le produit scalaire canonique de $\mR^3$ est appel\'ee {\em transform\'ee de Fourier} de $f$.  
\end{definition}

\noindent Bien \'evidemment, nous avons choisi ici les constantes de normalisation adapt\'ees \`a la situation qui nous int\'eresse. Alors, la transformation de Fourier a plusieurs propri\'et\'es qui nous seront utiles: \\

\noindent {\bf Inversion de la transform\'ee de Fourier: } Si $g$ est la transform\'ee de Fourier de $f$, alors on a aussi: 
\begin{eqnarray} \label{inverse_tf}
f(x)=   \frac{1}{(2 \pi \hbar)^{3/2}} \int_{\mR^3} g(p) e^{-\frac{i}{\hbar}  p\cdot x} dp.
\end{eqnarray} 

\noindent {\bf Th\'eor\`eme de Plancherel-Parseval: }
La transformation de Fourier est une isom\'etrie de $L^2(\mR^3)$. Autrement dit, si $g$ est la transform\'ee de Fourier de 
$f$, on a 
\begin{eqnarray} \label{plancherel}
\int_{\mR^3} |f(x)|^2 dx = \int_{\mR^3} |g(p)|^2 dp.
\end{eqnarray}



\chapter{Rappels sur les espaces de Hilbert} \label{hilbert}
La formulation de la m\'ecanique quantique s'appuie sur les espaces de Hilbert. Nous faisons ici quelques rappels sur
le sujet. 

\section{D\'efinitions, g\'en\'eralit\'es}
\subsection{\'Enonc\'es} 

\begin{definition} 
\begin{enumerate}
\item Soit $E$ un espace vectoriel complexe. Une forme $(\cdot,\cdot):\mC \times \mC \to \mC$ est appel\'ee {\em produit scalaire hermitien} si 
\begin{itemize}
\item pour tous $x,y,z \in E$, $\al \in \mC$, on a 
$(x,y) = \overline{(y,x)}$, $(x,\al y)=\al (x,y)$, $(x,y+z)= (x,y) +(x,z)$;
\item pour tout $x \in E$, $\| x\|^2:=(x,x) $ est un r\'eel positif ou nul, et nul si et seulement si $x=0$.
\end{itemize} 
Un espace de Hilbert (hermitien) est un espace vectoriel complexe $H$ muni d'un produit scalaire hermitien
 $(\cdot,\cdot)$ pour lequel $H$ est complet. 
\end{enumerate}
\end{definition}

\noindent Soit maintenant un espace de Hilbert $H$. On fixe $y \in H$. La forme lin\'eaire $x \to (y,x)$ est 
continue. Un th\'eor\`eme important est le r\'esultat suivant qui affirme que la r\'eciproque est vraie. 

\begin{theorem} \label{riesz} ({\bf (de repr\'esentation de Riesz) }\\
Soit $H$ un espace de Hilbert et $f$ une forme lin\'eaire continue sur $H$. Alors, il existe un unique $y \in H$ tel que pour tout $x \in H$, 
$f(x)= (y,x)$. Autrement dit, on a une isomorphisme canonique de $H$ sur l'espace dual $H'$ des formes lin\'eaires continues de $H$.
\end{theorem}

\noindent En raison de son importance pour la m\'ecanique quantique, nous donnons une d\'emonstration de ce th\'eor\`eme
dans le paragraphe \ref{demriesz}. Si chacun est familier avec la dimension finie, il n'est sans doute pas inutile de faire quelques rappels sur 
ce qui se passe en dimension infinie. En particulier, rappelons le r\'esultat fondamental suivant 
\begin{theorem} \label{sommedirecte} 
 SOit $H$ un espace de Hilbert et $F$ un sous-espace ferm\'e. Alors, $F \oplus F^\perp= H$. 
\end{theorem}

\noindent Nous ne donnons pas la preuve de ce r\'esultat. Nous l'utiliserons plus loin. 
\subsection{D\'emonstration du th\'eor\`eme \ref{riesz}} \label{demriesz} 
Soit $f:H \to H$ une application lin\'eaire continue. S'il existe $y \in H$ tel que $f=(y,\cdot)$ alors il est unique.
En effet, si $z \in H$ v\'erifie la m\^eme relation, alors $(y-z,\cdot)$ est l'application lin\'eaire nulle et donc 
$\|y-z\|^2= (y-z,y-z)=0$ ce qui implique $y=z$. 
Il reste donc \`a montrer l'existence. Si $f=0$, on prend $y=0$. Sinon, $Ker(f)$ est un espace vectoriel ferm\'e (image r\'eciproque 
d'un ferm\'e par $f$) distinct de $H$. D'apr\`es le th\'eor\`eme \ref{sommedirecte}, $Ker(f)^\perp$ n'est donc pas r\'eduit
\`a $\{ 0\}$. Soit $b \in Ker(f)^\perp$. Posons $p_x= x - \frac{f(x) }{f(b)}b$. On a $f(p_x)=0$ et donc $b \perp p_x$. 
Cela implique que $(b,p_x)= 0$ i.e. $f(x) = \left( \frac{f(b)}{\| b \|^2} b , x  \right)$. On obtient le $y$ cherch\'e 
en posant $y = \frac{f(b)}{\| b \|^2} b$. 
  
\section{Spectre des op\'erateurs}  
\subsection{D\'efinitions} 
\noindent Nous avons vu l'importance en m\'ecanique quantique de pouvoir \'ecrire  un vecteur comme  combinaison 
lin\'eaire de vecteur propres. Il faut rappeler d'abord quelques d\'efinitions. 

\begin{definition}
 Soient $H$ un espace de Hilbert, $T: H \to H$ un op\'erateur (un application lin\'eaire de $H$ dans lui-m\^eme) 
et $\lambda \in \mC$.
\begin{itemize}
 \item $\lambda$ est une {\em valeur propre} de $T$ s'il existe $x \in H$, $x \not=0$ tel que $T(x)= \lambda x$ autrement 
dit si $T-\lambda Id$ n'est pas injectif.  L'ensemble des valeurs propres est not\'e $Vp(T)$. L'ensemble $E_\lambda:=
ker(T-\lambda Id)$ est l'{\em espace propre} associ\'e \`a $\lambda$. 
\item $\lambda$ est une {\em valeur spectrale} de $T$ si $T - \lambda Id$ n'est pas inversible. L'ensemble des valeurs
spectrales est appel\'e {\em spectre} de $T$ et est not\'e $\sigma(T)$.  
\end{itemize}
  \end{definition}

\noindent En dimension finie, il est bien \'evident que ces deux d\'efinitions sont \'equivalentes tandis qu'en dimension
infinie, une valeur propre est aussi une valeur spectrale mais la r\'eciproque est fausse en g\'en\'eral.
Il est important de noter que 
\begin{prop} \label{vpreelle}
 Les valeurs propres d'un op\'erateur auto-adjoint (i.e. un op\'erateur tel que pour tous $x,y \in H$, $(Tx,y)=(x,Ty)$) 
sont r\'eelles.  
\end{prop}
\noindent C'est une proposition tr\`es importante: en m\'ecanique quantique, les mesures d'une quantit\'e 
physique sont 
des valeurs propres d'op\'erateurs auto-adjoints. Cette proposition assure que ces grandeurs sont r\'eelles. 
La preuve est \'evidente: si $\alpha$ est une valeur propre et $x$ un vecteur propre associ\'e, on a, puisque $T$ est auto-adjoint,
$$\alpha \|x\|^2= (x,Tx)=(Tx,x)= \overline{\alpha} \| x \|^2,$$
ce qui implique $\alpha=\overline {\alpha}$ et donc $\alpha \in \mR$. 
En travaillant un peu plus, on peut voir que  la proposition est 
vraie aussi avec les valeurs spectrales mais ce sont les valeurs propres qui seront importantes pour 
la m\'ecanique quantique. 

\subsection{D\'ecomposition spectrale d'un op\'erateur compact} 
L'un des princi--pes importants de la m\'ecanique quantique  est de repr\'esenter les quantit\'es 
physiques par des 
op\'erateurs d'un espace de Hilbert dont les vecteurs propres forment une base hilbertienne de $H$. 
Tous les op\'erateurs
ne poss\`edent pas cette propri\'et\'e. Un contre-exemple est donn\'e dans le paragraphe \ref{contreexemple}. 
Par contre, on a: 

\begin{theorem} \label{decomp}
 Soit $H$ un espace de Hilbert hermitien et $T$ un op\'erateur auto-adjoint compact (i.e. l'image de toute partie born\'ee est relativement 
compacte). Alors $H=\overline{\oplus_{\al \in Vp(T)} E_\al}$. De plus cette somme est orthogonale. 
\end{theorem}
 
\noindent Rappelons que par d\'efinition la somme directe infinie $\oplus_{i} E_i$ des espaces 
vectoriels $E_i$ est l'ensemble
des combinaisons lin\'eaires {\bf finies} de vecteurs de $E_i$.
Nous donnons une d\'emonstration de ce th\'eor\`eme
dans le paragraphe \ref{demds}.  Ce th\'eor\`eme implique  
\begin{corollary} \label{decompositionspectrale} 
  Soit $H$ un espace de Hilbert hermitien et $T$ un op\'erateur auto-adjoint compact, alors $Vp(T)$ est au plus d\'enombrable
 et tout vecteur $\psi$ de $H$ s'\'ecrit 
$$\psi = \sum_{\al \in Vp(T)} P_\al(\psi)$$
o\`u $P_\al$ est le projecteur orthogonal sur l'espace $E_\al$.  
\end{corollary}

\noindent Nous d\'emontrons ce corollaire dans le paragraphe \ref{demdecompositionspectrale}. Pour finir, en prenant 
une base hilbertienne dans chaque espace propre nous obtenons aussi
\begin{corollary} \label{decompositionspectrale2} 
 Soit $H$ un espace de Hilbert hermitien et $T$ un op\'erateur auto-adjoint compact. 
Alors, il existe une base 
hilbertienne orthonorm\'ee 
de $H$ (pas forc\'ement d\'enombrable) compos\'ee uniquement de vecteurs propres de $T$. 
\end{corollary}

\begin{remark}
 Remarquons que pour $\al \not=0$, l'espace propre $E_\al$ est de dimension finie. 
En effet, sinon, consid\'erons une 
suite orthonorm\'ee $(e_n)$ de $E_\al$ qui forme donc un sous-ensemble born\'e de $E_\al$. Puisque l'op\'erateur $T$ est compact,
$T(e_n)$ admet une sous-suite convergente. Mais, cela ne peut pas \^etre vrai car pour $n \not=m$, 
$$\|T(e_n)-T(e_m)\|^2 = \al^2 \| e_n - e_m\|^2 = 2 \al^2$$
la derni\`ere \'egalit\'e r\'esultant de Pythagore.  
\end{remark}

\subsubsection{Contre-exemple} \label{contreexemple}
Nous allons construire un espace de Hilbert $H$, un op\'e--rateur $T$ continu de $H$, auto-adjoint qui n'admet aucun vecteur 
propre. Nous d\'efinissons  $H= l^2(\mZ,\mC)$ comme l'espace des suites index\'ees par $\mZ$, \`a valeurs complexes, de 
carr\'e sommable, 
munies du produit hermitien canonique. 
Notons $A$ l'op\'erateur d\'efini par 
$A((x_i)_i)=(x_{i+1})_i$. Son adjoint est clairement donn\'e par 
$A^*((x_i)_i)= (x_{i-1})_i$.  
Notons donc $T = A + A^*$. On a donc $T((x_i)_i)=1/2(x_{i-1}+x_{i+1})$. Par construction, $T$ est auto-adjoint, born\'e 
(de norme $1$) donc continu. Par ailleurs, soit $x$ un vecteur propre associ\'e \`a une valeur propre $\la$.
On a donc une suite r\'ecurrente d\'efinie par $x_{i+1}-2\la x_i+x_{i-1}=0$. Notons $\al_1, \al_2$ les racines complexes 
de l'\'equation $y^2-2\la y+1=0$. Il est alors classique que:
\begin{itemize}
 \item si $\al_1 \not= \al_2$, la suite $x$ est de la forme $x_i= C_1 \al_1^i +C_2 \al_2^i$ et n'est donc pas dans $H$ 
sauf si elle est nulle (car elle n'est pas de carr\'e sommable). 
\item si $\al_1 =\al_2 =\al$,  la suite $x$ est de la forme $x_i= (C_1 i + C_2) \al^i$ qui n'est pas non plus de carr\'e 
sommable. 
\end{itemize}

\subsubsection{D\'emonstration du th\'eor\`eme \ref{decomp}} \label{demds}
Nous commen\c{c}ons par d\'emontrer le lemme suivant: 
\begin{lemma} \label{lemmevalp}
 La norme $\|T\|$ (norme d'op\'erateur) est valeur propre de $T$.  
\end{lemma}

\noindent {\bf Preuve du lemme \ref{lemmevalp}} \\ 
Notons $\al =\| T \|$. Si $\al=0$, le r\'esultat est \'evident. Supposons $\al \not=0$.  Par d\'efinition, il existe une suite $(x_n)$ de vecteurs unitaires telle que $\lim_n \| T(x_n) \|=\al$. 
Notons $A = \al Id -T$. On a par d\'efinition de $(x_n)$ 
\begin{eqnarray} \label{anto0}
 \lim_n (A(x_n),x_n) = 0.
\end{eqnarray}
L'op\'erateur $A$ est auto-adjoint et positif. En particulier, la forme 
$(u,v) \to B(u,v):=(Au,v)$ est bilin\'eaire sym\'etrique. Par Cauchy-Schwarz:  
\begin{eqnarray} \label{bx_n}
   \| A(x_n) \|^2  = B(x_n,A(x_n)) 
\leq B(x_n,x_n)  B(A(x_n),A(x_n)).  
  \end{eqnarray} 
D'apr\`es \eref{anto0}, $\lim_{n} B(x_n,x_n)=0$. De plus, 
$$B(A(x_n),A(x_n))= (A^2(x_n),A(x_n)) \leq \|A \| \|A(x_n)\| \leq \|A \|^2.$$ 
On d\'eduit de \eref{bx_n} que $\lim_n A(x_n) = 0$ et donc que $\lim_n (T(x_n) -\al x_n) = 0$. Par ailleurs, 
$T$ \'etant compact, $T(x_n)$ converge (quitte \`a prendre une sous-suite). Ainsi $(\al x_n)_n$ est convergente. Puisque $\al \not= 0$,
$(x_n)$ converge vers un $x$ non nul puisqu'unitaire. Par continuit\'e de $T$, on a 
$$T(x) - \al x = \lim_n (T(x_n) -\al x_n) = 0.$$
Cela d\'emontre le lemme. $\square$ \\

\noindent On peut maintenant passer \`a la \\

\noindent {\bf Preuve du th\'eor\`eme \ref{decomp}}
Notons $F=\overline{\oplus_{\al \in Vp(T)} E_\al}$. Supposons $F \not=H$. D'apr\`es le th\'eor\`eme \ref{sommedirecte}, l'espace
vectoriel $F^\perp$ n'est pas r\'eduit \`a $\{ 0\}$. De plus, il est stable par $T$. On d\'eduit du lemme \ref{lemmevalp}
qu'il existe $\al$ valeur propre de $T_{/F^\perp}$. Soit $x$ un vecteur propre non nul associ\'e. 
Puisque $\al$ est aussi valeur propre de $T$, $x\in F \cap F^\perp$. On en d\'eduit une contradiction puisque d'apr\`es le th\'eor\`eme 
\ref{sommedirecte}, $F \cap F^\perp= \{ 0\}$.$\square$

\subsubsection{D\'emonstration du corollaire \ref{decompositionspectrale}} \label{demdecompositionspectrale}
Remarquons d'abord que $P_\al$ est bien d\'efinie. En effet,
$E_\al = (T-\al Id)^{-1}(\{ 0 \})$ est ferm\'e en tant qu'image 
r\'eciproque d'un ferm\'e par une application continue. 
D'apr\`es le th\'eor\`eme \ref{sommedirecte}, on a $H=E_\al \oplus E_\al^\perp$ ce qui donne naturellement 
$P_\al$. 
Ensuite, le th\'eor\`eme dit que tout vecteur $\psi$ est limite de $\psi_n$ o\`u chaque $\psi_n = \sum_{i=1}^{k_n} \psi_{\al_i}^n$ o\`u
les $\psi_{\al^n_i}^n \in E_{\al^n_i}$ sont des vecteurs propres de $T$ deux \`a deux orthogonaux. Notons 
$\Om_n = \{ \al_1^n, \cdots, \al_k^n \}$ les valeurs propres impliqu\'ees dans $\psi_n$ et posons $\Om_n'= \cup_{k=0}^n \Om_n$ les valeurs
propres impliqu\'ees dans tous les $\psi_i$ pour $i \leq n$.  
Posons maintenant $u_n= \sum_{\al \in \Om_n'} P_\al(\psi)$. 
Il est facile de voir que 
$\| \psi -u_n\| \leq  \| \psi -\psi_n\|$ si bien que $(u_n)$ converge vers $\psi$. Par ailleurs, par d\'efinition,
$$\lim_n u_n = \sum_{\Om'} P_\al(\psi)$$
o\`u $\Om'= \cup_{ n \in \mN} \Om_n'$ est un ensemble au plus d\'enombrable. On en d\'eduit le corollaire. 

\subsection{D\'ecomposition spectrale dans le cadre de la m\'ecanique quantique} \label{discussion} 
Nous avons donn\'e dans le paragraphe pr\'ec\'edent le th\'eor\`eme de d\'ecomposi--tion spectrale des op\'erateurs compacts. 
Malheureusement, ces op\'erateurs sont rares et notamment, ceux dont nous aurons besoin ne le sont pas en g\'en\'eral. 
Rappelons  le principe $IIc)$ de m\'ecanique quantique \'enonc\'e dans le chapitre \ref{mecanique_quantique}. \\

{\em Principe de d\'ecomposition spectrale:} Notons $(a_\al)_{\al}$ 
les valeurs propres de l'observable $\hat{A}$. Soit $\psi$  le vecteur d'\'etat du syst\`eme. 
La probabilit\'e de trouver $a_\al$ en mesurant la grandeur physique $a_\al$ \`a l'instant $t$ est \'egale \`a 
$\| \psi_\al \|^2$ o\`u $\psi_\al$ est la projection orthogonale de $\psi(t)$ sur le 
sous-espace propre de $\hat{A}$ associ\'e \`a la valeur propre $a_\al$. \\

\noindent Ce principe sous-entend que l'op\'erateur $\hat{A}$ a une propri\'et\'e de d\'ecomposition spectrale telle que 
celle donn\'ee pour les op\'erateurs compacts. Or, on a vu dans le paragraphe \ref{contreexemple} que, m\^eme pour un op\'erateur
born\'e, ce r\'esultat est faux en g\'en\'eral. Que dire alors des op\'erateurs simples comme ceux rencontr\'es
en m\'ecanique ondulatoire (position et quantit\'e de mouvement) qui ne sont m\^eme pas born\'es ? 
Pour faire les choses rigoureusement, il convient de se placer dans un espace plus grand que $\ep_H$. Nous ne 
donnerons pas cette construction ici mais nous expliquons comment cela se passe lorsque 
$\hat{A}$ est l'observable position ou quantit\'e de mouvement de la m\'ecanique ondulatoire. En particulier,
on travaille dans l'espace de Hilbert $L^2$. 
Commen\c{c}ons par rappeler que si $x$ est le vecteur position de la particule, l'observable associ\'ee 
est la multiplication par $x$. Pour $x$ fix\'e, l'espace propre est de dimension $1$ et le projecteur associ\'e est 
la distribution de Dirac 
$\delta_{x}$ au point $x$.   
Cette fonction est nulle dans $L^2$ et donc $L^2$ n'est pas le bon espace de fonctions. N\'eanmoins 
chaque fonction de $L^2$ peut s'\'ecrire comme une somme infinie 
(en l'occurrence une int\'egrale)  de projections sur les espaces propres. 
Plus pr\'ecis\'ement, $Vp(\hat{x})= \mR^3$ et on a $\psi(y) = \int_{x \in \mR^3} \delta_y(x) \psi(x) dx$. \\

\noindent De m\^eme si on regarde la premi\`ere composante du vecteur quantit\'e de mouvement
$p_x$, l'op\'erateur associ\'e est $\hat{p}_x = -i \hbar \partial_x$ dont les vecteurs propres sont les fonctions 
$p \to e^{\frac{i}{\hbar} p x}$ associ\'ees \`a la valeur propre  $x \in \mR$. On peut d\'ecomposer $\psi$ dans 
cette ''base`` de vecteurs propres en utilisant l'inversion de la transform\'ee de Fourier. La projection 
de $\psi$ sur l'espace propre associ\'e \`a $x$ est $p \to e^{i\frac{i}{\hbar} p x} \phi(p)$ o\`u $\phi$ est la transform\'ee 
de Fourier de $\psi$. En int\'egrant sur l'ensemble des valeurs propres ($\mR$), on retrouve bien $\psi$. \\
 
\noindent Cette mani\`ere de proc\'eder se g\'en\'eralise de mani\`ere rigoureuse \`a tout op\'erateur auto-adjoint 
m\^eme non born\'e. Nous n'expliquerons pas ici la mani\`ere de s'y prendre.

\section{Diagonalisation simultan\'ee de deux op\'erateurs qui commutent} 
Consid\'erons deux op\'erateurs $A$ et $B$ auto-adjoints d'un espace de Hilbert $H$ satisfaisant le principe de d\'ecomposition
spectrale discut\'ee dans le paragraphe pr\'ec\'edent.
Alors 
\begin{prop} \label{diagsimul}
 Si $AB=BA$, il existe une base hilbertienne compos\'ee de vecteurs qui sont propres \`a la fois pour $A$ et $B$.  
\end{prop}
En effet, soit $x$ un vecteur propre de $A$ associ\'e \`a une valeur propre $\la$ alors $ABx= BAx = \la Bx$ et 
donc $Bx$ est aussi vecteur propre de $A$ associ\'e \`a $\la$. En particulier, les espaces propres de $A$ 
sont stables par $B$. 
Il ne reste plus qu'\`a appliquer le principe de d\'ecomposition spectrale 
\`a l'op\'erateur $B$ restreint \`a chacun des sous-espaces propres de $A$. 

\section{Produit tensoriel d'espaces de Hilbert} \label{tensoriel}
Soient $G,H$ deux espaces de Hilbert (cela pourrait \^etre seulement des espaces vectoriels). 
On peut montrer (ce que nous ne ferons pas ici) qu'\`a isomorphisme pr\`es, 
il existe un unique espace de Hilbert not\'e $G \otimes H$  appel\'e {\em produit tensoriel de $G$ par $H$} et 
une  application lin\'eaire $T: G \times H \to G \otimes H$ tels que 
\begin{itemize} 
 \item $T(G \times H)$ engendre $G \otimes H$;
\item Si $(g_i)_i$, $(h_j)_j$ sont des bases hilbertiennes de $G$ et $H$ alors $(T(g_i,h_j))_{i,j}$ est une base hilbertienne de 
$G \otimes H$. 
\end{itemize}
On notera g\'en\'eralement $u \otimes v$ \`a la place de $T(u,v)$. Si $G$ et $H$ sont de dimension finie, c'est aussi le 
cas de $G \otimes H$ et on a 
$$dim(G\otimes H )= dim(G) dim(H).$$
Le produit hermitien sur $G\otimes H$ est tel que pour $u,u' \in G$, $v,v' \in H$, 
$$(u\otimes v, u'\otimes v') = (u,u')_G(v,v')_H.$$
Pour finir, si $f,f'$ sont des op\'erateurs de $G$ et $H$, on obtient un op\'erateur $f \otimes f'$ de $G \otimes H$ 
en posant 
$$f\otimes f'(u\otimes v)= f(u)\otimes f(v).$$

\chapter{Les spineurs} \label{spineurs}
La construction pr\'ecise des spineurs en g\'en\'eral est compliqu\'ee. Une tr\`es bonne r\'ef\'erence est par exemple 
\cite{hijazi:99}. Nous nous pla\c{c}ons ici dans le cadre simple de la relativit\'e restreinte, c'est-\`a-dire sur $\mR^4$ munie
de la forme quadratique 
$\eta:= dx^2 - dt^2$. 
L'{\em alg\`ebre de Clifford} $Cl(\mR^4,\eta)$ est d\'efinie par 
$$Cl(\mR^4, \eta) := T(\mR^4) / I(\mR^4,\eta)$$ 
o\`u $T(\mR^4)$ est l'alg\`ebre tensorielle de $\mR^4$, autrement dit 
$$T(\mR^4)= \oplus_{k=0}^\infty (\otimes_{i=1}^k \mR^4)$$
et o\`u $I(\mR^4, \eta)$ est l'id\'eal engendr\'e par les \'el\'ements de la forme $x \otimes x +\eta(x,x) 1$. 
La loi de multiplication interne sur $Cl(\mR^4)$ (c'est-\`a-dire celle induite par $\otimes$ sur le quotient) est 
not\'ee ''$\cdot$`` et s'appelle {\em multiplication de Clifford}. On peut alors v\'erifier 
que $Cl(\mR^4,\eta)$ est l'alg\`ebre engendr\'ee par $\mR^4$ avec la relation
\begin{eqnarray} \label{defclif} 
 v \cdot w+ w\cdot v= -2\eta(v,w) 1
\end{eqnarray}
pour tous $v,w \in \mR^4$. En particulier, notons $(e_1,e_2,e_3,e_4)$ la base canonique de $\mR^4$, on a les relations  
$$e_i \cdot e_j+e_j \cdot e_i= 0 $$
si $i \not=j$ et 
 \[ e_i^2 = \left| \begin{array}{ccc}
                    -1 & \hbox{ si } & i \in \{1,2,3\}\\
1 & \hbox{ si } & i =4. 
                   \end{array} \right. \]

On peut v\'erifier que 
$$\{ e_{i_1} \cdot e_{i_2} \cdot ...\cdot e_{i_k} | 1 \leq i_1 < \cdots < i_k, 0 \leq k \leq 4 \}$$
est une  base de l'espace vectoriel $Cl(\mR^4,\eta)$. \\

\noindent En fait, ce sera la complexification de l'alg\`ebre de Clifford
$$\mC l(\mR^4,\eta) = Cl(\mR^4,\eta) \otimes_{\mR} \mC$$
qui nous int\'eressera. On peut maintenant regarder les \'el\'ements de  $\mC l(\mR^4,\eta)$ 
comme des endormophismes de l'espace vectoriel
$\mC l(\mR^4,\eta)$. En effet, $x \in \mC l(\mR^4,\eta)$ s'identifie \`a $y \to x \cdot y$. On peut montrer que sous cette action, il existe un espace vectoriel not\'e $\Si_4$ qui est stable
et de dimension minimale $4$. Si l'on avait travaill\'e dans $\mR^n$, la dimension aurait \'et\'e de $2^{[n/2]}$. Sans 
rentrer dans les d\'etails et notamment en oubliant les difficult\'ees li\'ees \`a la parit\'e de la dimension, cette 
dimension de $2^{[n/2]}$ n'est pas surprenante: l'espace des applications lin\'eaires d'un espace vectoriel de dimension 
$m=2^{n/2}$ est $2^m= 2^n$ qui est exactement en dimension paire du moins, la dimension complexe de  $\mC l(\mR^4,\eta)$.   
L'expression explicite de $\Si_4$ est compliqu\'ee
et inutile. Il suffit de retenir, en r\'e\'ecrivant les choses de mani\`ere plus habituelle, 
qu'il existe une repr\'esentation irr\'eductible
$$\rho: \mC l(\mR^4,\eta) \to \hbox{End}_{\mC } (\Si_4)$$
appel\'ee {\em repr\'esentation spinorielle complexe}. Pour $v \in \mR^4 \subset  \mC l(\mR^4,\eta) $ et $\psi \in \Si_4$, on notera
$v \cdot \psi$ plut\^ot que $\rho(v) \psi$. 
Soit maintenant $\psi: \mR^4 \to \Si_4$. En prenant une base de $\Si_4$, $\psi$ est donc une fonction \`a 
$4$ composantes complexes: elle correspond au vecteur d'\'etat qui appara\^{\i}t dans l'\'equation de Dirac dans le 
chapitre sur la m\'ecanique quantique relativiste. 
On d\'efinit maintenant l'op\'erateur 
$$D:= \sum_{i=1}^4 e_i \cdot \partial_i.$$ 
L'op\'erateur $D$ est appel\'e {\em op\'erateur de Dirac}. Les relations \eref{defclif} montrent imm\'ediatement que 
\begin{eqnarray} \label{forsl}
 D^2 =  -\square= -\Delta_{x} +\partial_{tt}. 
\end{eqnarray}
On peut aussi construire sur $\Si_4$ de mani\`ere naturelle un produit hermitien qui est tel que $D$ est auto-adjoint.

\providecommand{\bysame}{\leavevmode\hbox to3em{\hrulefill}\thinspace}

\end{document}